\documentclass{tmsce}

\usepackage{hyperref}
\hypersetup{colorlinks=true,citecolor=blue}

\usepackage{mathbbol}
\usepackage{amsmath, amsfonts, amsthm, mathtools}
\usepackage{array}
\usepackage{xcolor}
\usepackage{bm, bbm}
\usepackage{amssymb}
\usepackage{enumitem}
\usepackage{url}
\usepackage{graphicx}
\usepackage[round]{natbib}

\usepackage{cases}
\usepackage{subcaption}
\usepackage{float}
\usepackage{thmtools}

\usepackage{tikz}
\usetikzlibrary{decorations.pathreplacing}

\title{
    Convergence rates for truncated and discretised\\ signature expansions
}

\authors{Eduardo Abi Jaber$^{a}$, Dorinel Bastide$^{b}$, Ronan Pécheul$^{a,b}$, Clément Rey$^{a}$}
\affiliation{
    $^{a}$ CMAP, École Polytechnique, France\\ 
    $^{b}$ BNP Paribas RISK\footnote{This article represents the opinions of the authors, and it is not meant to represent the position or opinions of BNP Paribas or its members.}, France.
}

\newtheorem{theorem}{Theorem}[section]
\newtheorem{corollary}[theorem]{Corollary}
\newtheorem{lemma}[theorem]{Lemma}
\newtheorem{proposition}[theorem]{Proposition}
\newtheorem{remark}[theorem]{Remark}
\newtheorem{example}[theorem]{Example}
\newtheorem{definition}[theorem]{Definition}

\numberwithin{equation}{section}

\DeclareSymbolFontAlphabet{\amsmathbb}{AMSb}

\newcommand{\R}{\amsmathbb{R}}
\newcommand{\N}{\amsmathbb{N}}
\newcommand{\Nbar}{\overline{\N}}

\newcommand{\E}{\amsmathbb{E}}
\newcommand{\W}{\amsmathbb{W}}
\newcommand{\X}{\amsmathbb{X}}
\newcommand{\Q}{\amsmathbb{Q}}
\newcommand{\F}{\amsmathbb{F}}
\newcommand{\Prob}{\amsmathbb{P}}
\newcommand{\scriptG}{\mathcal{G}}
\newcommand{\VaR}{\mathrm{VaR}}
\newcommand{\ES}{\mathrm{ES}}

\newcommand{\tnorm}[1]{{\left\vert\kern-0.25ex\left\vert\kern-0.25ex\left\vert #1 \right\vert\kern-0.25ex\right\vert\kern-0.25ex\right\vert}}

\newcommand{\absE}[1]{\left| \E \left[ #1 \right] \right|}
\newcommand{\textabsE}[1]{\big| \E \big[ #1 \big] \big|}

\newcommand{\norm}[1]{\left\lVert #1 \right\rVert}
\newcommand{\textnorm}[1]{\big\lVert #1 \big\rVert}

\newcommand{\brackets}[2]{\left\langle #1, #2 \right\rangle}
\newcommand{\textbrackets}[2]{\langle #1, #2 \rangle}

\newcommand{\intset}[2]{\{#1, \dots, #2\}}

\newcommand{\Wh}{\widehat{W}}
\newcommand{\WH}{\widehat{\W}}
\newcommand{\Wt}{\widetilde{W}}
\newcommand{\WT}{\widetilde{\W}}
\newcommand{\nup}{\overline{n}}
\newcommand{\ndown}{\underline{n}}

\newcommand{\diff}{\mathrm{d}}

\newcommand{\indic}[1]{\mathbbm{1}_{#1}}

\newcommand{\kt}{\tilde{k}}
\newcommand{\pih}{\widehat{\pi}}
\newcommand{\pit}{\widetilde{\pi}}

\usepackage{shuffle}
\DeclareFontFamily{U}{bigshuffle}{}
\DeclareFontShape{U}{bigshuffle}{m}{n}{
  <5-8> s*[1.7] shuffle7
  <8->  s*[1.7] shuffle10
}{}
\DeclareSymbolFont{BigShuffle}{U}{bigshuffle}{m}{n}
\DeclareMathSymbol\bigshuffle{\mathop}{BigShuffle}{"001}

\definecolor{lightblue}{rgb}{0.25, 0.41, 0.88}
\newcommand{\lbold}{\boldsymbol{\ell}}
\newcommand{\phibold}{\boldsymbol{\varphi}}
\newcommand{\blold}[1]{\textcolor{lightblue}{\bm{\mathrm{#1}}}}
\newcommand{\empt}{\textcolor{lightblue}{\scriptstyle\bm{\varnothing}}}
\newcommand{\blue}[1]{\textcolor{lightblue}{#1}}

\newcommand{\udown}{\blold{\underline{u}}}
\newcommand{\uup}{\blold{\overline{u}}}
\newcommand{\updown}{\blold{\underline{\upsilon}}}
\newcommand{\upup}{\blold{\overline{\upsilon}}}

\begin{document}
\maketitle

\begin{abstract}
Linear expansions in Brownian Stratonovich signatures, which originate in stochastic Taylor expansions, are increasingly used to represent and model  processes with memory. Their numerical use requires truncating the signature at a level $M$ and discretising the driving path along $N$ steps. Under a simple  growth condition on the coefficients of the linear expansion, we establish strong and weak convergence rates for both approximations, uniformly in time. The strong truncation error decays like $(B^M/(M!)^{\varepsilon})^{a}$ for some constant $B$, uniformly in $N$, where $\varepsilon$ measures the growth of the coefficients and $a$ the regularity of the test function; the weak error does not improve on it. The  discretisation error is of order $N^{-a/2}$ in $L^p$, and $N^{-1}$ weakly for polynomial test functions, uniformly in $M$. All these rates are sharp for polynomial test functions, and $M$, $N$ can be chosen independently. We extend the estimates to errors in probability, on distribution functions and on risk measures, and illustrate them on stochastic delay equations, Volterra processes and option pricing in signature volatility models.
\vspace{-1em}
\paragraph{keywords.} signatures, weak and strong convergence rates, risk measures, non-Markovian processes, signature volatility models.
\vspace{-1em}
\end{abstract}

\tmsceendfrontmatter



\section{Introduction} \label{sec:introduction}

The signature of a path is the collection of its iterated integrals. Introduced by \citet*{ChenKuoTsai1954}, these integrals are central to stochastic Taylor expansions~\citep*{KloedenPlaten1992,BenArous1989} and to  rough path theory~\citep*{Lyons1998}. Signature coordinates play the role of monomials on path space: finite linear combinations are analogues of polynomials, while convergent infinite expansions are analogues of power series. Through the shuffle identity and suitable point-separation properties, they yield Stone--Weierstrass approximation theorems on compact sets~\citep*{LevinLyonsNi2013}. These approximation properties extend to broader path spaces~\citep*{CuchieroPrimaveraSvaluto2025}, weighted and $L^p$ settings~\citep*{CuchieroSchmockerTeichmann2023,CeylanPromel2025,CeylanPromel2026}, and progressively measurable controls~\citep*{BankBayerHagerRiedelNauen2025}.

Beyond approximation, suitable classes of Markovian and non-Markovian processes, including solutions to stochastic differential equations, convolution Volterra equations and delay equations, admit exact, possibly local, representations as infinite linear combinations of signature elements of the time-extended driving Brownian motion~\citep*{BenArous1989,AbiJaberGerardHuang2024}. Infinite Brownian signature expansions also arise in conditional expectations of suitable analytic functionals and in non-Markovian transform formulas~\citep*{CuchieroSvalutoFerroTeichmann2026,AbiJaberAttalSotnikov2026Affine}, as well as in stochastic control with path-dependent objectives~\citep*{AbiJaberAttalSotnikov2026Control}.

These properties motivate the use of signatures in stochastic modelling and statistical learning. In mathematical finance, finite linear combinations of signature coordinates underpin non-parametric pricing and hedging of exotic derivatives~\citep*{LyonsNejadPerezArribas2020}. They also lead to signature volatility models, in which volatility is a linear functional of the signature of time-extended Brownian motion. Such models can be calibrated to option prices~\citep*{CuchieroGazzaniSvalutoFerro2023}, priced and hedged with Fourier methods~\citep*{AbiJaberGerard2025}, and their sensitivities analysed using  Malliavin calculus~\citep*{AbiJaberReySotnikov2026}. Infinite signature expansions also yield tractable pricing and hedging formulas for path-dependent options with market impact~\citep*{AbiJaberHainautMotte2025}.

In practice, every numerical use of such a representation requires two approximations. The signature is an infinite object and must be truncated at some level $M$. The signature of a Brownian path cannot be simulated exactly, so the path is replaced by its piecewise linear interpolation along $N$ steps, whose signature is computed exactly through Chen's identity. This is precisely what standard signature libraries do~\citep*{ReizensteinGraham2020, KidgerLyons2021}: $M$ is the truncation level passed to the routine, and $N+1$ the number of points of the stream. The purpose of this paper is to quantify both errors, explicitly in $M$, in $N$ and in the regularity of the functional of interest.

\subsection{Problem formulation} \label{sec:problem-formulation}

Fix a horizon $T > 0$, let $W$ be a $d$-dimensional Brownian motion, and let
$\WH$ be the Stratonovich signature of the time-extended path $\Wh = (t, W_t)_{t\in[0,T]}$. Its coordinates are indexed by the words $\blold{v}$ over the alphabet $\{\blold{0}, \blold{1}, \dots, \blold{d}\}$, the letter $\blold{0}$ standing for integration in time and $\blold{1}, \dots, \blold{d}$ for integration with respect to Brownian motion. We consider processes of the form
\begin{equation}\label{eq:signature-representation}
    X_t = \textbrackets{\lbold_t}{\WH_t} = \textstyle{\sum_{\blold{v}}} \lbold^{\blold{v}}_t\, \WH^{\blold{v}}_t, \qquad t \in [0,T],
\end{equation}
for a deterministic coefficient family $\lbold = (\lbold_t)_{t\in[0,T]}$, the sum being a genuinely infinite series. Truncating at level $M$ discards the words longer than $M$, and replacing $\WH$ by the signature $\WT(N)$ of the interpolated path discretises it. Writing $X_t(M,N) \coloneqq \textbrackets{\lbold_t}{\WT^{\le M}_t(N)}$, with $M = \infty$ meaning no truncation and $N = \infty$ no discretisation, the target is $X_t = X_t(\infty,\infty)$ and the computable approximation is $X_t(M,N)$. For a test function $g$, we control the strong error $\norm{g(X_t) - g(X_t(M,N))}_p$, where $\norm{\,\cdot\,}_p \coloneqq \E[|\cdot|^p]^{1/p}$, and the weak error $\absE{g(X_t) - g(X_t(M,N))}$. Three questions organise the analysis.

\begin{enumerate}[itemsep=0em]
    \item Which coefficient families admit such a control? Since no differential equation generates them, the assumption has to bear on their growth alone.
    \item At which rates in $M$ and $N$ do the errors vanish, how do these rates depend on the horizon and on the regularity of $g$, and can $M$ and $N$ be chosen independently?
    \item Does the weak error improve on the strong one?
\end{enumerate}

\subsection{The classical benchmark} \label{sec:classical-benchmark}

When $X$ solves a stochastic differential equation $\diff X_t = a(X_t)\,\diff t + b(X_t)\,\diff W_t$ with smooth coefficients, the Itō--Taylor and Stratonovich--Taylor expansions~\citep*{PlatenWagner1982, KloedenPlaten1992} answer these questions. The coefficients are generated by iterating the differential operators associated with $a$ and $b$, and the remainder is controlled through the same structure. Truncating along the degree, in which time integration counts twice, yields a strong error in $L^2$ of order $(cT)^{(M+1)/2} / \sqrt{\lfloor (M+1)/2 \rfloor !}$~\citep*[Proposition~5.10.1]{KloedenPlaten1992} for some constant $c > 0$.
The present setting departs from this theory in three ways. First, we truncate along the word length $|\blold{v}|$, and we still obtain the factor $T^{(M+1)/2}$, together with a factorial decay $(M!)^{-\varepsilon}$ (Theorem~\ref{th:Lp-truncation-bound}). Here $\varepsilon \in (0, 1/2]$ measures how much of the factorial decay of the signature moments is left once the growth of the coefficients is accounted for. Second, and more importantly, the coefficients in~\eqref{eq:signature-representation} are \emph{given} rather than generated by an equation, so every estimate, including the factorial gain, must be derived from their size alone. The class~\eqref{eq:signature-representation} thus goes beyond the setting of \citet*{KloedenPlaten1992}: it contains non-Markovian processes, such as delay equations and Volterra processes, as shown by \citet*{AbiJaberGerardHuang2024}. The Ornstein--Uhlenbeck process belongs to both classes~\citep*[Theorem~4.2]{AbiJaberGerardHuang2024}. Third, the discretisation differs. \citet*[Section~5.8]{KloedenPlaten1992} approximate the multiple stochastic integrals through a truncated Fourier expansion of the Brownian bridge, which introduces a further approximation parameter. We instead replace the driving path by its piecewise linear interpolation on the grid, the classical approximation of Stratonovich integrals going back to \citet*{WongZakai1965}. Its signature is computed exactly through Chen's identity and remains group-like, so that the shuffle identity holds at every discretisation level, and it only depends on the grid increments of the Brownian motion. Within this class of approximations, the rate $N^{-1/2}$ cannot be improved (Remark~\ref{rem:optimality-statement}), and the discretisation does not interact with the truncation level, which yields estimates uniform in $M$.

\subsection{Contributions} \label{sec:contributions}

We introduce in Section~\ref{sec:main-results} the coefficient space $\mathcal{L}^1$ of tensors satisfying $|\lbold^{\blold{v}}| \lesssim (\norm{\blold{v}}!)^{1/2-\varepsilon}$ for some $\varepsilon \in (0,1/2]$. This condition is calibrated against the moments of the signature components, which decay like $(\norm{\blold{v}}!)^{-1/2}$ for both the signature and its discretisation (Proposition~\ref{prop:signature-moment-bound}). The space is stable under the usual algebraic operations (Proposition~\ref{prop:L1-stability}), contains the class $\mathcal{A}_{\exp}$ of \citet*{AbiJaberGerardHuang2024} (Remark~\ref{rem:Aexp}), and, on it, the series \eqref{eq:signature-representation} and its discretised counterpart converge absolutely, in every $L^p$, uniformly in time and in $N$ (Corollary~\ref{cor:well-posed}). Our main results are the following; throughout, $a \in (0,1]$ is the regularity exponent of $g$, equal to $1$ for polynomials and to the Hölder exponent otherwise. We write $\Nbar \coloneqq \N \cup \{\infty\}$ and $\Nbar^* := \Nbar \setminus \{0\}$.

\begin{itemize}
    \item \textbf{Truncation error.} For every polynomial or Hölder continuous $g$,
    \[
        \norm{g(X_t(\infty,N)) - g(X_t(M,N))}_p \lesssim \Big( \frac{B^M}{(M!)^{\varepsilon}} \Big)^a, \qquad t \in [0,T],\ M \in \N,\ N \in \Nbar^*
    \]
    (Theorem~\ref{th:Lp-truncation-bound}): the error decays like $\exp\{-\varepsilon a M\log M + O(M)\}$, faster than any geometric sequence, uniformly in time \emph{and in the discretisation level}. The proof is a direct moment computation with explicit constants (Lemma~\ref{lem:truncation}). The exponent $\varepsilon$ is sharp, even for the weak error and uniformly in $N$ (Proposition~\ref{prop:truncation-sharpness}): in contrast with the Itō setting of \citet*{KloedenPlaten1992}, in general, no weak gain is possible for Stratonovich signatures.
    \item \textbf{Discretisation error.} For every word, $\textnorm{\WH^{\blold{v}}_t - \WT^{\blold{v}}_t(N)}_p \lesssim N^{-1/2}$ (Proposition~\ref{prop:Lp-element-discretisation}), whence a rate $N^{-a/2}$ for $\norm{g(X_t(M,\infty)) - g(X_t(M,N))}_p$, uniformly in time \emph{and in the truncation level} $M \in \Nbar$ (Theorem~\ref{th:Lp-discretisation-bound}). The proof combines Chen's identity with Rosenthal's inequality to isolate the contributions of independent Brownian increments. The weak error improves to $N^{-1}$ for polynomial test functions (Proposition~\ref{prop:weak-element-discretisation} and Theorem~\ref{th:weak-discretisation-bound}), and, in dimension one, the strong rate improves to $N^{-1}$ at the grid points. All these rates are optimal (Remark~\ref{rem:optimality-statement}).
    \item \textbf{Total error.} Since each estimate is uniform in the parameter left fixed, they combine into $\norm{g(X_t) - g(X_t(M,N))}_p \lesssim (B^M/(M!)^\varepsilon)^a + N^{-a/2}$, uniformly in $t$, $M$ and $N$ (Corollary~\ref{cor:total-error}). The two parameters can therefore be chosen independently: $M$ is set by the target accuracy, and $N$ of order $(M!)^{2\varepsilon}B^{-2M}$ balances the two terms.
    \item \textbf{Extensions and applications.} The estimates transfer to errors in probability, on distribution functions, and on risk measures, namely coherent risk measures including Expected Shortfall, and Value-at-Risk; the density conditions needed for distribution functions and Value-at-Risk bear on the target only, never on its approximations (Section~\ref{sec:extensions-statements}). We verify our assumptions for linear delay equations and Volterra processes. The shifted Riemann--Liouville representation of fractional Brownian motion lies on the boundary $\varepsilon = 0$ of $\mathcal{L}^1$; a direct Gaussian argument shows that its truncation error still decays geometrically on the whole horizon (Proposition~\ref{prop:RL}), whereas another boundary family diverges beyond a finite horizon (Example~\ref{ex:boundary-case}). Finally, in a signature volatility model the price is not known to be a linear functional of the signature, but its logarithm is, which yields a numerical scheme for option prices up to an arbitrarily small loss in $\varepsilon$ (Sections~\ref{sec:pricing} and~\ref{sec:pricing-proof}).
\end{itemize}

\subsection{Related literature} \label{sec:related-literature}

Replacing a Brownian path by its piecewise linear interpolation goes back to \citet*{WongZakai1965}, and the convergence of the resulting iterated integrals, with rates, has since been studied extensively within rough path theory; see~\citet*{FrizVictoir2010}. Closest to our discretisation analysis, \citet*{FrizRiedel2014} establishes sharp almost sure rates for approximations of Gaussian rough paths, of order essentially $N^{-1/2}$ in the Brownian case. That result is pathwise and concerns a fixed level of the signature, whereas we work in $L^p$ with an infinite linear functional, so that truncation and discretisation must be handled together. On the weak side, \citet*{NiXu2015} computes the exact rate $N^{-1}$ for expected Brownian signatures; our Proposition~\ref{prop:weak-element-discretisation} gives the same order as an upper bound for all words, including those involving time. From a statistical perspective, \citet*{LucchesePakkanenVeraart2025} studies the convergence of the discretised expected signature as the mesh vanishes; our bounds provide explicit in-fill rates in the Brownian case.

On the qualitative side, \citet*{CeylanPromel2025} proves global $L^p$-universal approximation theorems for linear functionals of the time-extended signature under an exponential moment condition on the driving path. In particular, the solution of an Itō stochastic differential equation can be approximated by a truncated linear functional of the time-extended Stratonovich Brownian signature~\citep*[Proposition~4.6]{CeylanPromel2025}. That statement is existential: the truncation level carries no rate, the coefficients depend on the target accuracy, and the exact signature is used. We address the complementary quantitative question for a coefficient family given in advance, and add the discretisation error. The two approaches are dual in their hypotheses: \citet*{CeylanPromel2025} constrains the law of the driving path, whereas $\mathcal{L}^1$ constrains the coefficients. Similarly, \citet*{CeylanPromel2026} shows that signatures of piecewise linear paths are rich enough to approximate the functionals of interest, while we quantify how fast a given representation is recovered as $M$ and $N$ grow. Finally, Itō signatures enjoy a martingale structure, exploited by the weak expansions of~\citet*{KloedenPlaten1992}, and have recently been used as tradable hedging instruments~\citep*{GuoWangZhang2026}. Stratonovich signatures satisfy the shuffle identity (Proposition~\ref{prop:shuffle-identity}), on which the algebra of this paper relies, and they are the ones appearing in the representations of~\citet*{AbiJaberGerardHuang2024}. A universal approximation theorem for Itō signatures is established in~\citet*{CeylanKwossekPromel2026} by extending the path with its bracket terms; in the Brownian setting this extension is redundant, and both signatures perform comparably~\citep*[Remarks~4.6 and~4.11]{CeylanKwossekPromel2026}. The distinction matters for payoffs that depend explicitly on quadratic variation, which we do not consider, and for the weak truncation rate, as discussed above.

\subsection{Outline and notation} \label{sec:outline}

Section~\ref{sec:framework} recalls the algebraic framework, the discretisation scheme and the moment estimates. Section~\ref{sec:main-results} introduces $\mathcal{L}^1$ and states the main results, which Section~\ref{sec:example} extends to errors in probability and on risk measures, and illustrates on linear delay equations, Gaussian Volterra processes and a signature volatility model. The truncation and discretisation estimates are proved in Section~\ref{sec:theorem}, the extensions to probability errors and risk measures in Section~\ref{sec:extensions}, and the pricing application in Section~\ref{sec:pricing-proof}. The remaining proofs are collected in the
appendices.

\paragraph{Notation.} For $a, b \ge 0$, $a \lesssim b$ means that $a \le Kb$ for a constant $K > 0$ that may depend on the parameters held fixed, but not on the variables with respect to which the estimate is stated. For instance, ``$\norm{Z_t(M,N)}_p \lesssim N^{-1/2}$, $t \in [0,T]$, $M \in \N$, $N \in \N^*$'' means that for all fixed $T$ and $p$ there exists $K$, independent of $t$, $M$ and $N$, such that $\norm{Z_t(M,N)}_p \le KN^{-1/2}$.
\section{Framework} \label{sec:framework}

This section fixes the objects of the representation~\eqref{eq:signature-representation}: the words that index the coefficients, the signature $\WH$ of the time-extended Brownian motion, its discretised counterpart $\WT(N)$, and the moment estimates against which the coefficient space of Section~\ref{sec:main-results} is calibrated.

\subsection{Words and tensor series} \label{sec:words}

Fix $d \in \N^*$ and the \emph{alphabet} $\{\blold{0}, \blold{1}, \dots, \blold{d}\}$. Throughout, the letter $\blold{0}$ refers to time and the letters $\blold{1}, \dots, \blold{d}$ to the components of a $d$-dimensional Brownian motion. A \emph{word} is a finite sequence of letters $\blold{v} = \blold{i_1 \cdots i_m}$, and $\empt$ denotes the empty word; we write $V$ for the set of words and $V_m$ for those of length $m$. For $\blold{v} = \blold{i_1 \cdots i_m}$:
\begin{itemize}
    \item $|\blold{v}| \coloneqq m$ is its \emph{length}, and $n_i(\blold{v})$ the number of occurrences of the letter $\blold{i}$, $i \in \intset{0}{d}$;
    \item $s(\blold{v}) \coloneqq (i_1, \dots, i_m)$ is the sequence of its letters;
    \item $\norm{\blold{v}} \coloneqq |\blold{v}| + n_0(\blold{v})$ is its \emph{degree}, in which $\blold{0}$ counts twice. Since time scales like the square of Brownian motion, the degree measures the size of a word in the Brownian scaling.
\end{itemize}
The concatenation of $\blold{v}$ and $\blold{u} = \blold{j_1 \cdots j_n}$ is
$\blold{vu} := \blold{i_1 \cdots i_m \, j_1 \cdots j_n}$.

\medskip

A \emph{tensor series} is a formal linear combination of words,
\[
    \lbold = \sum_{\blold{v} \in V} \lbold^{\blold{v}}\, \blold{v}, \qquad \lbold^{\blold{v}} \in \R,
\]
and we write $\blold{v} \in \lbold$ when $\lbold^{\blold{v}} \ne 0$. The space of tensor series is denoted $T((\R^{d+1}))$, and the subspace of those with finitely many non-zero coefficients $T(\R^{d+1})$.\footnote{Identifying the word $\blold{i_1 \cdots i_m}$ with $e_{i_1} \otimes\cdots\otimes e_{i_m}$, where $(e_0,\dots,e_d)$ is the canonical basis of $\R^{d+1}$, $T((\R^{d+1}))$ is the extended tensor algebra $\bigtimes_{m \ge 0}(\R^{d+1})^{\otimes m}$.} We use the following operations.
\begin{itemize}
    \item \emph{Truncation and projection.} For $M \in \N$, $\lbold^{\le M} \coloneqq \sum_{|\blold{v}|\le M} \lbold^{\blold{v}}\blold{v}$, and $\ge$, $<$ or $>$ may replace $\le$. For $\blold{u} \in V$, $\lbold|_{\blold{u}} \coloneqq \sum_{\blold{v}\in V} \lbold^{\blold{vu}}\blold{v}$ retains the words of $\lbold$ ending in $\blold{u}$ and deletes that suffix.
    \item \emph{Tensor product.} It extends concatenation bilinearly: $(\lbold\otimes\X)^{\blold{v}} \coloneqq \sum_{\blold{v_1 v_2} = \blold{v}} \lbold^{\blold{v_1}}\X^{\blold{v_2}}$, a finite sum. We also write $\lbold\blold{v}$ for $\lbold\otimes\blold{v}$, $\lbold^{\otimes 0} \coloneqq \empt$, $\lbold^{\otimes(n+1)} \coloneqq \lbold^{\otimes n} \otimes \lbold$, and set $\exp_\otimes(\lbold) \coloneqq \sum_{n \ge 0} \lbold^{\otimes n}/n!$ whenever $\lbold^{\empt} = 0$, so that each coefficient is a finite sum.
    \item \emph{Bracket.} For $\lbold \in T(\R^{d+1})$ and $\X \in T((\R^{d+1}))$, $\brackets{\lbold}{\X} \coloneqq \sum_{\blold{v}} \lbold^{\blold{v}}\X^{\blold{v}}$, a finite sum. When both are infinite, the series may diverge; for $\X$ a signature, this is the object of Sections~\ref{sec:signature} and~\ref{sec:truncation}.
\end{itemize}

\begin{definition}[Shuffle product]\label{def:shuffle}
    The shuffle product is defined on words recursively by 
    \begin{align*}
        (\blold{vi}) \shuffle (\blold{wj}) &= \big(\blold{v} \shuffle (\blold{wj})\big)\blold{i} + \big((\blold{vi}) \shuffle \blold{w}\big)\blold{j},\\
        \blold{v} \shuffle \empt &= \empt \shuffle \blold{v} = \blold{v},
    \end{align*}
    for $\blold{v}, \blold{w} \in V$ and $\blold{i}, \blold{j}$ letters, and extended bilinearly to $T((\R^{d+1}))$.
\end{definition}

Concretely, $\blold{v} \shuffle \blold{w}$ sums, with multiplicity, all interleavings of $\blold{v}$ and $\blold{w}$ that preserve the order of the letters of each word; for instance, $\blold{0} \shuffle \blold{1} = \blold{01} + \blold{10}$ and $\blold{1} \shuffle \blold{1} = 2\,\blold{11}$. The shuffle product is commutative and associative, and adds lengths, so that each coefficient of $\lbold_1 \shuffle \lbold_2$ is a finite sum. We write
$\lbold^{\shuffle 0} \coloneqq \empt$, $\lbold^{\shuffle (n+1)} \coloneqq \lbold^{\shuffle n} \shuffle \lbold$ and $\exp_{\shuffle}(\lbold) \coloneqq \sum_{n \ge 0}\lbold^{\shuffle n}/n!$.

Finally, the \emph{absolute value} $|\lbold| \coloneqq \sum_{\blold{v}}|\lbold^{\blold{v}}|\,\blold{v}$ is taken coordinatewise, and we write $\lbold^1 \preceq \lbold^2$ if $|\lbold^{1,\blold{v}}| \le \lbold^{2,\blold{v}}$ for every $\blold{v} \in V$; in particular $\lbold^{2,\blold{v}} \ge 0$. Thus $\lbold \preceq |\lbold|$, and $\lbold^{\le M} \preceq \lbold$ whenever $\lbold$ has non-negative coefficients. Since the shuffle product has non-negative coefficients on words, domination is preserved by it: $\lbold_1 \preceq \lbold_1^*$ and $\lbold_2 \preceq \lbold_2^*$ imply $\lbold_1 \shuffle \lbold_2 \preceq \lbold_1^* \shuffle \lbold_2^*$.

\subsection{The signature and its linear functionals}\label{sec:signature}

Fix a filtered probability space $(\Omega, \mathcal{F}, \F, \Prob)$ satisfying the usual conditions. For a continuous semimartingale $X$ and a progressively measurable process $Y$ with $\int_0^T Y_s^2\,\diff [X]_s < \infty$ almost surely, we write $\int Y\,\diff X$ for the Itō integral. When $Y$ is also a continuous semimartingale, the Stratonovich integral is $\int Y\circ\diff X \coloneqq \int Y\,\diff X + \frac12[X,Y]$; it is again a continuous semimartingale.

\begin{definition}[Signature]\label{def:signatures}
    Let $X = (X^0, \dots, X^d)$ be a continuous semimartingale on $\R^{d+1}$. Its signature $\X_{s,t} \in T((\R^{d+1}))$, $0 \le s \le t \le T$, is defined coordinatewise by $\X^{\empt}_{s,t} \coloneqq 1$ and
    \begin{equation}\label{eq:signature}
        \X^{\blold{vi}}_{s,t} \coloneqq \int_s^t \X^{\blold{v}}_{s,r} \circ \diff X^i_r, \qquad \blold{v} \in V,\ i \in \intset{0}{d}.
    \end{equation}
\end{definition}

We write $\X_t \coloneqq \X_{0,t}$. By induction on $|\blold{v}|$, every $\X^{\blold{v}}_{s,\cdot}$ is a continuous semimartingale, so that~\eqref{eq:signature} is well defined. The coordinate $\X^{\blold{i_1 \cdots i_m}}_{s,t}$ is thus the iterated integral of $X^{i_1}, \dots, X^{i_m}$, in this order, over $s < r_1 < \cdots < r_m < t$. For a single letter, $\X^{\blold{i}^{\otimes m}}_{s,t} = (X^i_t - X^i_s)^m/m!$: linear functionals of the signature generalise power series in the increment to path-dependent quantities. Our driving path is the time-extended Brownian motion $\Wh = (t, W_t)$ of Section~\ref{sec:introduction}, whose signature we denote $\WH$. Its first coordinates are
\[
    \WH^{\blold{0}}_t = t, \quad
    \WH^{\blold{i}}_t = W^i_t, \quad
    \WH^{\blold{ii}}_t = \tfrac12 (W^i_t)^2, \quad
    \WH^{\blold{0i}}_t = \int_0^t s\,\diff W^i_s, \quad
    \WH^{\blold{i0}}_t = \int_0^t W^i_s\,\diff s,
\]
for $i \ge 1$, and more generally $\WH^{\blold{0}^{\otimes n}\blold{i}}_t = \int_0^t \frac{s^n}{n!}\,\mathrm dW^i_s$.

\medskip

\emph{Chen's identity}~\citep*{ChenKuoTsai1954} expresses the signature over $[s,t]$ through those over $[s,u]$ and $[u,t]$:
\begin{equation}\label{eq:chen-identity}
    \X_{s,t} = \X_{s,u} \otimes \X_{u,t}, \qquad\text{that is,}\qquad \X^{\blold{v}}_{s,t} = \textstyle{\sum_{\blold{v_1 v_2} = \blold{v}}} \X^{\blold{v_1}}_{s,u}\, \X^{\blold{v_2}}_{u,t}, \qquad s \le u \le t.
\end{equation}
For $\lbold \in T(\R^{d+1})$, $\brackets{\lbold}{\X_t}$ is a finite linear combination of iterated integrals. For an infinite $\lbold$, the series must converge; following \citet*{CuchieroSvalutoFerroTeichmann2026, AbiJaberGerardHuang2024}, we define the \emph{admissible set}
\[
    \mathcal{A}(\X) \coloneqq \Big\{ \lbold \in T((\R^{d+1})) : \norm{\lbold}^{\mathcal{A}(\X)}_t \coloneqq \sum_{m \ge 0} \Big| \sum_{\blold{v} \in V_m} \lbold^{\blold{v}}\, \X^{\blold{v}}_t \Big| < \infty \text{ for all } t \in [0,T] \text{ a.s.} \Big\},
\]
on which $\brackets{\lbold}{\X_t} \coloneqq \sum_{m \ge 0}\sum_{\blold{v}\in V_m} \lbold^{\blold{v}}\X^{\blold{v}}_t$ is well defined. It contains $T(\R^{d+1})$. The key algebraic property of signatures is that products of linear functionals are again linear functionals.

\begin{proposition}[Shuffle identity, {\citet*[Proposition~3.1]{AbiJaberGerardHuang2024}}]\label{prop:shuffle-identity}
    If $\lbold_1, \lbold_2 \in \mathcal{A}(\X)$, then 
    \[
        \lbold_1 \shuffle \lbold_2 \in \mathcal{A}(\X) \qquad \text{and} \qquad \brackets{\lbold_1}{\X_t}\, \brackets{\lbold_2}{\X_t} = \brackets{\lbold_1 \shuffle \lbold_2}{\X_t}.
    \]
\end{proposition}

For instance, $\blold{1} \shuffle \blold{1} = 2\,\blold{11}$ gives $W_t \cdot W_t = 2\cdot\frac12 W_t^2$, and $\blold{0} \shuffle \blold{1} = \blold{01} + \blold{10}$ is the integration by parts formula $tW_t = \int_0^t s\,\diff W_s + \int_0^t W_s\,\diff s$.

\subsection{Discretised signatures}\label{sec:discretised-signatures}

The signature of Brownian motion cannot be computed exactly, and in practice it is replaced by the signature of a discretised path. For $N \in \N^*$, let $h \coloneqq T/N$ and $t_k \coloneqq kh$, $k \in \intset{0}{N}$, and let $\Wt = \Wt(N)$ be the piecewise linear interpolation of $\Wh$ on this grid,
\begin{equation}\label{eq:discretised-signature}
    \Wt_t \coloneqq \Wh_{t_{k-1}} + \frac{t - t_{k-1}}{h} \big(\Wh_{t_k} - \Wh_{t_{k-1}}\big), \qquad t \in (t_{k-1}, t_k].
\end{equation}
The path $\Wt$ is continuous and of bounded variation, but not adapted: on $(t_{k-1},t_k)$ it uses $\Wh_{t_k}$. Its signature $\WT = \WT(N)$ is nevertheless well defined, pathwise, by~\eqref{eq:signature} with Riemann--Stieltjes integrals, and Chen's identity~\eqref{eq:chen-identity} holds for it as well. On a single step, the signature of an affine path is the tensor exponential of its increment, so that, for $t_{k-1} \le s \le t \le t_k$,
\[
    \WT_{s,t} = \exp_\otimes\Big\{ \frac{t-s}{h} \big(\Wh_{t_k} - \Wh_{t_{k-1}}\big) \Big\}, \qquad\text{that is,}\qquad \WT^{\blold{i_1 \cdots i_m}}_{s,t} = \frac{1}{m!}\prod_{l=1}^m \frac{t-s}{h}\big(\Wh^{i_l}_{t_k} - \Wh^{i_l}_{t_{k-1}}\big),
\]
and~\eqref{eq:chen-identity} yields $\WT_{s,t}$ for all $s \le t$ as a product of such exponentials. This is exactly what standard signature libraries, such as \texttt{iisignature}~\citep*{ReizensteinGraham2020}, \texttt{esig}\footnote{\url{https://pypi.org/project/esig/}} and \texttt{signatory}~\citep*{KidgerLyons2021}, compute: they take a finite stream of points, interpret it as the piecewise linear path through them, and return its signature truncated at a prescribed level. We set $\WT(\infty) \coloneqq \WH$.

As an illustration, for $d = 1$ and $N = 1$, the interpolation is the straight line from $0$ to $(T, W_T)$, so that $\WT^{\blold{11}}_T = \frac12 W_T^2 = \WH^{\blold{11}}_T$, whereas $\WT^{\blold{01}}_T = \frac12 T W_T$ differs from $\WH^{\blold{01}}_T = \int_0^T s\,\diff W_s$ by the centred Gaussian variable $\int_0^T (s - T/2)\,\diff W_s$, of variance $T^3/12$.

\begin{remark}[Shuffle identity for paths of bounded variation]\label{rem:shuffle-BV}
    For $\lbold_1, \lbold_2 \in T(\R^{d+1})$, the shuffle identity holds pathwise for the signature of any continuous path of bounded variation \citep*{Ree1958}; in particular, it holds for $\WT(N)$, $N \in \N^*$. Its extension to infinite coefficient families is given in Lemma~\ref{lem:shuffle-L1}.
\end{remark}

\subsection{Moment estimates}\label{sec:moments}

For $p \ge 1$, $t \ge 0$ and $\blold{v} \in V$, set $c_d \coloneqq 1 + d/2$,
\[
  \psi_t(\blold{v}) \coloneqq \indic{\norm{\blold{v}} \in 2\N} \, \frac{(c_d t)^{\norm{\blold{v}}/2}}{(\norm{\blold{v}}/2)!} \qquad \text{and} \qquad \psi_{p,t}(\blold{v}) \coloneqq (p+2)^{|\blold{v}|}\, \frac{(c_d\,t)^{\norm{\blold{v}}/2}}{\lfloor \norm{\blold{v}}/2 \rfloor!}.
\]
Up to geometric factors, a component of degree $n$ is thus of size $t^{n/2}/\sqrt{n!}$; this factorial decay is what allows infinite linear combinations to converge.

\begin{proposition}[Expected signature]\label{prop:expected-signature-bound}
    For every $N \in \Nbar^*$, $\blold{v} \in V$ and $0 \le s \le t \le T$,
    \[
        0 \le \E\big[\WT^{\blold{v}}_{s,t}(N)\big] \le \psi_{t-s}(\blold{v}).
    \]
\end{proposition}

\begin{proof}
    For $N = \infty$, Fawcett's formula \citep*{Fawcett2003, LyonsVictoir2004} gives $\E[\WH_{s,t}] = \exp_\otimes\{(t-s)(\blold{0} + \frac12\sum_{i \ge 1}\blold{ii})\}$. The tensor in the exponential has non-negative coefficients, carried by words of degree $2$, summing to $c_d(t-s)$. Hence $\E[\WH_{s,t}]$ has non-negative coefficients, vanishing in odd degree, and those of degree $2k$ sum to $(c_d(t-s))^k/k!$, which bounds each of them. For $N < \infty$, see Appendix~\ref{app:expected-signature-bound}.
\end{proof}

\begin{proposition}[Moments of the signature]\label{prop:signature-moment-bound}
    For every $p \ge 1$, $N \in \Nbar^*$, $\blold{v} \in V$ and $0 \le s \le t \le T$, 
    \[
        \textnorm{\WT^{\blold{v}}_{s,t}(N)}_p \le \psi_{p,t-s}(\blold{v}) \le \psi_{p,T}(\blold{v}).
    \]
\end{proposition}

\begin{proof}
    Set $q \coloneqq 2\lceil p/2 \rceil$, the smallest even integer not smaller than $p$, so that $q \le p+2$, and write $\X \coloneqq \WT(N)$. Since $q$ is even, $\norm{\X^{\blold{v}}_{s,t}}_p^q \le \E[(\X^{\blold{v}}_{s,t})^q]$. By the shuffle identity (Proposition~\ref{prop:shuffle-identity} for $N = \infty$, Remark~\ref{rem:shuffle-BV} for $N < \infty$), $\E[(\X^{\blold{v}}_{s,t})^q] = \sum_{\blold{u}} \brackets{\blold{u}}{\blold{v}^{\shuffle q}}\, \E[\X^{\blold{u}}_{s,t}]$, a finite sum with non-negative coefficients, over words $\blold{u}$ of degree $q\norm{\blold{v}}$. By Proposition~\ref{prop:expected-signature-bound}, each expectation is at most $\psi_{t-s}(\blold{u})$, and the coefficients sum to at most $q^{q|\blold{v}|}$ (Lemma~\ref{lem:count-v-shuffle_p-elements}). Hence $\norm{\X^{\blold{v}}_{s,t}}_p \le q^{|\blold{v}|}(c_d(t-s))^{\norm{\blold{v}}/2} [(q\norm{\blold{v}}/2)!]^{-1/q}$. Finally, with $a \coloneqq \lfloor\norm{\blold{v}}/2\rfloor$, $(q\norm{\blold{v}}/2)! \ge (qa)! \ge (a!)^q$, and $q \le p+2$ concludes.
\end{proof}

We will mostly use the following consequence. Since $n! \le 2^n(n/2+1)(\lfloor n/2\rfloor!)^2$, which follows from $\binom{2k}{k} \le 4^k$ for $n = 2k$, and from $2k+1 \le 2(n/2+1)$ for $n = 2k+1$,
\begin{equation}\label{eq:psi-bound}
    \psi_{p,t}(\blold{v}) \le (p+2)^{|\blold{v}|}\,(2c_d t)^{\norm{\blold{v}}/2} \sqrt{\frac{\norm{\blold{v}}/2+1}{\norm{\blold{v}}!}} .
\end{equation}
\section{Main results} \label{sec:main-results}

This section answers the three questions raised in Section \ref{sec:problem-formulation}. We first identify the coefficient families for which the representation~\eqref{eq:signature-representation} can be controlled, and then state the truncation and discretisation estimates, in the strong and in the weak sense.

\medskip

Throughout, we fix $T > 0$ and $D \in \N^*$, and we equip $\R^D$ with the Euclidean norm $\tnorm{\cdot}$. The test functions are taken in two classes: $\mathcal{P}(\R^D, \R)$, the class of polynomial functions of the form $x \in \R^D \mapsto \sum_{\alpha \in \N^D, |\alpha| \le m} c_\alpha x^\alpha$ for some $m \in \N$, and $\mathcal{H}(\R^D, \R)$, the class of Hölder continuous functions from $\R^D$ to $\R$. We introduce the unified class $\mathcal{G}^D \coloneqq \mathcal{P}(\R^D,\R) \cup \mathcal H(\R^D,\R)$. For $g \in \mathcal{G}^D$, we call \emph{regularity exponent} (abbreviated by r.e.) of $g$ the number $a = 1$ if $g \in \mathcal{P}(\R^D,\R)$, and any $a \in (0,1]$ such that $g$ is $a$-Hölder
otherwise.

\medskip

We now introduce the space of coefficients that will be used throughout the paper:
\[
    \mathcal{L}^1 \coloneqq \left\{ \lbold \in T((\R^{d+1})) \ : \ \exists \varepsilon \in (0, 1/2], \ |\lbold^{\blold{v}}| \lesssim \varphi_\varepsilon(\blold{v}), \, \blold{v} \in V \right\},
\]
where $\varphi_\varepsilon(\blold{v}) \coloneqq (\norm{\blold{v}}!)^{\frac{1}{2} - \varepsilon}$. For $\lbold \in \mathcal{L}^1$, we call any $\varepsilon \in (0, 1/2]$ for which $|\lbold^{\blold{v}}| \lesssim \varphi_\varepsilon(\blold{v})$, $\blold{v} \in V$ holds a \emph{summability exponent} of $\lbold$; since $\varphi_\varepsilon$ is non-increasing in $\varepsilon$, the larger the exponent, the stronger the condition. The purpose of this coefficient space is that the growth of $\lbold^{\blold{v}}$ is dominated by the decay of $\psi_{p,t}(\blold{v})$, so that the products $\lbold^{\blold{v}}\psi_{p,t}(\blold{v})$ remain summable. In addition, we define $\mathcal{L}^D \coloneqq [\mathcal{L}^1]^D$. Since the dimension $d$ is fixed, we do not indicate the underlying space $\R^{d+1}$ in this notation, and write $\mathcal{L}^1(\R^{d'+1})$ only when another dimension $d'$ is used, as in Section~\ref{sec:pricing}.

\begin{remark}[Exponentially bounded coefficients] \label{rem:Aexp}
    The space $\mathcal{L}^1$ is non-empty, as it contains all truncated tensors $T(\R^{d+1})$. It also includes non-truncated elements, for instance those in 
    \[
        \mathcal{A}_{\exp} \coloneqq \Big\{ \lbold \, : \, \exists C \ge 0, \, \lbold \preceq C \exp_{\shuffle}\big(C \, \textstyle{\sum_{i = 0}^d}\blold{i}\big) \Big\},
    \]
    which plays a key role in the signature expansions of Volterra processes and delay equations, as shown by \citet*{AbiJaberGerardHuang2024}.
    Indeed, let $\chi(C) \coloneqq \exp_{\shuffle}\big(C \sum_{\blold{i}}\blold{i}\big)$, we have $\chi^{\blold{v}}(C) = \brackets{\blold{v}}{\exp_{\shuffle}\left(C \textstyle\sum_{\blold{i}}\blold{i}\right)} = \sum_n C^n\brackets{\blold{v}}{\left(\Sigma_{\blold{i}}\blold{i} \right)^{\shuffle n}}/n! = C^{|\blold{v}|}\textbrackets{\blold{v}}{\left(\Sigma_{\blold{i}}\blold{i}\right)^{\shuffle |\blold{v}|}/|\blold{v}|!}$ and, using the identity $(\Sigma_{\blold{i}}\blold{i})^{\shuffle n} / n! = (\Sigma_{\blold{i}}\blold{i})^{\otimes n}$ (see the proof of Proposition 2.3 in \citet{AbiJaberGerardHuang2024}), 
    \[
        \chi^{\blold{v}}(C) = C^{|\blold{v}|}\sum_{\blold{i_1},\dots, \blold{i_{|v|}}}\brackets{\blold{v}}{\blold{i_1} \otimes \cdots \otimes \blold{i_{|v|}}} = C^{|\blold{v}|}.
    \]
    Since $C^{|\blold{v}|} \lesssim \varphi_{1/4}(\blold{v})$, $\blold{v} \in V$, we deduce that $\chi(C) \in \mathcal{L}^1$ and, consequently, $\mathcal{A}_{\exp} \subset \mathcal{L}^1$.
\end{remark}

The space $\mathcal{L}^1$ is stable under several natural operations; the proof is deferred to Appendix~\ref{app:L1-closeness}.

\begin{proposition}[Stability of $\mathcal{L}^1$] \label{prop:L1-stability}
    Let $\lbold, \lbold_1, \lbold_2 \in \mathcal{L}^1$. Then, for every $\blold{u} \in V$ and $\lambda_1, \lambda_2 \in \R$,
    \[
        \lambda_1 \lbold_1 + \lambda_2 \lbold_2, \quad \lbold_1 \shuffle \lbold_2, \quad \lbold_1 \otimes \lbold_2, \quad \lbold |_{\blold{u}} \quad \in \, \mathcal{L}^1.
    \]
\end{proposition}

We extend $\mathcal{L}^1$ to time-indexed families by defining
\[
    \mathcal{L}^1_T \coloneqq \left\{ \lbold \in T((\R^{d+1}))^{[0, T]} \, : \, \exists \lbold^* \in \mathcal{L}^1, \, \lbold_t \preceq \lbold^*, \, t \in [0, T] \right\}.
\]
Similarly, we define $\mathcal{L}^D_T \coloneqq \big[\mathcal{L}^1_T\big]^D$. Proposition \ref{prop:L1-stability} extends directly to $\mathcal{L}^1_T$.

\medskip

In this section, we fix $\lbold \in \mathcal{L}^D_T$, $M \in \Nbar$ and $N \in \Nbar^*$, and set
\[
    X_t(M,N) \coloneqq \Big( \brackets{\lbold^i_t}{\WT^{\le M}_t(N)} \Big)_{i \in \{1,\dots,D\}}, \qquad t \in [0,T],
\]
with $\WT^{\le \infty} \coloneqq \WT$. Thus $X_t \coloneqq X_t(\infty,\infty)$ is the target, $X_t(M) \coloneqq X_t(M,\infty)$ its truncation, $X_t(\infty,N)$ the full linear functional of the discretised signature, and $X_t(M,N)$ the computable approximation. For $M = \infty$ the bracket is an infinite series; Corollary~\ref{cor:well-posed} shows that it converges absolutely, almost surely and in every $L^p$, uniformly in $t$ and $N$, and Lemma~\ref{lem:shuffle-L1} extends the shuffle to $\mathcal{L}^1$.

\begin{figure}[H]
\centering
\begin{tikzpicture}[
    >=stealth,
    every node/.style={inner sep=4pt},
    lab/.style={align=center, font=\small}
  ]
  
  \node (A) at (0,0)  {$X_t(\infty,\infty)$};
  \node (B) at (7,0)  {$X_t(M,\infty)$};
  \node (C) at (0,-4) {$X_t(\infty,N)$};
  \node (D) at (7,-4) {$X_t(M,N)$};

  \draw[->] (B) -- node[above, lab]
    {Truncation\\ strong: $r_M^a$ ; weak: $r_M^a$} (A);
  \draw[->] (D) -- node[below, lab]
    {Truncation\\ strong: $r_M^a$ ; weak: $r_M^a$} (C);

  \draw[->] (C) -- node[left, lab]
    {Discretisation\\ strong: $N^{-a/2}$\\ weak: $N^{-1}$ ($g \in \mathcal P$)} (A);
  \draw[->] (D) -- node[right, lab]
    {Discretisation\\ strong: $N^{-a/2}$\\ weak: $N^{-1}$ ($g \in \mathcal P$)} (B);

  \draw[->] (D) -- node[fill=white, inner sep=2pt] {$(\star)$} (A);
\end{tikzpicture}
\caption{Truncation (horizontal) and discretisation (vertical) errors, where $a$ is a regularity exponent of $g$; strong: error in $L^p$, weak: weak error. Each arrow points from an approximation to its limit. Horizontal rates are uniform in $N \in \Nbar^*$, vertical rates in $M \in \Nbar$. For Hölder $g$, the weak discretisation rate is $N^{-a/2}$. $(\star)$: total error, Corollaries~\ref{cor:total-error} (strong) and~\ref{cor:total-error-weak} (weak).} \label{fig:error-square}
\end{figure}
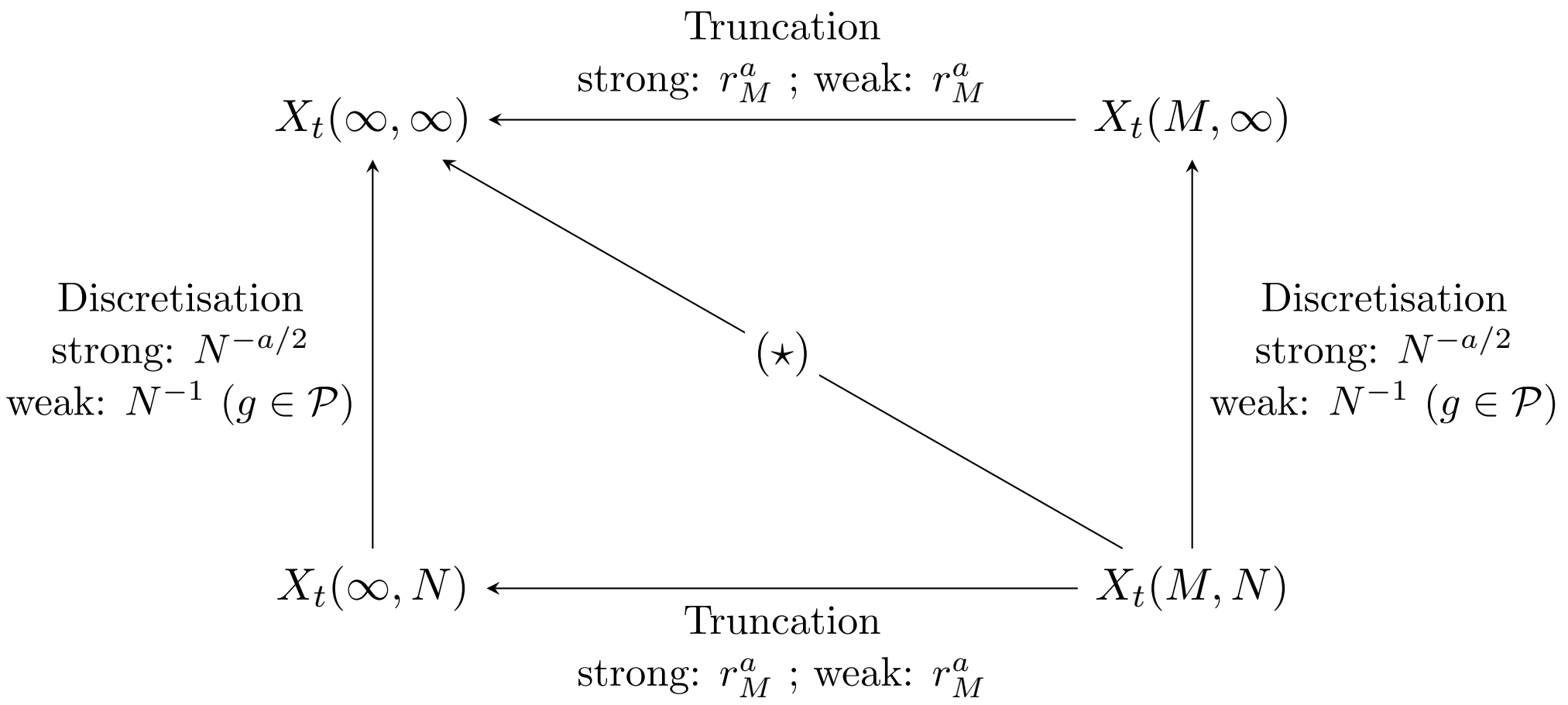

For a test function $g$, we call \emph{truncation error} the difference between $g(X_t(\infty,N))$ and $g(X_t(M,N))$, at a fixed discretisation level $N \in \Nbar^*$, and \emph{discretisation error} the difference between $g(X_t(M,\infty))$ and $g(X_t(M,N))$, at a fixed truncation level $M \in \Nbar$. Each is measured in the \emph{strong} sense, through $\norm{\, \cdot \,}_p$, and in the \emph{weak} sense, through $\absE{\, \cdot \,}$. Throughout, the truncation rate is written
\[
  r_M \coloneqq B^M \, (M!)^{-\varepsilon}, \quad M \in \N, \qquad r_\infty \coloneqq 0,
\]
where $\varepsilon$ is a summability exponent of the coefficients and $B > 0$ a constant given by each statement. Since factorials dominate geometric sequences, $\sup_{M} r_M < \infty$ and $r_M \to 0$ as $M \to \infty$. Figure~\ref{fig:error-square} summarises the rates obtained in this section.

\subsection{Strong errors} \label{sec:strong-errors}

We first control the strong errors, starting with truncation.

\begin{theorem}[Strong truncation error] \label{th:Lp-truncation-bound}
    Let $\lbold^\ast \in \mathcal{L}^1$ dominate the coefficients, that is, $\lbold^i_t \preceq \lbold^\ast$ for every $t \in [0, T]$ and $i \in \{1, \dots, D\}$, with $\varepsilon \in (0, 1/2]$ a summability exponent of $\lbold^\ast$. Let $g \in \scriptG^D$ with regularity exponent $a$ and let $p \geq 1$. Then, there exist $B, C > 0$, depending only on $d$, $D$, $T$, $p$, $g$, $\varepsilon$ and $\lbold^\ast$, such that
    \begin{equation} \label{eq:truncation-rate}
        \norm{g\big(X_t(\infty, N)\big) - g\big(X_t(M, N)\big)}_p \le C \, r_M^a, \qquad t \in [0, T], \ M \in \N, \ N \in \Nbar^*.
    \end{equation}
\end{theorem}

\begin{remark}
    As an illustration of Theorem~\ref{th:Lp-truncation-bound}, let $d = D = 1$, and set $f(x) \coloneqq \sum_{n} a_n x^n / n!$ with $(a_n)_{n \in \N}$ real coefficients. Since $\WH^{\blold{1}^{\otimes n}}_t = W^n_t / n!$, we have $f(W_t) = \textbrackets{\lbold}{\WH_t}$ and $\sum_{n \leq M} a_n W^n_t / n! = \textbrackets{\lbold^{\leq M}}{\WH_t}$, with $\lbold \coloneqq \sum_n a_n \blold{1}^{\otimes n}$. If $|a_n| \lesssim [n!]^{1/2 - \varepsilon}$, then $\lbold \in \mathcal{L}^1(\R^2)$, and Theorem~\ref{th:Lp-truncation-bound} bounds the truncation error of the power series along a Brownian path for a test function $g$ with regularity exponent $a$ by $C r_M^a$.
\end{remark}

The discretisation error is controlled uniformly in the truncation level, including $M = \infty$.

\begin{theorem}[Strong discretisation error]\label{th:Lp-discretisation-bound}
    Under the assumptions of Theorem~\ref{th:Lp-truncation-bound}, there exists $C > 0$, depending only on $d, D, T, p, g, \varepsilon$ and $\lbold^*$, such that
    \[
        \norm{ g\big(X_t(M,\infty)\big) - g\big(X_t(M,N)\big) }_p \le C\, N^{-a/2}, \qquad t \in [0,T],\ M \in \Nbar,\ N \in \Nbar^*.
    \]
    If moreover $d = 1$, the bound improves to $C N^{-a}$ at the grid points $t \in \{t_k : 0 \le k \le N\}$.
\end{theorem}

The two strong estimates combine along either path of Figure~\ref{fig:error-square}. For $a \in (0,1]$, $M \in \Nbar$ and $N \in \Nbar^*$, set $\mathcal{E}_a(M,N) \coloneqq r_M^a + N^{-a/2}$, with $N^{-a/2} \coloneqq 0$ if $N = \infty$. Since $u \mapsto u^\theta$ is subadditive on $\R_+$ for $\theta \in (0,1]$, $\mathcal{E}_a(M,N)^\theta \le \mathcal{E}_{\theta a}(M,N)$.

\begin{corollary}[Strong total error]\label{cor:total-error}
    Under the assumptions of Theorem~\ref{th:Lp-truncation-bound}, there exist $B, C > 0$, depending only on $d, D, T, p, g, \varepsilon$ and $\lbold^*$, such that
    \[
        \norm{g(X_t) - g\big(X_t(M,N)\big)}_p \le C\, \mathcal{E}_a(M,N), \qquad t \in [0,T],\ M \in \Nbar,\ N \in \Nbar^* .
    \]
    If $d = 1$, the term $N^{-a/2}$ can be replaced by $N^{-a}$ at the grid points.
\end{corollary}

\begin{proof}
    By the triangle inequality, it suffices to bound $g(X_t(\infty,\infty)) - g(X_t(\infty,N))$, controlled by Theorem~\ref{th:Lp-discretisation-bound} with $M = \infty$, and $g(X_t(\infty,N)) - g(X_t(M,N))$, controlled by Theorem~\ref{th:Lp-truncation-bound} uniformly in $N$.
\end{proof}

In particular, $M$ and $N$ can be chosen independently: the truncation level is set by the target accuracy alone, and $N$ is then taken of order $r_M^{-2}$ to balance the two terms in $L^p$, whatever the regularity exponent $a$.

\subsection{Weak errors} \label{sec:weak-errors}

We now turn to the third question of Section~\ref{sec:problem-formulation}. For the discretisation, the weak error improves on the strong one (Theorem~\ref{th:weak-discretisation-bound} below); for the truncation, it does not. This may surprise a reader familiar with~\citet*{KloedenPlaten1992}, where the weak Itō--Taylor expansions converge at twice the order of the strong ones. There, ``order'' refers to the power of the time horizon in the remainder: for an expansion truncated at the same level $M$, the weak remainder is of order $t^{M+1}$ whereas the strong one is of order $t^{(M+1)/2}$. The weak expansion therefore gains a factor $2$ in the exponent of $t$, that is, it is more accurate for small $t$. In our setting, the bound of Theorem~\ref{th:Lp-truncation-bound} (through the factor $t^{(M+1)/2}$ of Lemma~\ref{lem:truncation}) and the example below show that this gain does not occur: the weak truncation error has the same exponent in $t$ and the same factorial decay in $M$ as the strong one. The reason is that every component of the Itō signature of a Brownian motion involving at least one Brownian letter has zero mean (it is a martingale when its last letter is Brownian), so that cancellations separate the weak rate from the strong one. No such cancellation holds for Stratonovich signatures, and the remainder can have a constant sign, in which case the weak error coincides with the strong $L^1$ error. Let $d = D = 1$. The simplest example is
\[
    X_t = \cosh(W_t) = \brackets{\lbold}{\WH_t}, \qquad \lbold \coloneqq \sum_{n=0}^{\infty} a_n \blold{1}^{\otimes n}, \quad a_{2n} = 1, \ a_{2n+1} = 0,
\]
for which $X_t \geq X_t(M)$, hence $\E[X_t - X_t(M)] = \norm{X_t - X_t(M)}_1$, for every $t \in [0, T]$ and $M \in \N$. The same holds for the discretised signature: since $\WT^{\blold{1}^{\otimes n}}_t(N) = \Wt_t(N)^n/n!$, $X_t(\infty,N) = \cosh(\Wt^1_t(N)) \ge X_t(M,N)$ for every $N \in \Nbar^*$. The absence of a weak gain is therefore not an artefact of the continuous signature.

\begin{corollary}[Weak truncation error] \label{cor:weak-truncation-bound}
    Under the assumptions of Theorem~\ref{th:Lp-truncation-bound}, and with the same constants $B$ and $C$ (for $p = 1$),
    \[
        \absE{g\big(X_t(\infty, N)\big) - g\big(X_t(M, N)\big)} \le C \, r_M^a, \qquad t \in [0, T], \ M \in \N, \ N \in \Nbar^*.
    \]
\end{corollary}

The rate of Corollary~\ref{cor:weak-truncation-bound}, and hence that of Theorem~\ref{th:Lp-truncation-bound}, cannot be improved, exponent $\varepsilon$ included, and uniformly in $N$, as the following example shows.

\begin{proposition}[Sharpness of the truncation rate]\label{prop:truncation-sharpness}
    Let $d = D = 1$ and $\varepsilon \in (0,1/2]$, and consider
    \[
        \lbold^{(\varepsilon)} \coloneqq \textstyle{\sum_{n \ge 0}} (n!)^{1/2 - \varepsilon}\, \blold{1}^{\otimes n}.
    \]
    Since $\norm{\blold{1}^{\otimes n}} = n$, the coefficients are exactly $\varphi_\varepsilon(\blold{1}^{\otimes n})$, so that $\lbold^{(\varepsilon)} \in \mathcal L^1(\mathbb R^2)$ with summability exponent $\varepsilon$. For $N \in \Nbar^*$, let $X_t(\infty,N) \coloneqq \textbrackets{\lbold^{(\varepsilon)}}{\WT_t(N)}$ and $X_t(M,N) \coloneqq \textbrackets{\lbold^{(\varepsilon)}}{\WT^{\le M}_t(N)}$. Then, for every $t \in (0,T]$, $p \ge 1$, $M \ge 2$ and $N \in \Nbar^*$, writing $c_t \coloneqq 2 + |\log(t^2/(2T))|$, we have
    \begin{equation}\label{eq:sharpness}
        \norm{X_t(\infty,N) - X_t(M,N)}_p
        \;\ge\; \E\big[X_t(\infty,N) - X_t(M,N)\big]
        \;\ge\; \exp\big\{-\varepsilon M \log M - c_t M\big\}.
    \end{equation}
    Consequently, for this family and $g = \mathrm{Id}$, both the weak and the strong truncation errors satisfy, uniformly in $N \in \Nbar^*$, $\log \mathrm{error} = -\varepsilon M \log M + O(M)$, $M \to \infty$, and the exponent $\varepsilon$ in Theorem~\ref{th:Lp-truncation-bound} and Corollary~\ref{cor:weak-truncation-bound} cannot be improved, neither for the signature nor for its discretisation.
\end{proposition}

The weak discretisation error, in contrast, exhibits an improved rate for polynomial test functions. This follows from the refined weak discretisation estimates for the signature component established in Proposition~\ref{prop:weak-element-discretisation}.

\begin{theorem}[Weak discretisation error] \label{th:weak-discretisation-bound} 
    Under the assumptions of Theorem~\ref{th:Lp-truncation-bound}, there exists $C > 0$, depending only on $d, D, T, g, \varepsilon$ and $\lbold^*$, such that
    \[
        \Big| \E\big[g\big(X_t(M,\infty)\big) - g\big(X_t(M,N)\big)\big] \Big| \le \begin{cases}
            C\, N^{-1}, & g \in \mathcal{P}(\R^D, \R),\\
            C\, N^{-a/2}, & g \in \mathcal H(\R^D, \R),
        \end{cases} \qquad t \in [0,T],\ M \in \Nbar,\ N \in \N^* .
    \]
    If moreover $d = 1$, the Hölder bound improves to $C N^{-a}$ at the grid points.
\end{theorem}

\begin{remark}[Optimality of the discretisation rates]\label{rem:optimality-statement}
    The rates of Theorems~\ref{th:Lp-discretisation-bound} and~\ref{th:weak-discretisation-bound} are optimal for polynomial test functions, for every truncation level $M \ge 2$, and single words already attain them:
    \begin{enumerate}[label=(\roman*)]
        \item off the grid, the strong rate $N^{-1/2}$ and the weak rate $N^{-1}$ are attained, already for $d = 1$;
        \item at the grid points and for $d \ge 2$, the strong rate $N^{-1/2}$ cannot be improved by any approximation built on the grid increments alone \citep*{ClarkCameron1980};
        \item at the grid points and for $d = 1$, the improved strong rate $N^{-1}$ is attained.
    \end{enumerate}
    Explicit examples are given in Remark~\ref{rem:discretisation-optimality}. Together with Proposition~\ref{prop:truncation-sharpness}, this shows that all the rates of Sections~\ref{sec:strong-errors} and~\ref{sec:weak-errors} are sharp for polynomial test functions; for Hölder test functions with $a < 1$, sharpness is left open.
\end{remark}

Combining Corollary~\ref{cor:weak-truncation-bound} and Theorem~\ref{th:weak-discretisation-bound} gives the weak total error.

\begin{corollary}[Weak total error]\label{cor:total-error-weak}
    Under the assumptions of Theorem~\ref{th:Lp-truncation-bound}, if $g \in \mathcal{P}(\R^D,\R)$, there exist $B, C > 0$, depending only on $d, D, T, g, \varepsilon$ and $\lbold^*$, such that
    \[
        \Big| \E\big[ g(X_t) - g\big(X_t(M,N)\big) \big] \Big| \le C\big( r_M + N^{-1} \big), \qquad t \in [0,T],\ M \in \Nbar,\ N \in \Nbar^*.
    \]
\end{corollary}

\begin{proof}
    As for Corollary~\ref{cor:total-error}, with Theorem~\ref{th:weak-discretisation-bound} and Corollary~\ref{cor:weak-truncation-bound} in place of Theorems~\ref{th:Lp-discretisation-bound} and~\ref{th:Lp-truncation-bound}.
\end{proof}

For a polynomial test function, $N$ of order $r_M^{-1}$ thus suffices to balance the two terms of the weak error.

\section{Applications} \label{sec:example}

The estimates of Section~\ref{sec:main-results} apply to any process whose coefficient family lies in $\mathcal{L}^1_T$. We first transfer them to errors in probability, on distribution functions and on risk measures (Section~\ref{sec:extensions-statements}). We then illustrate them on three examples, ordered by increasing distance from that setting: linear delay equations, whose coefficients lie in $\mathcal{L}^1$ by construction; Gaussian Volterra processes, for which membership must be checked kernel by kernel and may fail; and a signature volatility model, where the estimates
apply to the log-price only.

\subsection{Errors in probability and on risk measures} \label{sec:extensions-statements}

The strong estimates of Section~\ref{sec:main-results} hold for every $p \ge 1$, and this freedom transfers them, almost without loss, to other notions of error; the proofs are given in Section~\ref{sec:extensions}. Errors in probability and on coherent risk measures follow directly. Distribution functions and Value-at-Risk correspond to the indicator test function, which is not Hölder continuous, and their control requires an assumption on the law of the target, never on the law of its approximations. Throughout, $g \in \mathcal G^D$ has $a$ as regularity exponent, and $\varepsilon$ is as in Theorem~\ref{th:Lp-truncation-bound}. We write $F^g_t$ and $F^{g,M,N}_t$ for the distribution functions of $g(X_t)$ and $g(X_t(M,N))$.

\begin{lemma}[Error in probability]\label{lem:probability-bound}
    For every $p \ge 1$, there exist $B, C > 0$, depending only on $d, D, T, p, g, \varepsilon$ and $\lbold^*$, such that
    \[
        \Prob\big( |g(X_t) - g(X_t(M,N))| \ge \eta \big) \le C\, \eta^{-p}\, \mathcal{E}_a(M,N)^p, \qquad t \in [0,T],\ M \in \Nbar,\ N \in \Nbar^*,\ \eta > 0.
    \]
\end{lemma}

\begin{lemma}[Error on distribution functions]\label{lem:cdf-bound}
    Let $I \subset \R$ be open, $J \subset I$ with $\eta_0 \coloneqq \operatorname{dist}(J, \R \setminus I) > 0$, and $\theta \in (0,1)$. Assume that, for every $t \in (0,T]$, $g(X_t)$ admits a density $f^g_t$ on $I$ with $\norm{f^g}_I \coloneqq \sup_{t\in(0,T],\, x\in I} f^g_t(x) < \infty$. Then there exist $B, C > 0$, depending only on $d, D, T, g, \varepsilon, \lbold^*, \theta, \eta_0$ and $\norm{f^g}_I$, such that
    \[
        \sup_{x\in J} \big| F^g_t(x) - F^{g,M,N}_t(x) \big| \le C\, \mathcal{E}_{\theta a}(M,N), \qquad t \in [0,T],\ M \in \Nbar,\ N \in \Nbar^*.
    \]
\end{lemma}

Throughout the rest of this subsection, random variables are \emph{losses}: large values are the unfavourable ones. Two questions must be kept apart: the risk carried by the approximation error, $\rho(Z_1 - Z_2)$, and the error committed on the risk itself, $\rho(Z_1) - \rho(Z_2)$. A risk measure $\rho$ on a space of integrable losses is \emph{coherent} \citep*{ArtznerDelbaenEberHeath1999} if it is monotone, cash-additive, positively homogeneous and sub-additive. We consider coherent, lower semicontinuous risk measures with robust representation $\rho(Z) = \sup_{\Q \in \mathcal{Q}_\rho} \E^{\Q}[Z]$, where $\mathcal{Q}_\rho$ is a set of probability measures absolutely continuous with respect to $\Prob$, and we write $|\rho|[Z] \coloneqq |\rho(Z)| \vee |\rho(-Z)|$ and $(\rho\circ g)(Z) \coloneqq \rho(g(Z))$.

\begin{theorem}[Coherent risk measures]\label{th:coherent-bound}
    Assume that there exist $q \in (1,\infty]$ and $C_\rho > 0$ such that $\norm{\diff \Q / \diff \Prob}_q \le C_\rho$ for every $\Q \in \mathcal{Q}_\rho$. Then there exist $B, C > 0$, depending only on $d, D, T, g, \varepsilon, \lbold^*, q$ and $C_\rho$, such that, for $t \in [0,T]$, $M \in \Nbar$ and $N \in \Nbar^*$,
    \[
        \big| (\rho\circ g)(X_t) - (\rho\circ g)(X_t(M,N)) \big| \le |\rho|\big[ g(X_t) - g(X_t(M,N)) \big] \le C\, \mathcal{E}_a(M,N).
    \]
\end{theorem}

For $\alpha \in (0,1)$, the \emph{Value-at-Risk} of a loss $Z$ at level $\alpha$ is its $\alpha$-quantile, $\VaR_\alpha[Z] \coloneqq \inf\{x \in \R : \Prob(Z \le x) \ge \alpha\}$, and its \emph{Expected Shortfall} is $\ES_\alpha[Z] \coloneqq (1-\alpha)^{-1}\int_\alpha^1 \VaR_u[Z]\,\diff u$. The Expected Shortfall is coherent \citep*{AcerbiTasche2002} and admits the representation $\ES_\alpha[Z] = \sup_{\Q \in \mathcal{Q}_\alpha} \E^{\Q}[Z]$, with $\mathcal{Q}_\alpha \coloneqq \{\Q \ll \Prob : \diff \Q / \diff\Prob \le (1-\alpha)^{-1}\}$ \citep*{FollmerSchied2016, RuszczynskiShapiro2006}. Theorem~\ref{th:coherent-bound} therefore applies with $q = \infty$ and $C_\rho = (1-\alpha)^{-1}$. Value-at-Risk is not sub-additive, and the two questions above separate.

\begin{theorem}[Value-at-Risk]\label{th:var-bound}
    Let $\alpha \in (0,1)$. 
    \begin{enumerate}[label=(\roman*)]
    \item If $\alpha \in (1/2,1)$, there exist $B, C > 0$, depending only on $d, D, T, g, \varepsilon$ and $\lbold^*$, such that
    \[
        |\VaR_\alpha|\big[ g(X_t) - g(X_t(M,N)) \big] \le C\, \mathcal{E}_a(M,N), \qquad t \in [0,T],\ M \in \Nbar,\ N \in \Nbar^*.
    \]
    \item Let $\theta \in (0,1)$, and assume that there exist $\gamma > 0$, an interval $\mathcal{B} \subset \R$ and $c > 0$ such that $\alpha \pm \gamma \in (0,1)$ and, for every $t \in (0,T]$, $\VaR_{\alpha-\gamma}[g(X_t)], \VaR_{\alpha+\gamma}[g(X_t)] \in \mathcal{B}$ and the law of $g(X_t)$ admits a density bounded below by $c$ on $\mathcal{B}$. Then there exist $B, C > 0$, depending only on $d, D, T, g, \varepsilon, \lbold^*, \theta, \gamma$ and $c$, such that
    \[
        \big| \VaR_\alpha[g(X_t)] - \VaR_\alpha[g(X_t(M,N))] \big| \le C\, \mathcal{E}_{\theta a}(M,N), \qquad t \in [0,T],\ M \in \Nbar,\ N \in \Nbar^*.
    \]
\end{enumerate}
\end{theorem}

\begin{remark}\label{rem:law}
    Lemma~\ref{lem:cdf-bound} and Theorem~\ref{th:var-bound}(ii) are the only places where the law of the process enters, and only through the target $g(X_t)$: no assumption is made on the law of the approximations $g(X_t(M,N))$. Obtaining explicit lower bounds for densities is a genuine difficulty in its own right; we refer to \citet*{NourdinViens2009} for Malliavin-calculus techniques and to \citet*{Kohatsu2003, Bally2006} for lower bounds under ellipticity conditions. The constant $c$ enters only through a multiplicative factor of order $c^{-1}$, so a crude lower bound degrades the constant without affecting the rate. Lemma~\ref{lem:probability-bound} and Theorems~\ref{th:coherent-bound} and~\ref{th:var-bound}(i), which need no assumption on the law, also hold verbatim for the discretisation error at a fixed truncation level, between $X_t(M,\infty)$ and $X_t(M,N)$, with $\mathcal{E}_a(M,N)$ replaced by $N^{-a/2}$, uniformly in $M \in \Nbar$. In dimension $d = 1$ and at the grid points, $N^{-a/2}$ can be replaced by $N^{-a}$ in all the statements of this subsection.
\end{remark}

\subsection{Linear delay equations} \label{sec:linear-delay}

In Sections~\ref{sec:linear-delay} and \ref{sec:gaussian-volterra-processes}, errors are estimated by Monte Carlo over $20\,000$ paths of the time-extended Brownian motion, sampled on a uniform grid of $10^4$ steps. On this grid, the reference process $X$ is computed by an Euler scheme applied to its defining dynamics, independently of any signature. The approximations are computed from the signature of the piecewise linear interpolation of the same path: along the full grid for the truncated process, which is thus $X(M, 10^4)$, and along a coarser grid of $N$ steps for $X(M,N)$. The levels range over $M \in \intset{2}{7}$ and $N \in \{10, 50, 100, 500, 1\,000, 5\,000\}$. Strong errors are computed with $p = 3$, and the test function is $g = \mathrm{Id}$ unless stated otherwise. Since Theorem~\ref{th:Lp-truncation-bound} holds uniformly in $N$, it applies to $X(M,10^4)$ as stated; strictly speaking, the errors labelled ``truncation'' are total errors at $N = 10^4$, of order $\mathcal{E}_a(M,10^4)$ by Corollary~\ref{cor:total-error}.

\medskip

Consider the stochastic delay equation
\[
  \diff X_t = \sum_{i=0}^d \Big( a_i + b_i X_t + \int_0^t K_i(t-s)\, X_s\,\diff s \Big) \diff \Wh^i_t, \qquad X_0 = x, \qquad K_i(t) = \sum_{k=1}^{n_i} \alpha_i^k e^{\beta_i^k t},
\]
with real parameters. By \citet*[Section~4.2]{AbiJaberGerardHuang2024}, the process $X$ admits an infinite signature series representation $X_t = \textbrackets{\lbold^{\mathrm{DE}}}{\WH_t}$ for every $t \in [0,T]$, almost surely, where 
\[
    \lbold^{\mathrm{DE}} \coloneqq \boldsymbol{p}\,\big(\empt - \boldsymbol{q}\big)^{-1}, \qquad \boldsymbol{p} \coloneqq x\,\empt + \sum_{i=0}^d a_i\,\blold{i} - \frac12\Big(\sum_{i=1}^d a_i b_i\Big)\blold{0}, \qquad \boldsymbol{q} \coloneqq \sum_{i=0}^d \boldsymbol{c}_i\,\blold{i} - \frac12\Big(\sum_{i=1}^d b_i\,\boldsymbol{c}_i\Big)\blold{0},
\]
with $\boldsymbol{c}_i \coloneqq b_i\,\empt + \sum_{k=1}^{n_i} \alpha_i^k\,\blold{0}\,\exp_{\shuffle}(\beta_i^k\blold{0})$ and $(\empt - \boldsymbol{q})^{-1} \coloneqq \sum_{n \ge 0} \boldsymbol{q}^{\otimes n}$, a finite sum coordinatewise since $\boldsymbol{q}^{\empt} = 0$. Moreover, $\lbold^{\mathrm{DE}} \in \mathcal{A}_{\exp}$ \citep*[Section~4.2]{AbiJaberGerardHuang2024}, so that $\lbold^{\mathrm{DE}} \in \mathcal{L}^1$ (Remark~\ref{rem:Aexp}) and its coefficients grow at most geometrically, $|\lbold^{\mathrm{DE},\blold{v}}| \le C^{|\blold{v}|+1}$. Such a geometric growth only changes the constant $B$ in Theorem~\ref{th:Lp-truncation-bound}, which therefore holds with the optimal exponent $\varepsilon = 1/2$; this is proved in Remark~\ref{rem:geometric-coefficients}.

\medskip

\begin{figure}[p]
  \centering
  \captionsetup[sub]{skip=1pt, font=small}
  \captionsetup{skip=4pt}

  \begin{subfigure}[t]{0.49\textwidth}
    \centering
    \caption{Strong truncation error}
    \includegraphics[width=\linewidth]{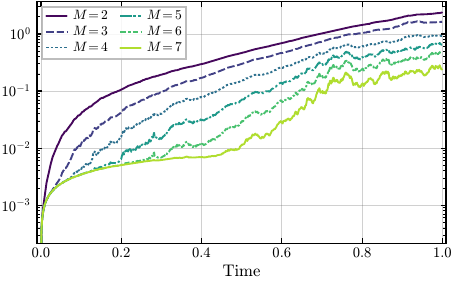}
  \end{subfigure}
  \hfill
  \begin{subfigure}[t]{0.49\textwidth}
    \centering
    \caption{Weak truncation error}
    \includegraphics[width=\linewidth]{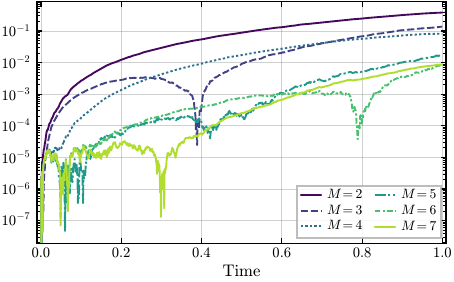}
  \end{subfigure}

  \vspace{0.2em}

  \begin{subfigure}[t]{0.49\textwidth}
    \centering
    \caption{Strong discretisation error, $M = 5$}
    \includegraphics[width=\linewidth]{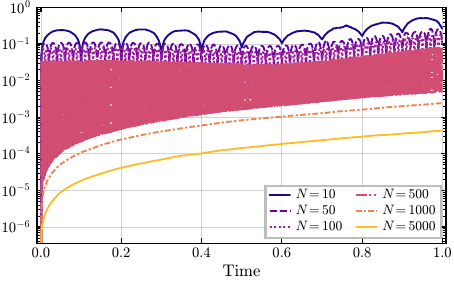}
  \end{subfigure}
  \hfill
  \begin{subfigure}[t]{0.49\textwidth}
    \centering
    \caption{Weak discretisation error, $M = 5$}
    \includegraphics[width=\linewidth]{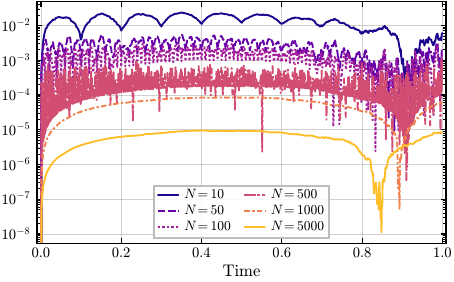}
  \end{subfigure}

  \caption{Strong (left) and weak (right) truncation (top) and discretisation (bottom, $M = 5$) errors as functions of time, for the delay equation of
  Section~\ref{sec:linear-delay}. Vertical scales are logarithmic; the downward cusps of the weak errors are sign changes of the Monte Carlo estimate of the mean.}
  \label{fig:linear-delay-errors-full}

  \vspace{1em}

  \begin{subfigure}[t]{0.49\textwidth}
    \centering
    \caption{Strong truncation error, fit $CB^M/(M!)^{1/2}$}
    \includegraphics[width=\linewidth]{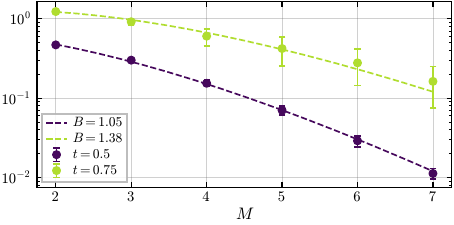}
  \end{subfigure}
  \hfill
  \begin{subfigure}[t]{0.49\textwidth}
    \centering
    \caption{Weak truncation error, fit $CB^M/(M!)^{1/2}$}
    \includegraphics[width=\linewidth]{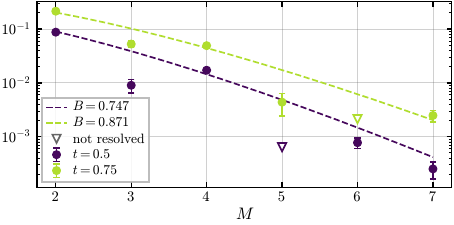}
  \end{subfigure}

  \vspace{0.2em}

  \begin{subfigure}[t]{0.49\textwidth}
    \centering
    \caption{Strong discretisation error multiplied by $\sqrt N$}
    \includegraphics[width=\linewidth]{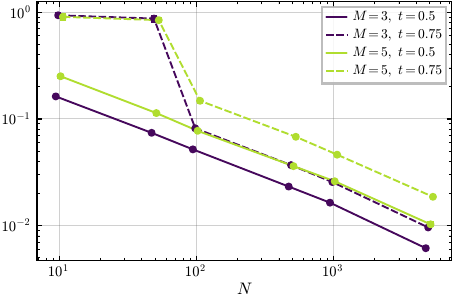}
  \end{subfigure}
  \hfill
  \begin{subfigure}[t]{0.49\textwidth}
    \centering
    \caption{Weak discretisation error multiplied by $N$}
    \includegraphics[width=\linewidth]{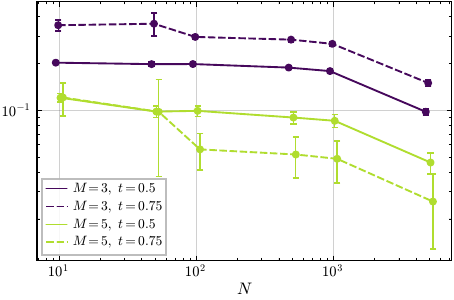}
  \end{subfigure}

  \caption{Dependence on $M$ (top) and $N$ (bottom) of the errors of
  Figure~\ref{fig:linear-delay-errors-full}, at $t = 0.5$ and $t = 0.75$, with $\pm 2$ standard errors; open triangles mark estimates not resolved at this level, drawn at the upper bound. Top: fit of $CB^M/(M!)^{1/2}$ (dashed). Bottom: errors multiplied by $\sqrt N$ (strong) and $N$ (weak), for $M \in \{3,5\}$; a bounded curve means that the exponent is attained.}
  \label{fig:linear-delay-errors-focus}
\end{figure}

We take $d = 1$, $T = 1$, $x = 0.25$, $a = (0.25, 1.0)$, $b = (0, 1.0)$, $n = (1,1)$, $\alpha^1 = (-4.0, -2.0)$ and $\beta^1 = (1.0, 1.5)$, each vector being indexed by $i \in \{0,1\}$.

\emph{Truncation.} The strong errors are ordered in $M$ over the whole horizon and grow with $t$ (Figure~\ref{fig:linear-delay-errors-full}(a)). Over $M \le 7$, a factorial decay is hard to tell apart from a geometric one by eye, but the two-parameter law $CB^M/(M!)^{1/2}$ of Theorem~\ref{th:Lp-truncation-bound} fits the data within the Monte Carlo error at both dates (Figure~\ref{fig:linear-delay-errors-focus}(a)). The fitted values $B \approx 1.05$ and $B \approx 1.38$ have ratio $0.76$, of the order of the value $\sqrt{0.5/0.75} \approx 0.82$ suggested by the time dependence of the truncation bound: at time $t$, the proof of Theorem~\ref{th:Lp-truncation-bound} gives a bound proportional to $(B\sqrt t)^{M+1}/((M+1)!)^{\varepsilon}$ (Lemma~\ref{lem:truncation}). The weak errors are smaller but not monotone in $M$, as the contributions of successive levels partially cancel (Figure~\ref{fig:linear-delay-errors-focus}(b)). At $t = 0.5$, the weak error for $M = 3$ lies just after a sign change of $\E[X_t - X_t(M)]$, visible as a cusp near $t \approx 0.4$ in Figure~\ref{fig:linear-delay-errors-full}(b), which explains why it falls well below the fitted envelope; for $M = 5$, the mean is too close to $0$ to be resolved by the Monte Carlo estimate. The gain lies in the constant, not in the rate, in agreement with Corollary~\ref{cor:weak-truncation-bound} and Proposition~\ref{prop:truncation-sharpness}.

\emph{Discretisation.} Each curve of Figure~\ref{fig:linear-delay-errors-full}(c)--(d) is scalloped: the interpolation coincides with the path at the grid points, so the level-one error vanishes there and the total error drops to the smaller order $N^{-1}$ predicted by Theorem~\ref{th:Lp-discretisation-bound} in dimension $d = 1$ (Proposition~\ref{prop:Lp-discretisation-d1}), while it peaks, at order $N^{-1/2}$, inside each step. The upper envelopes decrease from about $3\cdot10^{-1}$ at $N = 10$ to $2\cdot10^{-3}$ at $N = 5\,000$ in the strong case, and from $2\cdot10^{-2}$ to $3\cdot10^{-5}$ in the weak case. Rescaled by $N$, the weak errors are essentially flat (Figure~\ref{fig:linear-delay-errors-focus}(d)): the rate $N^{-1}$ of Theorem~\ref{th:weak-discretisation-bound} is attained. Rescaled by $\sqrt N$, the strong errors show two regimes (Figure~\ref{fig:linear-delay-errors-focus}(c)). For $N \in \{10, 50\}$, $t = 0.75$ lies inside a step and the rescaled error is flat: this is the rate $N^{-1/2}$, optimal off the grid (Remark~\ref{rem:optimality-statement}(i)). For $N \ge 100$, $t = 0.75$ is a grid point, as $t = 0.5$ is for every $N$, and the error decays like $N^{-1}$, as predicted by Theorem~\ref{th:Lp-discretisation-bound} in dimension $d = 1$. The curves for $M = 3$ and $M = 5$ are of the same order, in agreement with the uniformity in $M$ of these estimates.

\subsection{Gaussian Volterra processes} \label{sec:gaussian-volterra-processes}

For delay equations, the coefficients lie in $\mathcal{L}^1$ by construction. For Gaussian Volterra processes, membership depends on the kernel and must be checked case by case. By \citet*[Theorem~4.5]{AbiJaberGerardHuang2024}, a Gaussian Volterra process $X_t = \int_0^t K(t-s)\,\diff W_s$, with $K$ as in \citet*[Section~4.3]{AbiJaberGerardHuang2024}, admits the time-dependent representation
\[
    X_t = \brackets{\lbold^{\mathrm{GV}}_t}{\WH_t}, \qquad \lbold^{\mathrm{GV}}_t = \indic{t>0} \sum_{n \ge 0} K^{(n)}(t)\,(-1)^n\, \blold{0}^{\otimes n}\blold{1}.
\]
Since $\norm{\blold{0}^{\otimes n}\blold{1}} = 2n+1$, our estimates apply as soon as, for some $\varepsilon \in (0,1/2]$,
\begin{equation}\label{eq:kernel-condition}
    \sup_{t\in[0,T]} |K^{(n)}(t)| \lesssim \big[(2n+1)!\big]^{\frac12-\varepsilon}, \qquad n \in \N.
\end{equation}

\paragraph{Shifted Riemann--Liouville fractional Brownian motion}
For $H \in (0,1)$ and $\mu > 0$, the process $X_t = \Gamma(H+\frac12)^{-1}\int_0^t (t+\mu-s)^{H-\frac12}\,\diff W_s$ has, by
\citet*[Section~4.4.1]{AbiJaberGerardHuang2024}, the coefficients
\[
    \lbold^{\mathrm{RL}}_t = \sum_{n \ge 0} c_n(t)\, \blold{0}^{\otimes n}\blold{1}, \qquad c_n(t) \coloneqq \frac{(t+\mu)^{H-\frac12}}{\Gamma(H+\frac12)}\, \frac{(\frac12-H)^{\nup}}{(t+\mu)^n},
\]
where $x^{\nup} \coloneqq \prod_{k=0}^{n-1}(x+k)$ is the rising factorial. For $H \ne 1/2$, the asymptotics $\Gamma(n+x) \sim n^x \Gamma(n)$ give $|c_n(t)| \sim C_{t}\, n!\, n^{-\frac12-H} (t+\mu)^{-n}$, so that the ratio $|c_n(t)|/[(2n+1)!]^{\frac12-\varepsilon}$ grows like $(n^{2\varepsilon}/C'_t)^n$ and \eqref{eq:kernel-condition} fails for every $\varepsilon > 0$, whereas $|c_n(t)| \lesssim (2\mu)^{-n}\sqrt{(2n+1)!}$ uniformly in $t$. These coefficients therefore lie on the boundary $\varepsilon = 0$ of $\mathcal{L}^1$. There, the factorial gain is lost, and the moment argument behind Theorem~\ref{th:Lp-truncation-bound} only yields a geometric rate on a short horizon (Remark~\ref{rem:boundary-case}, with $c_1 = 1$ and $c_2 = (2\mu)^{-1} \vee 1$). The restriction comes from the method, not from the process: the generic bound accounts for every word of a given length, through the factor $(p+2)^{|\blold{v}|}$ of Proposition~\ref{prop:signature-moment-bound}, whereas $\lbold^{\mathrm{RL}}$ charges a single word per length, whose signature component is Gaussian.

\begin{proposition}[Riemann–-Liouville truncation error]\label{prop:RL}
    Let $d = D = 1$, $H \in (0,1)$, $\mu > 0$, $q_t \coloneqq t/(t+\mu)$, and $\kappa_p \coloneqq \E[|Z|^p]^{1/p}$ for $Z \sim \mathcal N(0,1)$. For every $N \in \Nbar^*$ and $t \in (0,T]$, the series $X_t(\infty,N) \coloneqq \textbrackets{\lbold^{\mathrm{RL}}_t}{\WT_t(N)}$ converges absolutely almost surely and in every $L^p$, and, for every $p \ge 1$ and $M \in \N^*$:
    \begin{enumerate}[label=(\roman*)]
        \item uniformly in $N \in \Nbar^*$, $\displaystyle \textnorm{X_t(\infty,N) - X_t(M,N)}_p \le \frac{\kappa_p (t+\mu)^{H+\frac12}\, t^{\frac12}}{\mu\,\Gamma(H+\frac12)}\, \frac{q_t^M}{\sqrt{2M+1}}$;
        \item if $H \ne 1/2$, $\displaystyle \textnorm{X_t - X_t(M)}_p \ge \frac{\kappa_p\, |\frac12-H|\, (t+\mu)^{H-\frac12}\, t^{\frac12}}{2\,\Gamma(H+\frac12)}\, \frac{q_t^M}{M^{3/2}\sqrt{2M+1}}$.
    \end{enumerate}
    In particular, $\log\norm{X_t - X_t(M)}_p = M\log q_t + O(\log M)$: the truncation error decays geometrically, with ratio $q_t < 1$, on the whole horizon.
\end{proposition}

\begin{proof}
    Write $a \coloneqq \frac12 - H$, $f_n(s) \coloneqq s^n/n!$ and $V_n(t) \coloneqq \int_0^t f_n^2 = t^{2n+1}/((n!)^2(2n+1))$. The truncation keeps the words of length $n+1 \le M$, so that $R_t \coloneqq X_t(\infty,N) - X_t(M,N) = \sum_{n\ge M} c_n(t)\, \WT^{\blold{0}^{\otimes n}\blold{1}}_t(N)$. The time coordinate is interpolated exactly, so $\WT^{\blold{0}^{\otimes n}}_s(N) = f_n(s)$ for every $N$. Hence $\WH^{\blold{0}^{\otimes n}\blold{1}}_t = \int_0^t f_n\,\diff W$ and, for $N < \infty$, with $I_k \coloneqq (t_{k-1},t_k]\cap(0,t]$, $\WT^{\blold{0}^{\otimes n}\blold{1}}_t(N) = \sum_k h^{-1}(W_{t_k} - W_{t_{k-1}})\int_{I_k} f_n$. Both are centred Gaussian with variance at most $V_n(t)$, by the Cauchy--Schwarz inequality $h^{-1}(\int_{I_k} f_n)^2 \le \int_{I_k} f_n^2$ in the second case.

    \emph{(i)} Since $|a+k| \le k+1$, $|a^{\nup}| \le n!$, so that $|c_n(t)|\sqrt{V_n(t)} \le (t+\mu)^{H-\frac12} t^{\frac12}\, q_t^n /(\Gamma(H+\frac12)\sqrt{2n+1})$. The sum of these bounds over $n \ge 0$ is finite, which gives the absolute convergence, almost surely and in $L^2$. Hence $R_t$ is centred Gaussian and $\norm{R_t}_p = \kappa_p\norm{R_t}_2$, and summing over $n \ge M$ with $\sum_{n\ge M} q_t^n = q_t^M (t+\mu)/\mu$ gives (i).

    \emph{(ii)} For $N = \infty$, all $c_n(t)$, $n \ge 1$, have the sign of $a$, and $f_n \ge 0$. By the Itō isometry, $\norm{R_t}_2^2 = \int_0^t (\sum_{n\ge M} c_n(t) f_n)^2 \ge c_M(t)^2 V_M(t)$. Finally, $|a^{\overline{M}}| \ge |a|(M-1)!\prod_{k=1}^{M-1}(1-\frac1{2k}) \ge |a|(M-1)!/(2\sqrt{M})$, since $\prod_{k=1}^m(1-\frac1{2k}) = \binom{2m}{m}4^{-m} \ge (2\sqrt{m})^{-1}$, and $\norm{R_t}_p = \kappa_p\norm{R_t}_2$ concludes.
\end{proof}

Other families on this boundary behave differently. $\exp(c^2 W_t^2/2) = \textbrackets{\sum_{k \ge 0} c^{2k}(2k-1)!!\,\blold{1}^{\otimes 2k}}{\WH_t}$, for instance, is not integrable for $t \ge 1/c^2$, so that its truncation error diverges beyond this horizon (Example~\ref{ex:boundary-case}). The growth of the coefficients alone therefore does not decide the boundary case.

\begin{figure}[H]
  \centering

  \begin{subfigure}[b]{0.49\textwidth}
    \centering
    \caption{Strong truncation error}
    \includegraphics[width=\linewidth]{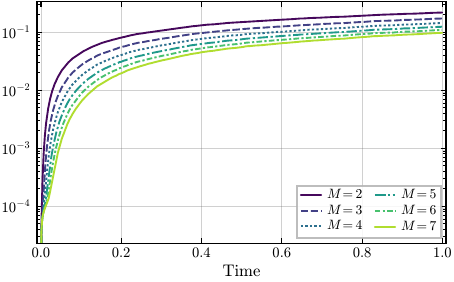}
  \end{subfigure}
  \hfill
  \begin{subfigure}[b]{0.49\textwidth}
    \centering
    \caption{Weak truncation error}
    \includegraphics[width=\linewidth]{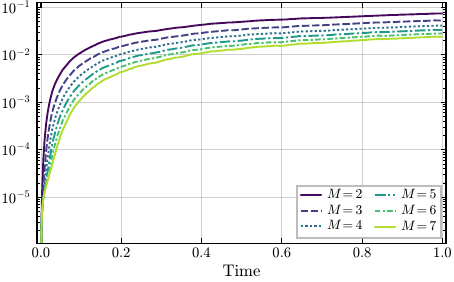}
  \end{subfigure}

  \caption{Strong (left) and weak (right) truncation errors as functions of time, for the shifted Riemann--Liouville process with $H = 0.3$ and $\mu = 1/52$; the weak errors are computed with $g = \cos$.} \label{fig:fBm-errors}
\end{figure}

For the numerical illustration, we take $\mu = 1/52$ and $T = 1$. Since every
$\WH^{\blold{0}^{\otimes n}\blold{1}}_t$ is centred, the weak error vanishes for $g = \mathrm{Id}$, and we use $g = \cos$, which is Lipschitz. The levels are ordered in $M$ over the whole horizon (Figure~\ref{fig:fBm-errors}), but the decay in $M$ is far slower than for the delay equation, as Proposition~\ref{prop:RL} predicts: at $t = 1$, $q_t = 52/53 \approx 0.98$, so that over $M \in \intset{2}{7}$ the geometric factor contributes only about $0.91$, and the separation between levels comes mostly from the factor $(2M+1)^{-1/2}$. Near $t = 0$, $q_t$ is much smaller and the levels separate quickly. The case $H = 0.7$ gives the same picture, with errors about half as large.

\subsection{Pricing with a signature volatility model} \label{sec:pricing}

Following \citet*{AbiJaberGerard2025}, let $d = 1$ and consider the stock price
\[
  \diff S_t = S_t\,\sigma_t\,\diff W^S_t, \qquad \sigma_t = \brackets{\lbold}{\WH_t}, \qquad t \in [0,T],
\]
with $S_0 > 0$, $\lbold \in \mathcal{L}^1(\R^2)$, zero interest rates, and
$W^S = \rho W + \sqrt{1-\rho^2}\,W^\perp$ for a Brownian motion $W^\perp$ independent of $W$. We price the put with payoff $\Pi(S) \coloneqq (K - S)^+$. The price $S$ itself need not be a linear functional of a signature, but its logarithm is: by Itō's formula and the shuffle identity, $X_t \coloneqq \log S_t = \textbrackets{\Phi(\lbold)}{\WH^*_t}$, where $\WH^*$ is the signature of $(t, W_t, W^\perp_t)$ and $\Phi$, an explicit combination of shuffle products and concatenations with letters given in Section~\ref{sec:pricing-proof}, maps $\mathcal{L}^1(\R^3)$ into itself. Extending $\lbold$ by zero on words containing the letter $2$, we view $\lbold \in \mathcal{L}^1(\R^3)$, and $\sigma_t = \textbrackets{\lbold}{\WH^*_t}$. Moreover, none of the operations in $\Phi$ shortens a word, so that $\Phi(\lbold)^{\le M} = \Phi(\lbold^{\le M})^{\le M}$ (Lemma~\ref{lem:phi-trunc}). The scheme
\[
  \check{X}_t(M,N) \coloneqq \brackets{\Phi(\lbold^{\le M})^{\le M}}{\WT^*_t(N)}
\]
therefore only requires $\lbold^{\le M}$, and coincides with $X_t(M,N)$ for the coefficient family $\Phi(\lbold)$. Let $\varepsilon$ be a summability exponent of $\lbold$ and $\varepsilon' < \varepsilon$. Since $\Phi(\lbold)$ admits $\varepsilon'$ as a summability exponent and $\Pi \circ \exp$ is Lipschitz, Corollary~\ref{cor:total-error} applied to $\Phi(\lbold)$ gives (Proposition~\ref{prop:pricing})
\[
    \Big| \E\, \Pi(S_t) - \E\, \Pi\big(e^{\check X_t(M,N)}\big) \Big| \le C\big( r'_M + N^{-1/2} \big), \qquad r'_M \coloneqq \frac{B^M}{(M!)^{\varepsilon'}}, \qquad t \in [0,T],\ M \in \N,\ N \in \N^*,
\]
with $B, C$ independent of $t$, $M$ and $N$. As after Corollary~\ref{cor:total-error}, $M$ and $N$ can thus be chosen independently.

\section{Proofs of the truncation and discretisation estimates} \label{sec:theorem}

This section proves the estimates of Section~\ref{sec:main-results}. We first control the truncation error in $L^p$ (Section~\ref{sec:truncation}), then the discretisation error, in $L^p$ (Section~\ref{sec:strong-discretisation}) and in the weak sense (Section~\ref{sec:weak-discretisation}). The $L^p$ estimates are also the starting point of the extensions of Section~\ref{sec:extensions}.

\subsection{Truncation error} \label{sec:truncation}

The truncation error relies almost entirely on the next lemma.

\begin{lemma}[Tail of a linear functional] \label{lem:truncation}
    Let $\lbold : \N \times [0, T] \to T((\R^{d+1}))$, $q \in \N^*$ and $p \geq 1$. Define
    \[
        \Delta^q_p(\lbold, t, M, N) \coloneqq \norm{\brackets{\lbold(M, t)}{\WT^{>M}_t(N)}^q}_p, \qquad t \in [0, T], \ M \in \N, \ N \in \Nbar^*.
    \]
    Assume that there exists $\lbold^\ast \in \mathcal{L}^1$ such that $\lbold^{\ast,\blold{v}} \le C_\ast \varphi_\varepsilon(\blold{v})$ for every $\blold{v} \in V$ for some $\varepsilon \in (0, 1/2]$ and $C_\ast > 0$, and $\lbold(M, t) \preceq \lbold^\ast$ for every $M \in \N$ and $t \in [0, T]$. Set
    \[
        B_r \coloneqq 2(r+2)\sqrt{c_d}\,\big(d + 2\sqrt{c_d T}\big), \quad r \ge 1, \qquad E_\varepsilon(y) \coloneqq \sum_{j \ge 0} \frac{y^j}{(j!)^\varepsilon}, \quad y \ge 0.
    \]
    Then, for every $t \in [0,T]$, $M \in \N$ and $N \in \Nbar^*$,
    \[
        \Delta^q_p(\lbold, t, M, N) \le \left( C_*\, E_\varepsilon\big(B_{pq}\sqrt T\big)\, \frac{\big(B_{pq}\sqrt t\big)^{M+1}}{\big((M+1)!\big)^{\varepsilon}} \right)^{q}.
    \]
    In particular, $\Delta^q_p(\lbold,t,M,N) \lesssim t^{q\frac{M}{2}} B_{pq}^{qM} (M!)^{-q\varepsilon}$, which decays like $\exp\{-q\varepsilon M\log M + O(M)\}$ uniformly in $N$.
\end{lemma}
 
\begin{proof}
    Fix $N \in \Nbar^*$ and write $\X \coloneqq \WT(N)$. Set $r \coloneqq pq$. Since $\Delta^q_p = (\Delta^1_{r})^q$, it suffices to bound $\Delta \coloneqq \Delta^1_r(\lbold,t,M,N)$. For a word $\blold{v}$, write $m \coloneqq |\blold{v}|$, $n_0 \coloneqq n_0(\blold{v})$ and $n \coloneqq \norm{\blold{v}} = m + n_0$. By Proposition~\ref{prop:signature-moment-bound} and Inequality~\eqref{eq:psi-bound}, and using $n/2 + 1 \le 2^n$,
    \[
        \lbold^{\ast,\blold{v}}\, \norm{\X^{\blold{v}}_t}_r \le C_*\, (n!)^{\frac12-\varepsilon}\, (r+2)^m (2c_d t)^{n/2}\sqrt{\frac{n/2+1}{n!}} \le C_*\, (r+2)^m\, \big(4 c_d t\big)^{n/2}\, (m!)^{-\varepsilon},
    \]
    where the last step also uses $n \ge m$. By Lemma~\ref{lem:count-length-m-k-zeros}, there are $\binom{m}{n_0} d^{\,m-n_0}$ words of length $m$ with $n_0$ occurrences of $0$, so the binomial theorem gives
    \begin{equation}\label{eq:per-length}
        \sum_{|\blold{v}| = m} \lbold^{\ast,\blold{v}}\, \norm{\X^{\blold{v}}_t}_r \le C_*\, \frac{x_t^m}{(m!)^{\varepsilon}}, \qquad x_t \coloneqq 2(r+2)\sqrt{c_d t}\,\big(d + 2\sqrt{c_d t}\big) \le B_r \sqrt{t}.
    \end{equation}
    In particular $\sum_{\blold{v}} \lbold^{\ast,\blold{v}} \norm{\X^{\blold{v}}_t}_r < \infty$, so the series $\brackets{\lbold(M,t)}{\X^{>M}_t}$ converges absolutely almost surely and in $L^r$, and Minkowski's inequality together with $\lbold(M,t) \preceq \lbold^\ast$ yields
    \[
        \Delta \le \sum_{|\blold{v}| > M} \lbold^{\ast,\blold{v}}\, \norm{\X^{\blold{v}}_t}_r \le C_* \sum_{m \ge M+1} \frac{x_t^m}{(m!)^\varepsilon}.
    \]
    Writing $m = M+1+j$ and using $(M+1+j)! \ge (M+1)!\, j!$,
    \[
        \Delta \le C_*\, \frac{x_t^{M+1}}{\big((M+1)!\big)^\varepsilon} \sum_{j\ge0} \frac{x_t^j}{(j!)^\varepsilon} \le C_*\, E_\varepsilon\big(B_r\sqrt T\big)\, \frac{\big(B_r\sqrt t\big)^{M+1}}{\big((M+1)!\big)^\varepsilon},
    \]
    the series $E_\varepsilon$ being finite for every $\varepsilon > 0$. Raising to the power $q$ concludes the proof.
\end{proof}

\begin{corollary}[Well-posedness on $\mathcal{L}^1$]\label{cor:well-posed}
    Let $\lbold \in \mathcal L^1_T$ be dominated by $\lbold^\ast \in \mathcal{L}^1$ with $\lbold^{\ast,\blold{v}} \le C_* \varphi_\varepsilon(\blold{v})$, and $p \ge 1$. Then
    \[
        K_p \coloneqq \sup_{N \in \Nbar^*} \ \sup_{t \in [0,T]} \ \sum_{\blold{v} \in V} \lbold^{\ast,\blold{v}} \, \textnorm{\WT^{\blold{v}}_t(N)}_p \le C_*\, E_\varepsilon\big(B_p \sqrt T\big) < \infty.
    \]
    Consequently, for every $t \in [0,T]$ and $N \in \Nbar^*$:
    \begin{enumerate}[label=(\roman*)]
        \item the series $\textbrackets{\lbold_t}{\WT_t(N)} = \sum_{\blold{v}} \lbold^{\blold{v}}_t \, \WT^{\blold{v}}_t(N)$ converges absolutely almost surely, so its value does not depend on the order of summation, and it converges in $L^p$;
        \item $\sup_{M \in \Nbar} \ \sup_{N} \ \sup_{t} \textnorm{\textbrackets{\lbold_t}{\WT^{\le M}_t(N)}}_p \le K_p$;
        \item there exists $B > 0$, depending only on $d$, $T$ and $p$, such that $\textnorm{\textbrackets{\lbold_t}{\WT^{> M}_t(N)}}_p \lesssim r_M$, uniformly in $t \in [0,T]$ and $N \in \Nbar^*$.
    \end{enumerate}
\end{corollary}

\begin{proof}
    The bound on $K_p$ is \eqref{eq:per-length} with $r = p$, summed over $m \ge 0$. Since $\E\big[\sum_{\blold{v}} |\lbold^{\blold{v}}_t| \, |\WT^{\blold{v}}_t(N)|\big] \le K_1 < \infty$, the series converges absolutely almost surely, and by monotone convergence and Minkowski's inequality, $\textnorm{\sum_{\blold{v}} |\lbold^{\blold{v}}_t\, \WT^{\blold{v}}_t(N)|}_p \le K_p$. This gives (i) and (ii), the bound on any sub-sum being at most the bound on the full one. Point (iii) is Lemma~\ref{lem:truncation} with $q = 1$.
\end{proof}

\begin{remark}[Admissibility of $\mathcal{L}^1$]\label{rem:L1-admissible}
    Membership of $\mathcal{A}(\WH)$ requires a single null set for all $t \in [0,T]$, whereas Corollary~\ref{cor:well-posed} provides one for each fixed $t$. The supremum bound of \citet*[Proposition~3.5]{AbiJaberGerardHuang2024} closes the gap: it gives $\textnorm{ \sup_{t \le T} |\WH^{\blold{v}}_t| }_2 \lesssim B^{\norm{\blold{v}}}/\sqrt{(\norm{\blold{v}} - 1)!}$ for some $B > 0$, which has the same decay as the bound~\eqref{eq:psi-bound}, for $\blold{v} \ne \empt$, while $\WH^{\empt} \equiv 1$. For $\lbold \in \mathcal{L}^1$ with summability exponent $\varepsilon$, the summation of Lemma~\ref{lem:truncation} therefore yields $\E\big[\sum_{\blold{v}} |\lbold^{\blold{v}}| \sup_{t\le T} |\WH^{\blold{v}}_t|\big] < \infty$. Hence, almost surely, $\norm{\lbold}^{\mathcal{A}(\WH)}_t \le \sum_{\blold{v}} |\lbold^{\blold{v}}|\,|\WH^{\blold{v}}_t| < \infty$ for all $t \in [0,T]$ simultaneously, that is, $\mathcal{L}^1 \subset \mathcal{A}(\WH)$. In particular, Proposition~\ref{prop:shuffle-identity} applies to $\lbold \in \mathcal{L}^1$.
\end{remark}

\begin{lemma}[Shuffle identity on $\mathcal{L}^1$]\label{lem:shuffle-L1}
Let $\lbold_1, \lbold_2 \in \mathcal{L}^1$ and $N \in \Nbar^*$. Then $\lbold_1 \shuffle \lbold_2 \in \mathcal{L}^1$ and, for every $t \in [0,T]$, almost surely,
\[
  \brackets{\lbold_1}{\WT_t(N)} \brackets{\lbold_2}{\WT_t(N)} = \brackets{\lbold_1 \shuffle \lbold_2}{\WT_t(N)}.
\]
\end{lemma}

\begin{proof}
    By Proposition~\ref{prop:L1-stability}, applied to $\lbold_1, \lbold_2$ and to $|\lbold_1|, |\lbold_2|$, both $\lbold_1 \shuffle \lbold_2$ and $\lbold^* \coloneqq |\lbold_1| \shuffle |\lbold_2|$ belong to $\mathcal{L}^1$. Write $\X \coloneqq \WT(N)$ and consider the family, indexed by triples of words,
    \[
        a_{\blold{v},\blold{w},\blold{u}} \coloneqq \lbold_1^{\blold{v}}\, \lbold_2^{\blold{w}}\, \brackets{\blold{u}}{\blold{v} \shuffle \blold{w}} \, \X^{\blold{u}}_t.
    \]
    For each $\blold{u}$, only finitely many pairs $(\blold{v}, \blold{w})$ contribute, and the coefficients $\brackets{\blold{u}}{\blold{v} \shuffle \blold{w}}$ are non-negative, so that $\sum_{\blold{v},\blold{w},\blold{u}} |a_{\blold{v},\blold{w},\blold{u}}| = \sum_{\blold{u}} \lbold^{*,\blold{u}}\, |\X^{\blold{u}}_t|$, which is finite almost surely by Corollary~\ref{cor:well-posed}(i) applied to $\lbold^*$. On this event, the family may be summed in any order. Summing first over $(\blold{v}, \blold{w})$ gives $\brackets{\lbold_1 \shuffle \lbold_2}{\X_t}$. Summing first over $\blold{u}$ and using the shuffle identity for words (Remark~\ref{rem:shuffle-BV} for $N < \infty$, Proposition~\ref{prop:shuffle-identity} for $N = \infty$) gives $\sum_{\blold{v},\blold{w}} \lbold_1^{\blold{v}}\lbold_2^{\blold{w}}\,\X^{\blold{v}}_t \X^{\blold{w}}_t$, which is the product of the two absolutely convergent series $\brackets{\lbold_1}{\X_t}$ and $\brackets{\lbold_2}{\X_t}$.
\end{proof}

\begin{remark}[Geometric growth of the coefficients] \label{rem:geometric-coefficients}
    The assumption of Lemma~\ref{lem:truncation} may be relaxed to
    \[
        \lbold^{\ast,\blold{v}} \le C_*\, c_1^{|\blold{v}|} c_2^{n_0(\blold{v})}\, \varphi_\varepsilon(\blold{v}), \qquad \blold{v} \in V,
    \]
    for some $c_1 > 0$, $c_2 \ge 1$ and $\varepsilon \in (0, 1/2]$. Indeed, the extra factor only modifies the per-length bound~\eqref{eq:per-length}: the binomial sum over $n_0$ now carries the weight $\big(c_2\sqrt{4c_d t}\big)^{n_0}$ and the common factor $c_1^m$, so that $x_t$ is replaced by
    \[
        2c_1(r+2)\sqrt{c_d t}\,\big(d + 2c_2\sqrt{c_d t}\big) \le c_1 c_2\, B_r \sqrt t ,
    \]
    and the rest of the proof is unchanged. The conclusion of Lemma~\ref{lem:truncation} therefore holds with $B_r$ replaced by $c_1c_2B_r$, uniformly in $N \in \Nbar^*$: the factorial gain $\big((M+1)!\big)^{-\varepsilon}$ is kept, and only the geometric constant deteriorates. For $\varepsilon < 1/2$, such coefficients already lie in $\mathcal{L}^1$ with any exponent $\varepsilon' < \varepsilon$: since $|\blold{v}| \le \norm{\blold{v}}$ and $n_0(\blold{v}) \le \norm{\blold{v}}$, $c_1^{|\blold{v}|} c_2^{n_0(\blold{v})} \le ((c_1 \vee 1)c_2)^{\norm{\blold{v}}} \lesssim (\norm{\blold{v}}!)^{\varepsilon-\varepsilon'}$. The remark keeps the exponent $\varepsilon$ itself, which matters most at the endpoint $\varepsilon = 1/2$, where $\varphi_{1/2} \equiv 1$. By Remark~\ref{rem:Aexp}, elements of $\mathcal A_{\exp}$ satisfy $|\lbold^{\blold{v}}| \le C^{|\blold{v}|+1}$, which is not $\lesssim \varphi_{1/2}(\blold{v})$ when $C > 1$, but is of the above form with $\varepsilon = 1/2$, $C_* = C$, $c_1 = C$ and $c_2 = 1$. Hence, for $\lbold \in \mathcal A_{\exp}$, Theorem~\ref{th:Lp-truncation-bound} holds with the optimal exponent $\varepsilon = 1/2$, for every $N \in \Nbar^*$, the constant $B$ being multiplied by $C$. This is the value used in the numerical fits of Section~\ref{sec:linear-delay}.
\end{remark}

The coefficients enter the rate of Lemma~\ref{lem:truncation} only through $\varepsilon$, and through $c_1c_2$ in the setting of Remark~\ref{rem:geometric-coefficients}. The extreme regime $\varepsilon = 0$ is worth recording, as it costs the horizon.

\begin{remark}[The boundary of $\mathcal L^1$]\label{rem:boundary-case}
    Suppose only that $\lbold^{\ast,\blold{v}} \le C_*\, c_1^{|\blold{v}|} c_2^{n_0(\blold{v})}\varphi_0(\blold{v})$, $\blold{v} \in V$, for some $c_1, c_2 \ge 1$, where $\varphi_0(\blold{v}) = \sqrt{\norm{\blold{v}}!}$. With $m \coloneqq |\blold{v}|$ and $n \coloneqq \norm{\blold{v}} \le 2m$, Proposition~\ref{prop:signature-moment-bound} and Inequality~\eqref{eq:psi-bound} give
    \[
        \lbold^{\ast,\blold{v}}\,\textnorm{\WT^{\blold{v}}_t(N)}_r \le C_*\, \sqrt{m+1}\; \big(c_1(r+2)\big)^{m}\, c_2^{\,n_0(\blold{v})}\,(2c_d t)^{n/2},
    \]
    so the per-length bound~\eqref{eq:per-length} is replaced by $C_*\sqrt{m+1}\, y_t^m$, with $y_t \coloneqq c_1(r+2)\sqrt{2c_d t}\,\big(d + c_2\sqrt{2c_d t}\big)$. The factorial gain disappears, and the tail series converges only if $y_T < 1$. In that case, using $\sqrt{M+2+j} \le \sqrt{M+2}\,\sqrt{j+1}$ and $\sum_{j\ge0}(j+1)y^j = (1-y)^{-2}$,
    \[
        \Delta^1_r(\lbold, t, M, N) \le \frac{C_*\sqrt{M+2}}{(1-y_T)^2}\; y_t^{M+1}, \qquad t\in[0,T], \ M \in \N, \ N \in \Nbar^*.
    \]
    The estimates of Section~\ref{sec:main-results} therefore survive, with a geometric rate $y^M$ for any $y \in (y_T, 1)$ and uniformly in $N$, on any horizon $T_0$ such that $c_1(r+2)\sqrt{2c_d T_0}\,\big(d + c_2\sqrt{2c_d T_0}\big) < 1$, and only there. This horizon shrinks as $r = pq$ grows. See Section~\ref{sec:gaussian-volterra-processes} and Example~\ref{ex:boundary-case} for two examples.

    Summing the same per-length bound over all $m \ge 0$ gives $\sum_{\blold{v}} \lbold^{\ast,\blold{v}}\,\textnorm{\WT^{\blold{v}}_t(N)}_r \le C_*(1-y_T)^{-2}$ for every $t \in [0,T_0]$ and $N \in \Nbar^*$. As in Corollary~\ref{cor:well-posed}, the series $\textbrackets{\lbold_t}{\WT_t(N)}$ therefore converges absolutely almost surely and in $L^r$ on $[0,T_0]$, uniformly in $N$, and its truncations are bounded in $L^r$ uniformly in $M$, $N$ and $t$; this replaces Corollary~\ref{cor:well-posed}, which does not apply outside $\mathcal L^1$. Since $y_T$ increases with $r$, almost sure convergence ($r = 1$) holds on the largest horizon, and the $L^r$ statements on nested smaller ones. The horizon $T_0$ must be computed with the exponent $r$ used in the proof.
\end{remark}

\begin{example}[A boundary family with finite horizon]\label{ex:boundary-case}
    Let $d = 1$, $c > 0$, and $\lbold^{(c)} \coloneqq \sum_{k\ge0} c^{2k}(2k-1)!!\,\blold{1}^{\otimes 2k}$. Since $\WH^{\blold{1}^{\otimes 2k}}_t = W_t^{2k}/(2k)!$ and $(2k-1)!!/(2k)! = 1/(2^k k!)$, 
    \[
        X_t \coloneqq \brackets{\lbold^{(c)}}{\WH_t} = \exp\big(c^2 W_t^2/2\big).
    \]
    As $((2k-1)!!)^2 \le (2k-1)!!\,(2k)!! = (2k)!$, the coefficients are bounded by $c^{|\blold{v}|}\varphi_0(\blold{v})$, which is the setting of Remark~\ref{rem:boundary-case}. Using the bounds
    \begin{equation}\label{eq:central-binomial}
        \frac{1}{2\sqrt k} \le \binom{2k}{k} 4^{-k} \le \frac{1}{\sqrt{\pi k}}, \qquad k \ge 1,
    \end{equation}
    we get $(2k-1)!!/\sqrt{(2k)!} = (\binom{2k}{k}4^{-k})^{1/2} \ge (4k)^{-1/4}$, so the coefficients are not $\lesssim \varphi_\varepsilon$ for any $\varepsilon > 0$, and $\lbold^{(c)} \notin \mathcal{L}^1$. The restriction to a short horizon in Remark~\ref{rem:boundary-case} cannot be removed: since $\E[X_t^p] = (1 - pc^2t)^{-1/2}$ if $pc^2t < 1$ and $+\infty$ otherwise, $X_t$ belongs to $L^p$ exactly for $t < 1/(pc^2)$, a horizon that shrinks as $p$ grows. On this horizon, the truncation error converges at a geometric rate, without factorial gain. All terms being non-negative, and $\E[W_t^{2k}] = t^k(2k-1)!!$, $\E\big[X_t - X_t(M)\big] = \sum_{k \ge K} \binom{2k}{k} 4^{-k} (c^2 t)^k$, with $K \coloneqq \lfloor M/2 \rfloor + 1$. If $c^2 t \ge 1$, the terms are bounded below by $1/(2\sqrt k)$ by~\eqref{eq:central-binomial}, and the series diverges. If $c^2 t < 1$, keeping only the first term for the lower bound, summing the geometric series for the upper bound, and using $2K \in \{M+1, M+2\}$ and $(M+1)/2 \le K \le M/2 + 1$, we obtain, for every $M \in \N$,
    \[
        \frac{(c\sqrt t)^{M+2}}{\sqrt{2M+4}} \;\le\; \E\big[X_t - X_t(M)\big] \;\le\; \sqrt{\frac{2}{\pi(M+1)}}\; \frac{(c\sqrt t)^{M+1}}{1 - c^2 t}.
    \]
    The same holds for $\WT(N)$, uniformly in $N$, with $t$ replaced by $\sigma_N(t)^2 \le t$ (proof of Proposition~\ref{prop:truncation-sharpness}). Remark~\ref{rem:boundary-case} recovers this structure, a horizon of order $c^{-2}$ that shrinks with $p$ and a geometric ratio of order $c\sqrt t$, up to the constants lost in Proposition~\ref{prop:signature-moment-bound}.
\end{example}

\begin{remark}[Truncating by degree] \label{rem:modified-length-truncation}
    We truncate along the word length $\{|\blold{v}| \le M\}$, so that $X_t(M,N)$ only involves the signature up to level $M$. Since $\textnorm{\WT^{\blold{v}}_t(N)}_p \le \psi_{p, t}(\blold{v})$ is of order $t^{\norm{\blold{v}}/2}$ (Proposition~\ref{prop:signature-moment-bound}), truncating along the degree, $\{\norm{\blold{v}} \le M\}$, would instead order the terms by their size, as in the Stratonovich--Taylor expansions of~\citet*{KloedenPlaten1992}. The proof of Lemma~\ref{lem:truncation} applies to this truncation, grouping the words by degree instead of length. There are at most $(d+1)^n$ words of degree $n$ (Lemma~\ref{lem:count-degree}), and $(r+2)^{|\blold{v}|} \le (r+2)^{\norm{\blold{v}}}$, so the per-degree bound reads $C_*\big(2(d+1)(r+2)\sqrt{c_d t}\big)^n (n!)^{-\varepsilon}$. This gives a bound of the same form $(C\sqrt t)^{M+1}/\big((M+1)!\big)^{\varepsilon}$, uniformly in $N \in \Nbar^*$, where the power of $t$ now matches the truncation level exactly. Since $\norm{\blold{v}} \ge |\blold{v}|$, the words of degree at most $M$ have length at most $M$, so the degree-truncated functional is obtained from the signature up to level $M$ by discarding the words of degree larger than $M$. We do not pursue this, since the truncation by length is sufficient for our purposes.
\end{remark}

The processes of Section~\ref{sec:main-results} are $\R^D$-valued, so we extend Lemma~\ref{lem:truncation} to coefficient families in $\mathcal{L}^D_T$.

\begin{lemma}[Products of tails]\label{lem:multiD-truncation}
    Let $\lbold \in \mathcal{L}^D_T$, and let $\varepsilon$ be a common summability exponent of the dominating families $\lbold^{i, \ast}$, $i \in \{1, \dots, D\}$. For $M \in \N$ and $N \in \Nbar^*$, define $x(M,N) \coloneqq x = (x^1, \dots, x^D)$ by $x^i_t \coloneqq \textbrackets{\lbold^i_t}{\WT^{>M}_t(N)}$. Then, for every $\alpha = (\alpha_1, \dots, \alpha_D) \in \N^D$ and $p \geq 1$, setting $B \coloneqq B_{pD|\alpha|}\sqrt{T}$ and $r_M$ accordingly, with $B_r$ as in Lemma~\ref{lem:truncation},
    \[
        \norm{x_t(M,N)^\alpha}_p \lesssim r_M^{|\alpha|}, \qquad t \in [0, T], \ M \in \N, \ N \in \Nbar^*,
    \]
    where $|\alpha| = \sum_{1 \le i \le D} \alpha_i$ and $x_t^\alpha = \prod_{1 \le i \le D} (x^i_t)^{\alpha_i}$.
\end{lemma}
 
\begin{proof}
    The case $\alpha = 0$ is trivial, so assume $|\alpha| \ge 1$. By Hölder's inequality,
    \[
        \|x^\alpha_t\|_p  \le \prod_{i=1}^D \big\| (x^i_t)^{\alpha_i} \big\|_{pD} = \prod_{\substack{1 \le i \le D \\ \alpha_i \ne 0}} \Delta^{\alpha_i}_{pD}(\lbold^i, t, M, N),
    \]
    with $\lbold^i(M,t) \coloneqq \lbold^i_t$ in the notation of Lemma~\ref{lem:truncation}. That lemma, applied with $p$ replaced by $pD$ and $q = \alpha_i$, gives
    \[
        \Delta^{\alpha_i}_{pD}(\lbold^i, t, M, N) \le \left( C^i_*\, E_\varepsilon\big(B_{pD\alpha_i}\sqrt T\big)\, \frac{\big(B_{pD\alpha_i}\sqrt t\big)^{M+1}}{\big((M+1)!\big)^{\varepsilon}}\right)^{\alpha_i}.
    \]
    Since $r \mapsto B_r$ and $E_\varepsilon$ are non-decreasing, $\alpha_i \le |\alpha|$, $t \le T$, and $(M+1)! \ge M!$, each factor is at most $\big( C^i_* E_\varepsilon(B)\, B^{M+1} (M!)^{-\varepsilon} \big)^{\alpha_i}$. Multiplying over $i$ gives the result.
\end{proof}

\begin{proof}[Proof of Theorem~\ref{th:Lp-truncation-bound}]
    Fix $p \ge 1$, $M \in \N$, $N \in \Nbar^*$ and $t \in [0,T]$, and write
    \[
        X \coloneqq X_t(\infty,N), \qquad Y \coloneqq X_t(M,N), \qquad x \coloneqq X - Y = x_t(M,N),
    \]
    with $x_t(M,N)$ as in Lemma~\ref{lem:multiD-truncation}. The constant $B$ in $r_M$ is chosen below. All the implicit constants below are independent of $t$, $M$ and $N$.

    \emph{Polynomial test functions.} Let $g \in \mathcal{P}(\R^D, \R)$ whose r.e. is $a = 1$, and write $g(z) = \sum_{|\alpha| \le m} c_\alpha z^\alpha$ for some $m \in \N$. By the triangle inequality, it suffices to treat a monomial $z \mapsto z^\alpha$ with $1 \le |\alpha| \le m$. Expanding $X = Y + x$ coordinatewise, $X^\alpha - Y^\alpha = \sum_{0 < \beta \le \alpha} \binom{\alpha}{\beta}\,Y^{\alpha-\beta}\, x^\beta$, where $\beta \le \alpha$ means $\beta_i \le \alpha_i$ for every $i$, $0 < \beta$ means $\beta \ne 0$, and $\binom{\alpha}{\beta} = \prod_i \binom{\alpha_i}{\beta_i}$. By the Cauchy--Schwarz inequality, $\textnorm{Y^{\alpha-\beta} x^\beta}_p \le \textnorm{Y^{\alpha-\beta}}_{2p}\, \textnorm{x^\beta}_{2p}$. First, by Hölder's inequality and Corollary~\ref{cor:well-posed}(ii),
    \[
        \norm{Y^{\alpha-\beta}}_{2p} \le \prod_{i=1}^D \norm{Y^i}_{2pD(\alpha_i-\beta_i)}^{\alpha_i-\beta_i} \le \prod_{i=1}^D \big(K^{(i)}_{2pDm}\big)^{\alpha_i-\beta_i},
    \]
    where $K^{(i)}_q$ is the constant of Corollary~\ref{cor:well-posed} for $\lbold^{i}$, and $\norm{ \, \cdot \, }_{q} \le \norm{ \, \cdot \, }_{2pDm}$ for $q \le 2pDm$ (for $\alpha_i = \beta_i$ the factor equals $1$). This bound is uniform in $M$, $N$ and $t$. Second, Lemma~\ref{lem:multiD-truncation} with $2p$ in place of $p$ gives $\textnorm{x^\beta}_{2p} \lesssim \big(B_{2pD|\beta|}^M\, T^{M/2}\, (M!)^{-\varepsilon}\big)^{|\beta|}$. Set $B \coloneqq B_{2pDm}\sqrt T$. Since $r \mapsto B_r$ is non-decreasing, $|\beta| \ge 1$ and $(r_M)_M$ is bounded, we get $\|x^\beta\|_{2p} \lesssim r_M^{|\beta|} \lesssim r_M$. Summing over the finitely many $\beta$ and $\alpha$ yields $\|g(X) - g(Y)\|_p \lesssim r_M$, which is~\eqref{eq:truncation-rate} with $a = 1$.

    \emph{Hölder test functions.} Let $g \in \mathcal{H}(\R^D, \R)$ with r.e. $a$, so that $|g(z) - g(z')| \le K_g \tnorm{z - z'}^a$ for all $z, z' \in \R^D$. Since $\tnorm{z} \le \sum_{i=1}^D |z^i|$ and $u \mapsto u^a$ is subadditive on $\R_+$ for $a \in (0,1]$, $|g(X) - g(Y)| \le K_g \sum_{i=1}^D |x^i|^a$. By Minkowski's inequality and $\big\||x^i|^a\big\|_p = \|x^i\|_{ap}^a$, together with the monotonicity of $L^q$ norms ($ap \le p$),
    \[
        \|g(X) - g(Y)\|_p \le K_g \sum_{i=1}^D \|x^i\|_{ap}^a \le K_g \sum_{i=1}^D \|x^i\|_p^a .
    \]
    Lemma~\ref{lem:truncation} with $q = 1$ gives $\|x^i\|_p \lesssim r_M$ with $B \coloneqq B_p\sqrt T$, uniformly in $t$ and $N$, which is~\eqref{eq:truncation-rate}.
\end{proof}

We now prove that the truncation rate is sharp.

\begin{proof}[Proof of Proposition~\ref{prop:truncation-sharpness}]
    We simplify $W^1$ into $W$. We write the expected error as a series of Gaussian moments with positive terms, bound it below by its first term, and compare the result with the upper bound of Theorem~\ref{th:Lp-truncation-bound}. Fix $N \in \Nbar^*$ and $t \in (0,T]$, and write $\Wt_t(\infty) = W_t$. The iterated integrals of a single coordinate of a path are the powers of its increment, so that, pathwise, $\WT^{\blold{1}^{\otimes n}}_t(N) = \Wt_t^{n}/n!$. Moreover, $\Wt_t$ is centred Gaussian with variance $\sigma^2 \coloneqq \sigma_N(t)^2$ satisfying
    \begin{equation}\label{eq:sharpness-variance}
        \frac{t^2}{T} \le \sigma^2 \le t .
    \end{equation}
    For $N = \infty$, $\sigma^2 = t \ge t^2/T$. For $N < \infty$, write $t = t_k + \lambda h$ with $\lambda \in (0,1]$. Then $\widetilde W_t = (1-\lambda) W_{t_k} + \lambda W_{t_{k+1}}$, so that $\sigma^2 = t_k + \lambda^2 h$. Hence $t - \sigma^2 = \lambda(1-\lambda)h \ge 0$, and
    \[
        \sigma^2 t_{k+1} - t^2 = (t_k + \lambda^2 h)(t_k + h) - (t_k + \lambda h)^2 = t_k h (1-\lambda)^2 \ge 0,
    \]
    which, together with $t_{k+1} \le T$, gives $\sigma^2 \ge t^2/(t_k+h) \ge t^2/T$. By Corollary~\ref{cor:well-posed} (finiteness of $K_1$), $\sum_n \lbold^{(\varepsilon),\blold{1}^{\otimes n}} \textnorm{\WT^{\blold{1}^{\otimes n}}_t(N)}_1 < \infty$, so that expectation and summation may be interchanged. Since $\E[\Wt_t^{n}]$ vanishes for odd $n$ and equals $\sigma^{2j}(2j-1)!!$ for $n = 2j$, we obtain, using $(2j)! = (2j-1)!!\,2^j j!$,
    \[
        \E\big[X_t(\infty,N) - X_t(M,N)\big] = \sum_{2j > M} \big[(2j)!\big]^{\frac12-\varepsilon}\,\frac{(2j-1)!!}{(2j)!}\,\sigma^{2j} = \sum_{2j > M} a_j, \qquad a_j \coloneqq \big[(2j)!\big]^{\frac12-\varepsilon}\,\frac{(\sigma^2/2)^j}{j!}.
    \]
    Every $a_j$ is positive, so the sum is bounded below by its first term, of index $J \coloneqq \lfloor M/2 \rfloor + 1$: $\E\big[X_t(\infty,N) - X_t(M,N)\big] \ge a_J$. We use the elementary bounds $n! \ge (n/e)^n$ and $n! \le n^n$, valid for $n \ge 1$. The first gives $\log\big((2J)!\big) \ge 2J\log(2J) - 2J$, and the second gives $\log J! \le J\log J$. Since $1/2 - \varepsilon \ge 0$,
    \begin{align*}
        \log a_J
        &\ge (1-2\varepsilon)\,J\log(2J) - (1-2\varepsilon)\,J - J\log J + J\log(\sigma^2/2) \\
        &= -2\varepsilon J\log J + J\Big[(1-2\varepsilon)(\log 2 - 1) + \log(\sigma^2/2)\Big].
    \end{align*}
    Now $(1-2\varepsilon)(\log 2 - 1) \ge -1$ and, by~\eqref{eq:sharpness-variance}, $\log(\sigma^2/2) \ge \log(t^2/(2T)) \ge -|\log(t^2/(2T))|$. Hence \begin{equation}\label{eq:sharpness-aJ} 
        \log a_J \ge -2\varepsilon J\log J - \big(c_t - 1\big)J.
    \end{equation}
    Since $M \ge 2$, we have $J \le M/2 + 1 \le M$, so $J \le M$ and $2J \le M+2$. Therefore
    \[
        2\varepsilon J\log J \le \varepsilon(M+2)\log M = \varepsilon M\log M + 2\varepsilon\log M \le \varepsilon M\log M + M,
    \]
    using $2\varepsilon \le 1$ and $\log M \le M$. Combining with \eqref{eq:sharpness-aJ} yields $\log a_J \ge -\varepsilon M\log M - c_t M$, which is the second inequality in~\eqref{eq:sharpness}. The constant $c_t$ does not depend on $N$. The first inequality in~\eqref{eq:sharpness} is $\E Z \le |\E Z| \le \norm{Z}_1 \le \norm{Z}_p$. Conversely, Theorem~\ref{th:Lp-truncation-bound} with $g = \mathrm{Id}$ gives $\norm{X_t(\infty,N) - X_t(M,N)}_p \le C r_M$ uniformly in $N$, and $\log M! \ge M\log M - M$ yields $\log\norm{X_t(\infty,N) - X_t(M,N)}_p \le -\varepsilon M\log M + O(M)$. The two bounds match up to $O(M)$, uniformly in $N$.
\end{proof}

\begin{remark} \label{rem:sharpness-comments}
    \begin{enumerate}[label=(\roman*)]
        \item[]
        \item For $\varepsilon = 1/2$, keeping only the even coefficients recovers $\cosh(\Wt^1_t(N))$ exactly, and $\E[X_t(\infty,N) - X_t(M,N)] = \sum_{2j > M} (\sigma_N(t)^2/2)^j/j! = \exp\{-\tfrac12 M\log M + O(M)\}$.
        \item The constant $c_t$ blows up as $t \downarrow 0$, as it must: $X_0 - X_0(M) = 0$. The lower bound is uniform on $[t_0, T]$ for every $t_0 > 0$.
        \item The proposition concerns $g = \mathrm{Id}$, hence polynomial test functions. For an $a$-Hölder $g$, the upper bound of Theorem~\ref{th:Lp-truncation-bound} carries the exponent $\varepsilon a$; whether this is sharp for $a < 1$ is left open.
    \end{enumerate}
\end{remark}

\subsection{Strong discretisation error} \label{sec:strong-discretisation}

We first prove the convergence componentwise, word by word, before passing to linear combinations. Throughout this section, we fix $N \in \N^*$, recall that $h \coloneqq T/N$ and $t_k \coloneqq kh$, and write $\WT \coloneqq \WT(N)$ for the signature of the interpolation~\eqref{eq:discretised-signature}, making the dependence on $N$ explicit only when several values are involved. For $t \in (0,T]$, we let $n \in \intset{0}{N-1}$ be such that $t \in (t_n, t_{n+1}]$, and define the modified times $\tilde t_i \coloneqq t_i$ for $i \le n$ and $\tilde t_i \coloneqq t$ for $n < i \le N$. For $p \ge 1$, we write $\psi \coloneqq \psi_{p,T}$, and we set $\Delta^{\blold{v}}(t) \coloneqq \WH^{\blold{v}}_t - \WT^{\blold{v}}_t$, $\blold{v} \in V$. The error is first decomposed into local errors, one per step of the grid (Lemma~\ref{lem:error-decomposition}). A direct moment bound controls each of them, except for the Lévy-area words $\blold{ij}$, $i \ne j$ (Lemma~\ref{lem:local-errors}), for which the cancellations between the steps must be exploited (Lemma~\ref{lem:levy-area-words}). By Proposition~\ref{prop:signature-moment-bound}, $\|\WH^{\blold{v}}_{s,u}\|_p \vee \|\WT^{\blold{v}}_{s,u}\|_p \le \psi_{p,u-s}(\blold{v}) \le h^{\norm{\blold{v}}/2}\,\psi(\blold{v})/T^{\norm{\blold{v}}/2}$ whenever $0 \le u - s \le h$: on a single step, a component of degree $n$ is of order $h^{n/2}$. We now decompose the error along the grid.

\begin{lemma}[Chen decomposition of the error]\label{lem:error-decomposition}
    For every $\blold{v} \in V$,
    \[
        \Delta^{\blold{v}}(t) = \sum_{\blold{v_1 v_2 v_3} = \blold{v}} \Delta^{\blold{v_1}, \blold{v_2}, \blold{v_3}}(t), \qquad \Delta^{\blold{v_1}, \blold{v_2}, \blold{v_3}}(t) \coloneqq \sum_{k = 1}^{n+1}\WH_{\tilde{t}_{k-1}}^{\blold{v_1}} \left(\WH_{\tilde{t}_{k-1}, \tilde{t}_k}^{\blold{v_2}} - \WT_{\tilde{t}_{k-1}, \tilde{t}_k}^{\blold{v_2}}\right) \WT_{\tilde{t}_{k}, t}^{\blold{v_3}}.
    \]
    Moreover, for each $k$, the three factors of the $k$-th summand depend on the Brownian increments over disjoint intervals, and are therefore independent.
\end{lemma}

\begin{proof}
    Write $\boldsymbol{a}_k \coloneqq \WH_{\tilde t_{k-1}, \tilde t_k}$ and $\boldsymbol{b}_k \coloneqq \WT_{\tilde t_{k-1}, \tilde t_k}$, $k \in \intset{1}{n+1}$. By Chen's identity~\eqref{eq:chen-identity}, which holds for both paths, $\WH_t = \boldsymbol{a}_1 \otimes \cdots \otimes \boldsymbol{a}_{n+1}$ and $\WT_t = \boldsymbol{b}_1 \otimes \cdots \otimes \boldsymbol{b}_{n+1}$. The telescoping identity
    \[
        \bigotimes_{l=1}^{n+1} \boldsymbol{a}_l - \bigotimes_{l=1}^{n+1} \boldsymbol{b}_l = \sum_{k=1}^{n+1} \Big[\bigotimes_{l<k} \boldsymbol{a}_l\Big] \otimes (\boldsymbol{a}_k - \boldsymbol{b}_k) \otimes \Big[\bigotimes_{l>k} \boldsymbol{b}_l\Big]
    \]
    and Chen's identity again give $\WH_t - \WT_t = \sum_{k=1}^{n+1} \WH_{\tilde t_{k-1}} \otimes (\boldsymbol{a}_k - \boldsymbol{b}_k) \otimes \WT_{\tilde t_k, t}$, and extracting the $\blold{v}$-coordinate yields the decomposition. The three factors are functionals of the Brownian path on $[0, \tilde t_{k-1}]$, $[t_{k-1}, t_k]$ and $[t_k, t_{n+1}]$ respectively, hence independent.
\end{proof}

A direct moment bound controls each local error, except for the Lévy-area words.

\begin{lemma}[Local errors]\label{lem:local-errors}
    Let $p \ge 1$ and $\blold{v_1}, \blold{v_2}, \blold{v_3} \in V$. If $\blold{v_2} = \empt$, then $\Delta^{\blold{v_1}, \blold{v_2}, \blold{v_3}}(t) = 0$. If $\blold{v_2} \notin \{\blold{ij} : i \ne j \in \intset{1}{d}\}$, then
    \[
        \textnorm{\Delta^{\blold{v_1}, \blold{v_2}, \blold{v_3}}(t)}_p \le 2\,\psi(\blold{v_1})\,\psi(\blold{v_2})\,\psi(\blold{v_3})\,N^{-1/2}.
    \]
\end{lemma}

\begin{proof}
    If $\blold{v_2} = \empt$, every local error vanishes. Otherwise, by the independence in Lemma~\ref{lem:error-decomposition} and Proposition~\ref{prop:signature-moment-bound},
    \begin{equation}\label{eq:decomposition_1}
        \textnorm{\Delta^{\blold{v_1}, \blold{v_2}, \blold{v_3}}(t)}_p \le \psi(\blold{v_1})\,\psi(\blold{v_3}) \sum_{k=1}^{n+1} \textnorm{\WH^{\blold{v_2}}_{\tilde t_{k-1}, \tilde t_k} - \WT^{\blold{v_2}}_{\tilde t_{k-1}, \tilde t_k}}_p .
    \end{equation}
    Since $\tilde t_k - \tilde t_{k-1} \le h$ and $\psi_{p,s}(\blold{v}) \le N^{-\norm{\blold{v}}/2}\psi(\blold{v})$ for $s \le h$, Proposition~\ref{prop:signature-moment-bound} bounds each local error by $2\psi(\blold{v_2})N^{-\norm{\blold{v_2}}/2}$. If $\blold{v_2} \in \{\blold{i}, \blold{ii} : 0 \le i \le d\}$, the coordinate $\X^{\blold{v_2}}_{s,u}$ only depends on the increment $X^i_u - X^i_s$ (Section~\ref{sec:signature}); since $\Wt$ and $\Wh$ coincide at the grid points, the local errors vanish for $k \le n$, and only the last one remains: the sum is at most $2\psi(\blold{v_2})N^{-1/2}$, as $\norm{\blold{v_2}} \ge 1$. Otherwise, since $\blold{v_2}$ is not of the form $\blold{ij}$, $i \ne j$, it has degree at least $3$, and the $n + 1 \le N$ local errors sum to at most $2\psi(\blold{v_2})N^{1 - 3/2}$.
\end{proof}

When $i \ne j$, the argument of Lemma~\ref{lem:local-errors} only gives $\|\Delta^{\blold{v_1}, \blold{ij}, \blold{v_3}}(t)\|_p \lesssim 1$. The local errors $\WH^{\blold{ij}}_{t_{k-1}, t_k} - \WT^{\blold{ij}}_{t_{k-1}, t_k}$ are, however, the Lévy areas of the successive steps, which are centred, symmetric and independent. Expanding the $p$-th moment, only the terms in which every Lévy area appears an even number of times contribute at leading order, and counting them with Rosenthal's inequality gives the rate $N^{-1/2}$. The proof is given in Appendix~\ref{app:Lp-element-discretisation}.

\begin{lemma}[Lévy-area words]\label{lem:levy-area-words}
    Let $p \in 2\N^*$, $i \ne j \in \intset{1}{d}$ and $\blold{v_1}, \blold{v_3} \in V$, and set $\blold{v} \coloneqq \blold{v_1}\,\blold{ij}\,\blold{v_3}$. There exists a polynomial $P$, depending only on $p$, such that
    \[
        \textnorm{\Delta^{\blold{v_1}, \blold{ij}, \blold{v_3}}(t)}_p \le P(\norm{\blold{v}})\,(2p+3)^{\norm{\blold{v}}/2}\,\psi(\blold{v})\,N^{-1/2}.
    \]
\end{lemma}

We can now prove the word-by-word estimate.

\begin{proposition}[Strong error on words]\label{prop:Lp-element-discretisation}
    Let $p \ge 1$. Then, there exists $C > 0$, depending only on $d$ and $p$ such that
    \[
        \norm{\WH^{\blold{v}}_t - \WT^{\blold{v}}_t(N)}_p \le C\,e^{C\norm{\blold{v}}}\,\frac{T^{\norm{\blold{v}}/2}}{\lfloor\norm{\blold{v}}/2\rfloor!}\,N^{-1/2}, \qquad \blold{v} \in V, \ t \in [0, T], \ N \in \N^*.
    \]
\end{proposition}

\begin{proof}
    By the monotonicity of $L^p$ norms, it suffices to consider $p \in 2\N^*$, and the case $t = 0$ is immediate. By Lemma~\ref{lem:error-decomposition} and the triangle inequality,
    \[
        \textnorm{\Delta^{\blold{v}}(t)}_p \le \sum_{\blold{v_1 v_2 v_3} = \blold{v}} \textnorm{\Delta^{\blold{v_1}, \blold{v_2}, \blold{v_3}}(t)}_p .
    \]
    If $\blold{v_2}$ is not of the form $\blold{ij}$, $i \ne j$, Lemmas~\ref{lem:local-errors} and~\ref{lem:prod-psi} give
    \[
        \textnorm{\Delta^{\blold{v_1}, \blold{v_2}, \blold{v_3}}(t)}_p \le 2 \cdot 3^{\norm{\blold{v}}/2} \big(1 + \norm{\blold{v}}/2\big)^{3/2}\, \psi(\blold{v})\, N^{-1/2},
    \]
    and otherwise Lemma~\ref{lem:levy-area-words} applies. Summing over the $(|\blold{v}|+1)(|\blold{v}|+2)/2$ decompositions $\blold{v_1v_2v_3}=\blold{v}$ and using $|\blold{v}|\le\norm{\blold{v}}$, all polynomial factors, as well as $(2p+3)^{\norm{\blold{v}}/2}$, $(p+2)^{|\blold{v}|}$ and $(1+d/2)^{\norm{\blold{v}}/2}$, are absorbed in $C \, \exp (C\norm{\blold{v}})$: $\textnorm{\Delta^{\blold{v}}(t)}_p\le C\,\exp(C\norm{\blold{v}}) \, T^{\norm{\blold{v}}/2} \, N^{-1/2}/\lfloor\norm{\blold{v}}/2\rfloor!$, with $C$ depending only on $d$ and $p$.
\end{proof}

In dimension one, and at the grid points, the rate of Proposition~\ref{prop:Lp-element-discretisation} improves.

\begin{proposition}[Grid points in dimension one]\label{prop:Lp-discretisation-d1}
    Assume $d = 1$, and let $p \geq 1$. Then, there exists $C > 0$, depending only on $p$, such that, for every $\blold{v} \in V$,
    \[
        \norm{\WH^{\blold{v}}_t - \WT^{\blold{v}}_t(N)}_p \le C\,e^{C\norm{\blold{v}}}\,\frac{T^{\norm{\blold{v}}/2}}{\lfloor\norm{\blold{v}}/2\rfloor!}\,N^{-1}, \qquad N \in \N^*, \ t \in \{t_k : 0 \leq k \leq N\}.
    \]
    In particular, for a fixed $t \in (0, T]$ with $t/T \in \Q$, the rate $N^{-1}$ holds along the meshes that contain $t$, that is, along $\{N \in \N^* : tN/T \in \N\}$.
\end{proposition}
 
\begin{proof}
    We use Lemma~\ref{lem:error-decomposition} and only indicate the changes with respect to Lemmas~\ref{lem:local-errors} and~\ref{lem:levy-area-words}. Since $t = t_{n+1}$ is a grid point, the modified times satisfy $\tilde{t}_k = t_k$ for every $k \leq n + 1$. The last interval is therefore a full step, so that all the local errors in the decomposition of Lemma~\ref{lem:error-decomposition} are of the form $\WH^{\blold{v_2}}_{t_{k-1}, t_k} - \WT^{\blold{v_2}}_{t_{k-1}, t_k}$, $1 \leq k \leq n + 1$. With $d = 1$, the alphabet is $\{\blold{0}, \blold{1}\}$, and there is no Lévy area between Brownian components: the critical case $\blold{v_2} = \blold{ij}$, $i \ne j$, of Lemma~\ref{lem:levy-area-words} does not occur. As in the proof of Lemma~\ref{lem:local-errors}, the local errors vanish for every $k$ if $\blold{v_2}$ consists of a single repeated letter, since the two paths coincide at the grid point. The remaining words fall into two classes.
    \begin{itemize}
        \item If $\norm{\blold{v_2}} \geq 4$, \eqref{eq:decomposition_1} and the bound $2\psi(\blold{v_2})N^{-\norm{\blold{v_2}}/2}$ on each local error (proof of Lemma~\ref{lem:local-errors}) gives $\textnorm{\Delta^{\blold{v_1}, \blold{v_2}, \blold{v_3}}(t)}_p \le 2 \, \psi(\blold{v_1})\psi(\blold{v_2})\psi(\blold{v_3}) \, N^{-1} $.
        \item If $\norm{\blold{v_2}} = 3$, then $\blold{v_2} \in \{\blold{01}, \blold{10}\}$, and the crude bound only gives $N \cdot N^{-3/2} = N^{-1/2}$. In that case the local errors are centred and independent, and we argue as in the proof of Lemma~\ref{lem:levy-area-words}. On $[t_{k-1}, t_k]$, with midpoint $t_{k-1/2} \coloneqq (t_{k-1}+t_{k})/2$,
        \[
            Z_k \coloneqq \WH^{\blold{01}}_{t_{k-1}, t_k} - \WT^{\blold{01}}_{t_{k-1}, t_k} = \int_{t_{k-1}}^{t_k} (s - t_{k-1/2}) \, \diff W_s = -\big(\WH^{\blold{10}}_{t_{k-1}, t_k} - \WT^{\blold{10}}_{t_{k-1}, t_k}\big),
        \]
        where the last identity follows from the shuffle identity $\blold{0} \shuffle \blold{1} = \blold{01} + \blold{10}$, which is reproduced exactly by the scheme, i.e. $\WH^{\blold{01}} + \WH^{\blold{10}} = \WH^{\blold{0}}\WH^{\blold{1}} = 2\WT^{\blold{01}}$. The variables $Z_k$, $1 \leq k \leq n+1$, are independent, each $Z_k$ depends only on the Brownian path on $[t_{k-1}, t_k]$, and $Z_k$ is centred Gaussian with variance $h^3/12$, so that $\E[Z_k^{2q}] \lesssim h^{3q}$ and $\E[Z_k^{2q+1}] = 0$. These are the only properties of the Lévy area $A_{\blold{ij}}(t_{k-1}, t_k)$ used in the proof of Lemma~\ref{lem:levy-area-words}, with the order $h$ of the Lévy areas replaced by the order $h^{3/2}$ of $Z_k$. The same argument therefore gives
        \[
            \norm{\sum_{k=1}^{n+1} \WH^{\blold{v_1}}_{t_{k-1}} \, Z_k \, \WT^{\blold{v_3}}_{t_k, t} }_p^p \lesssim N^{-3p/2} \, S_{n+1, N} \lesssim N^{-p},
        \]
        by Lemma~\ref{lem:combinatorial-sum} (note $n+1 \le N$). That is, $\textnorm{\Delta^{\blold{v_1}, \blold{v_2}, \blold{v_3}}(t)}_p \lesssim N^{-1}$.
    \end{itemize}
    The constants are tracked as in the proofs of Lemma~\ref{lem:levy-area-words} and Proposition~\ref{prop:Lp-element-discretisation}, with $h^{3/2}$ in place of $h$, and summing over the finitely many decompositions $\blold{v_1}\blold{v_2}\blold{v_3} = \blold{v}$ concludes the proof.
\end{proof}

We now pass from words to linear functionals.

\begin{lemma}[Strong discretisation error of linear functionals]\label{lem:Lp-discretisation}
    Let $\lbold : \Nbar \times [0,T] \to T((\R^{d+1}))$ and $\lbold^\ast \in \mathcal L^1$ with $\lbold^{\ast,\blold{v}} \le C_* \varphi_\varepsilon(\blold{v})$, $\blold{v} \in V$, for some $\varepsilon \in (0,1/2]$ and $C_* > 0$, and assume that $\lbold(M,t) \preceq \lbold^\ast$ for every $M \in \Nbar$ and $t \in [0,T]$. Set $\tilde{r}_N \coloneqq N^{-1/2}$. Then, for every $p \ge 1$, there exists $C > 0$, depending only on $d, T, p, \varepsilon$ and $C_*$, such that
    \begin{equation}\label{eq:Lp-discretisation}
        \Big\| \brackets{\lbold(M,t)}{\WH_t} - \brackets{\lbold(M,t)}{\WT_t(N)} \Big\|_p \le C\, \tilde{r}_N, \qquad t \in [0,T],\ M \in \Nbar,\ N \in \N^* .
    \end{equation}
    If $d = 1$, \eqref{eq:Lp-discretisation} also holds with $\tilde{r}_N \coloneqq N^{-1}$ at the grid points $t \in \{t_k : 0 \le k \le N\}$.
\end{lemma}

\begin{proof}
    Fix $t$, $M$, $N$ and write $\Delta^{\blold{v}} \coloneqq \Delta^{\blold{v}}(t)$ and $n \coloneqq \norm{\blold{v}}$. By Corollary~\ref{cor:well-posed}(i), applied with $N = \infty$ and with $N$, both series $\textbrackets{\lbold(M,t)}{\WH_t}$ and $\textbrackets{\lbold(M,t)}{\WT_t}$ converge absolutely almost surely, so that their difference equals $\sum_{\blold{v}} \lbold^{\blold{v}}(M,t)\, \Delta^{\blold{v}}$, the series converging absolutely almost surely. By monotone convergence and Minkowski's inequality, $\textnorm{\textstyle{\sum_{\blold{v}}} \lbold^{\blold{v}}(M,t)\, \Delta^{\blold{v}}}_p \le \textstyle{\sum_{\blold{v}}} \lbold^{\ast,\blold{v}}\, \textnorm{\Delta^{\blold{v}}}_p$, and it remains to show that the right-hand side is at most $C N^{-1/2}$. By Proposition~\ref{prop:Lp-element-discretisation}, there is $C_1 > 0$, depending only on $d$ and $p$, such that $\norm{\Delta^{\blold{v}}}_p \le C_1 e^{C_1 n}\, T^{n/2}\, N^{-1/2} / \lfloor n/2 \rfloor !$. The inequality $n! \le 2^n (n/2+1) (\lfloor n/2\rfloor!)^2$, used for~\eqref{eq:psi-bound}, gives $1/\lfloor n/2 \rfloor ! \le 2^{n/2} \sqrt{n/2+1}\,/\sqrt{n!}$. Together with $\lbold^{\ast,\blold{v}} \le C_* (n!)^{\frac12-\varepsilon}$, this yields 
    \[
        \lbold^{\ast,\blold{v}}\, \norm{\Delta^{\blold{v}}}_p \le C_* C_1\, \big(e^{C_1}\sqrt{2T}\big)^{n}\, \sqrt{n/2+1}\; (n!)^{-\varepsilon}\, N^{-1/2}.
    \]
    There are at most $(d+1)^n$ words of degree $n$ (Lemma~\ref{lem:count-degree}). Setting $A \coloneqq (d+1)\, e^{C_1}\sqrt{2T}$, we obtain 
    \[
        \textstyle{\sum_{\blold{v}}} \lbold^{\ast,\blold{v}}\, \norm{\Delta^{\blold{v}}}_p \le C_* C_1 \Big( \sum_{n \ge 0} A^n \sqrt{n/2+1}\; (n!)^{-\varepsilon} \Big) N^{-1/2},
    \]
    the series being finite since $\varepsilon > 0$. This is the first part of the claim, with a constant depending only on $d, T, p, \varepsilon$ and $C_*$. If $d = 1$ and $t$ is a grid point, Proposition~\ref{prop:Lp-discretisation-d1} gives the same bound on $\norm{\Delta^{\blold{v}}}_p$ with $N^{-1}$ in place of $N^{-1/2}$, and the summation above is unchanged.
\end{proof}

As for the truncation error, Theorem~\ref{th:Lp-discretisation-bound} follows by treating the two classes of test functions.

\begin{proof}[Proof of Theorem~\ref{th:Lp-discretisation-bound}]
    Fix $p \ge 1$, $M \in \Nbar$, $N \in \Nbar^*$ and $t \in [0,T]$, write $Y \coloneqq X_t(M,\infty)$ and $Y^N \coloneqq X_t(M,N)$, and let $\tilde{r}_N$ be as in Lemma~\ref{lem:Lp-discretisation}. All the constants below are independent of $t$, $M$ and $N$. \emph{Polynomial test functions.} By the triangle inequality, it suffices to treat a monomial $z \mapsto z^\alpha$, $\alpha \in \N^D$. Set 
    \[
        \lbold^\alpha(M,t) \coloneqq \bigshuffle_{i=1}^D \big(\lbold^i_{t,\le M}\big)^{\shuffle \alpha_i}, \qquad \lbold^{\alpha,*} \coloneqq \bigshuffle_{i=1}^D \big(\lbold^{*}\big)^{\shuffle \alpha_i}.
    \]
    By Section~\ref{sec:words}, $\lbold^\alpha(M,t) \preceq \lbold^{\alpha,*}$ for every $M$ and $t$, and $\lbold^{\alpha,*} \in \mathcal{L}^1$ with some summability exponent $\varepsilon_\alpha > 0$ by Proposition~\ref{prop:L1-stability}. By the shuffle identity for both signatures (Lemma~\ref{lem:shuffle-L1}), $Y^\alpha = \textbrackets{\lbold^\alpha(M,t)}{\WH_t}$ and $(Y^N)^\alpha = \textbrackets{\lbold^\alpha(M,t)}{\WT_t(N)}$. Lemma~\ref{lem:Lp-discretisation}, applied with $\lbold^{\alpha,*}$, gives $\|Y^\alpha - (Y^N)^\alpha\|_p \lesssim \tilde{r}_N$, and summing over the finitely many $\alpha$ yields $\|g(Y) - g(Y^N)\|_p \lesssim \tilde{r}_N$. \emph{Hölder test functions.} As in the proof of Theorem~\ref{th:Lp-truncation-bound}, $\|g(Y) - g(Y^N)\|_p \le K_g \sum_{1 \le i \le D} \|Y^i - Y^{N,i}\|_p^a$. Lemma~\ref{lem:Lp-discretisation}, applied with $\lbold(M,t) \coloneqq \lbold^i_{t,\le M} \preceq \lbold^\ast$, gives $\|Y^i - Y^{N,i}\|_p \lesssim \tilde{r}_N$, hence $\|g(Y) - g(Y^N)\|_p \lesssim \tilde{r}_N^a$. Taking $\tilde{r}_N = N^{-1/2}$ gives the first claim, with $a = 1$ for polynomial $g$; taking $\tilde{r}_N = N^{-1}$ when $d = 1$ and $t$ is a grid point gives the second.
\end{proof}

\subsection{Weak discretisation error} \label{sec:weak-discretisation}

The weak truncation error does not improve on the strong one (Proposition~\ref{prop:truncation-sharpness}). The weak discretisation error does: by Proposition~\ref{prop:expected-signature-bound}, the expected signatures of $\WH$ and of $\WT(N)$ are non-negative, vanish in odd degree, and coincide in degree two on every full step of the grid. The local errors therefore start at degree four, and the rate doubles.

\begin{proposition}[Weak error on words]\label{prop:weak-element-discretisation}
    For every $\blold{v} \in V$, $t \in [0,T]$ and $N \in \N^*$,
    \[
        \Big| \E\big[\WH^{\blold{v}}_t - \WT^{\blold{v}}_t(N)\big] \Big| \le \indic{\norm{\blold{v}} \in 2\N}\; 2(|\blold{v}|+1)(|\blold{v}|+2)\, \frac{(3c_dT)^{\norm{\blold{v}}/2}}{(\norm{\blold{v}}/2)!}\; N^{-1}.
    \]
\end{proposition}

\begin{proof}
    In this proof only, write $\psi \coloneqq \psi_T$, the bound of Proposition~\ref{prop:expected-signature-bound}, instead of $\psi_{p,T}$. With the notation of Lemma~\ref{lem:error-decomposition},
    \[
        \E\big[\Delta^{\blold{v}}(t)\big] = \sum_{\blold{v_1 v_2 v_3} = \blold{v}} \sum_{k=1}^{n+1} \E\big[\WH^{\blold{v_1}}_{\tilde{t}_{k-1}}\big]\, \E\big[\WH^{\blold{v_2}}_{\tilde{t}_{k-1},\tilde{t}_k} - \WT^{\blold{v_2}}_{\tilde{t}_{k-1},\tilde{t}_k}\big]\, \E\big[\WT^{\blold{v_3}}_{\tilde{t}_k,t}\big],
    \]
    the three factors depending on the Brownian increments over disjoint intervals. By Proposition~\ref{prop:expected-signature-bound}, the first and third expectations are bounded by $\psi(\blold{v_1})$ and $\psi(\blold{v_3})$. Since both expected signatures in the middle factor are non-negative, the same proposition bounds it by $2\,(c_d(\tilde{t_k} - \tilde{t}_{k-1}))^{\norm{\blold{v_2}}/2}/(\norm{\blold{v_2}}/2)!$, and shows that it vanishes when $\norm{\blold{v_2}}$ is odd. Moreover:
    \begin{itemize}
        \item if $\blold{v_2} = \empt$, the middle factor vanishes;
        \item if $k \le n$, the step is full and the middle factor also vanishes when $\norm{\blold{v_2}} = 2$, that is, for $\blold{v_2} \in \{\blold{0}, \blold{ii}, \blold{ij}\}$: both expectations equal $h$, $h/2$ and $0$ respectively. Hence it is non-zero only if $\norm{\blold{v_2}} \ge 4$, and then it is at most $2 (c_d h)^{\norm{\blold{v_2}}/2}/(\norm{\blold{v_2}}/2)! \le 2 \psi(\blold{v_2})\, N^{-2}$;
        \item if $k = n+1$, it is non-zero only if $\norm{\blold{v_2}} \ge 2$, and then it is at most $2 \psi(\blold{v_2})\, N^{-1}$, since $\tilde t_{n+1} - \tilde t_n \le h$.
    \end{itemize}
    Using $n \le N$, we obtain
    \[
        \big|\E[\Delta^{\blold{v}}(t)]\big| \le \sum_{\blold{v_1 v_2 v_3} = \blold{v}} \psi(\blold{v_1}) \psi(\blold{v_2}) \psi(\blold{v_3})\,\big(2nN^{-2} + 2N^{-1}\big) \le 4N^{-1} \sum_{\blold{v_1 v_2 v_3} = \blold{v}} \psi(\blold{v_1}) \psi(\blold{v_2}) \psi(\blold{v_3}).
    \]
    A summand is non-zero only if the three degrees are even, in which case $\norm{\blold{v}}$ is even and, writing $a_i \coloneqq \norm{\blold{v_i}}/2$ and $A \coloneqq a_1 + a_2 + a_3 = \norm{\blold{v}}/2$, the multinomial bound $A!/(a_1!a_2!a_3!) \le 3^A$ gives $\psi(\blold{v_1}) \psi(\blold{v_2}) \psi(\blold{v_3}) \le (3c_dT)^{A}/A!$. There are $(|\blold{v}|+1)(|\blold{v}|+2)/2$ decompositions $\blold{v_1 v_2 v_3} = \blold{v}$, which concludes the proof.
\end{proof}

\begin{remark}\label{rem:NX15}
    For Brownian words, that is, words without the letter $\blold{0}$, \citet*{NiXu2015} shows that $|\E[\WH^{\blold{v}}_T - \WT^{\blold{v}}_T(N)]| = \Theta(N^{-1})$ when $|\blold{v}| \in \{4, 6, \dots\}$, while the error vanishes in odd length and in length two. Proposition~\ref{prop:weak-element-discretisation} gives the same order as an upper bound for all words, uniformly in time, together with the vanishing in odd degree.
\end{remark}

As in Section~\ref{sec:strong-discretisation}, we pass to linear functionals.

\begin{lemma}[Weak discretisation error of linear functionals]\label{lem:weak-discretisation}
    Under the assumptions of Lemma~\ref{lem:Lp-discretisation}, there exists $C > 0$, depending only on $d$, $T$, $\varepsilon$ and $C_*$, such that
    \[
        \Big| \E \brackets{\lbold(M,t)}{\WH_t} - \E\brackets{\lbold(M,t)}{\WT_t(N)} \Big| \le C\, N^{-1}, \qquad t \in [0,T],\ M \in \Nbar,\ N \in \N^* .
    \]
\end{lemma}

\begin{proof}
    By Corollary~\ref{cor:well-posed}, applied with $N = \infty$ and with $N$, $\sum_{\blold{v}} \lbold^{\ast,\blold{v}}\, \E|\Delta^{\blold{v}}(t)| < \infty$, so Fubini's theorem gives $\E\textbrackets{\lbold(M,t)}{\WH_t} - \E\textbrackets{\lbold(M,t)}{\WT_t(N)} = \sum_{\blold{v}} \lbold^{\blold{v}}(M,t)\, \E[\Delta^{\blold{v}}(t)]$. By Proposition~\ref{prop:weak-element-discretisation}, only words of even degree $2k$ contribute. For them, $\lbold^{\ast,\blold{v}} \le C_* ((2k)!)^{\frac12-\varepsilon}$, $|\blold{v}| \le 2k$, and $\sqrt{(2k)!} \le 2^k k!$, since $\binom{2k}{k} \le 4^k$. As there are at most $(d+1)^{2k}$ words of degree $2k$ (Lemma~\ref{lem:count-degree}), 
    \[
        \sum_{\blold{v}} \lbold^{\ast,\blold{v}}\, \big|\E[\Delta^{\blold{v}}(t)]\big| \le 2C_* \Big( \sum_{k \ge 0} (2k+1)(2k+2)\, \big(6(d+1)^2 c_d T\big)^k \big((2k)!\big)^{-\varepsilon} \Big) N^{-1},
    \]
    the series being finite since $\varepsilon > 0$.
\end{proof}

\begin{proof}[Proof of Theorem~\ref{th:weak-discretisation-bound}]
    Fix $M \in \Nbar$, $N \in \N^*$ and $t \in [0,T]$, and write $Y \coloneqq X_t(M,\infty)$ and $Y^N \coloneqq X_t(M,N)$. \emph{Polynomial test functions.} By linearity, it suffices to treat a monomial $z \mapsto z^\alpha$. With $\lbold^\alpha(M,t) \preceq \lbold^{\alpha,*}$ as in the proof of Theorem~\ref{th:Lp-discretisation-bound}, the shuffle identity (Lemma~\ref{lem:shuffle-L1}) gives $\E[Y^\alpha] - \E[(Y^N)^\alpha] = \E\textbrackets{\lbold^\alpha(M,t)}{\WH_t} - \E\textbrackets{\lbold^\alpha(M,t)}{\WT_t(N)}$, and Lemma~\ref{lem:weak-discretisation} bounds it by $C N^{-1}$. \emph{Hölder test functions.} Since $|\E[g(Y) - g(Y^N)]| \le \|g(Y) - g(Y^N)\|_1$, Theorem~\ref{th:Lp-discretisation-bound} with $p = 1$ concludes.
\end{proof}

\begin{remark}[Proof of the optimality of the discretisation rates]\label{rem:discretisation-optimality}
    We give the examples behind Remark~\ref{rem:optimality-statement}. The rates of Propositions~\ref{prop:Lp-element-discretisation} and~\ref{prop:weak-element-discretisation} and of Proposition~\ref{prop:Lp-discretisation-d1} cannot be improved, and single words suffice to see it. In particular, Theorems~\ref{th:Lp-discretisation-bound} and~\ref{th:weak-discretisation-bound} are sharp for polynomial test functions, for every $M \ge 2$.
    \begin{itemize}
        \item \emph{Off the grid.} Let $t = t_k + \lambda h$ with $\lambda \in (0,1)$. By the proof of Proposition~\ref{prop:truncation-sharpness}, $\WH^{\blold{1}}_t - \WT^{\blold{1}}_t$ is centred Gaussian with variance $\lambda(1-\lambda)h$, and $\E\big[\WH^{\blold{11}}_t - \WT^{\blold{11}}_t\big] = \tfrac12\big(t - \sigma_N(t)^2\big) = \tfrac12\lambda(1-\lambda)h$. At midpoints ($\lambda = 1/2$), this gives a strong error of exact order $N^{-1/2}$ and a weak error of exact order $N^{-1}$, uniformly in $t$. In particular, the restriction to grid points in Proposition~\ref{prop:Lp-discretisation-d1} cannot be removed.
        \item \emph{At the grid points, $d \ge 2$.} The interpolation reproduces the Brownian path, and only higher-order components carry an error. The difference $\WH^{\blold{12}}_T - \WT^{\blold{12}}_T$ is the sum of the Lévy areas $A_{12}(t_{k-1}, t_k)$, and no approximation built on the grid increments alone can approximate it faster than $N^{-1/2}$ in $L^2$ \citep*{ClarkCameron1980}. This is the same obstruction that limits the Euler scheme to strong order $1/2$ for non-commutative noise \citep*[Chapter 10]{KloedenPlaten1992}. The weak rate $N^{-1}$ is attained for Brownian words of even length at least four \citep*{NiXu2015}, and it matches the weak order one of the Euler scheme \citep*[Chapter 14]{KloedenPlaten1992}.
        \item \emph{At the grid points, $d = 1$.} There is no Lévy area between Brownian components. At $t = T$, Chen's identity and the exactness of the level-one increments give $\WH^{\blold{01}}_T - \WT^{\blold{01}}_T = \sum_{k=1}^N Z_k$, with $Z_k$ as in the proof of Proposition~\ref{prop:Lp-discretisation-d1}, independent, centred and Gaussian with variance $h^3/12$. Hence $\textnorm{\WH^{\blold{01}}_T - \WT^{\blold{01}}_T}_2 = T^{3/2}/(\sqrt{12}\,N)$, and the rate $N^{-1}$ of Proposition~\ref{prop:Lp-discretisation-d1} and Theorem~\ref{th:Lp-discretisation-bound} is sharp.
    \end{itemize}
    For Hölder test functions with $a < 1$, whether the exponent $a$ is sharp is left open, as for the truncation error (Remark~\ref{rem:sharpness-comments}(iii)).
\end{remark}

\section{Proofs of the extensions}\label{sec:extensions}

All the statements of Section~\ref{sec:extensions-statements} follow from the strong estimate of Corollary~\ref{cor:total-error}, applied at a suitable order $p$. We fix $t$, $M$ and $N$, and write $Z_1 \coloneqq g(X_t)$, $Z_2 \coloneqq g(X_t(M,N))$ and $\varrho \coloneqq \norm{Z_1 - Z_2}_p \le C\,\mathcal{E}_a(M,N)$, where $p$ is specified in each proof. Since $(r_M)_M$ is bounded, we have $R \coloneqq \sup_{t,M,N} \varrho < \infty$. Since $\WH_0 = \WT_0(N) = \empt$, we have $X_0(M,N) = \lbold_0^{\empt} = X_0$ for all $M \in \Nbar$, $N \in \Nbar^*$, so all the errors of Section~\ref{sec:extensions-statements} vanish at $t=0$. In the proofs of Lemma~\ref{lem:cdf-bound} and Theorem~\ref{th:var-bound}(ii) we may therefore assume $t\in(0,T]$.

\begin{proof}[Proof of Lemma~\ref{lem:probability-bound}]
    Markov's inequality gives $\Prob(|Z_1 - Z_2| \ge \eta) \le \eta^{-p}\varrho^p$.
\end{proof}

\begin{proof}[Proof of Lemma~\ref{lem:cdf-bound}]
    Fix $p \ge 1$ with $p/(p+1) \ge \theta$. Let $x \in J$ and $\eta \in (0,\eta_0)$, so that $[x-\eta, x+\eta] \subset I$. Lemma~\ref{lem:cdf-comparison} and Markov's inequality give
    \begin{equation}\label{eq:cdf-split}
        \big|F^g_t(x) - F^{g,M,N}_t(x)\big| \le 2\eta\,\|f^g\|_I + \eta^{-p}\varrho^p .
    \end{equation}
    If $\varrho \le 1$ and $\varrho^{p/(p+1)} < \eta_0$, the choice $\eta \coloneqq \varrho^{p/(p+1)}$ (or any $\eta \to 0$ if $\varrho = 0$) yields the bound $(2\norm{f^g}_I + 1)\varrho^{p/(p+1)} \le (2\norm{f^g}_I + 1)\varrho^\theta$. Otherwise, either $\varrho > 1$ and $1 \le \varrho^\theta$, or $\varrho \le 1$ and $1 \le \eta_0^{-1}\varrho^{p/(p+1)} \le \eta_0^{-1}\varrho^\theta$; since distribution functions differ by at most $1$, the left-hand side of~\eqref{eq:cdf-split} is then at most $(1 \vee \eta_0^{-1})\varrho^\theta$. In all cases, it is at most $(2\norm{f^g}_I + 1 + \eta_0^{-1})\varrho^\theta \le C\,\mathcal{E}_a(M,N)^\theta \le C\,\mathcal{E}_{\theta a}(M,N)$.
\end{proof}

\begin{lemma}[Hölder bound for coherent risk measures] \label{lem:coherent-Lp-control}
    Assume that there exist $q \in (1, \infty]$ and $C > 0$ such that $\norm{\diff \Q / \diff \Prob}_q \le C$ for every $\Q \in \mathcal{Q}_\rho$, and let $q^* \coloneqq q/(q-1) \in [1, \infty)$ be the conjugate exponent. Then, for all $Z, Z_1, Z_2 \in L^{q^*}$,
    \[
        |\rho|[Z] \le C \norm{Z}_{q^*}
        \qquad \text{and} \qquad
        |\rho(Z_1) - \rho(Z_2)| \le |\rho|[Z_1 - Z_2] \le C \norm{Z_1 - Z_2}_{q^*}.
    \]
\end{lemma}

\begin{proof}
    For $\Q \in \mathcal{Q}_\rho$, Hölder's inequality gives $|\E^{\Q}[Z]| = |\E[Z \, \diff\Q / \diff\Prob]| \le C \norm{Z}_{q^*}$; taking the supremum over $\mathcal{Q}_\rho$, and applying the same to $-Z$, yields the first bound. By sub-additivity, $\rho(Z_1) \le \rho(Z_2) + \rho(Z_1 - Z_2)$ and $\rho(Z_2) \le \rho(Z_1) + \rho(Z_2 - Z_1)$, so $|\rho(Z_1) - \rho(Z_2)| \le |\rho|[Z_1 - Z_2]$.
\end{proof}

\begin{proof}[Proof of Theorem~\ref{th:coherent-bound}]
    Both $Z_1$ and $Z_2$ belong to $L^{q^*}$: for polynomial $g$ by Corollary~\ref{cor:well-posed}(ii), and for Hölder $g$ since $|g(z)| \le |g(0)| + K_g\tnorm{z}^a$. Lemma~\ref{lem:coherent-Lp-control} and Corollary~\ref{cor:total-error} with $p = q^*$ conclude.
\end{proof}

\begin{lemma}[Value-at-Risk and Expected Shortfall]\label{lem:VaR-dominated-by-ES}
    Let $\alpha \in (1/2, 1)$. Then $|\VaR_\alpha|[Z] \le |\ES_\alpha|[Z]$ for every integrable random variable $Z$.
\end{lemma}

\begin{proof}
    Since $u \mapsto \VaR_u[Z]$ is non-decreasing, $\VaR_\alpha[Z] \le \ES_\alpha[Z]$. Writing $\VaR^+_u[Z] \coloneqq \inf\{x \in \R : \Prob(Z \le x) > u\}$ for the upper quantile, the reflection identity $\VaR_\alpha[Z] = - \VaR^+_{1-\alpha}[-Z]$ and the inequality $\VaR^+_u \le \VaR_v$ for $u < v$ give, since $1 - \alpha < \alpha$,
    \[
        -\VaR_\alpha[Z] = \VaR^+_{1-\alpha}[-Z] \le \VaR_\alpha[-Z] \le \ES_\alpha[-Z].
    \]
    Hence $|\VaR_\alpha[Z]| \le \ES_\alpha[Z] \vee \ES_\alpha[-Z] \le |\ES_\alpha|[Z]$, and the same applied to $-Z$ yields the claim.
\end{proof}

\begin{lemma}[Stability of quantiles]\label{lem:VaR-stability}
    Let $Z_1, Z_2$ be real random variables and $\alpha \in (0,1)$. Assume that there exist $\gamma > 0$, an interval $\mathcal{B} \subset \R$ and $c > 0$ such that $\alpha \pm \gamma \in (0,1)$, $\VaR_{\alpha - \gamma}[Z_1], \VaR_{\alpha + \gamma}[Z_1] \in \mathcal{B}$, and the law of $Z_1$ admits a density $f \ge c$ on $\mathcal{B}$. Then, for every $\eta > 0$ such that $\nu \coloneqq \Prob(|Z_1 - Z_2| \ge \eta)$ satisfies $2\nu \le \gamma$,
    \[
        \big|\VaR_\alpha[Z_1] - \VaR_\alpha[Z_2]\big| \le \eta + \frac{2\nu}{c}.
    \]
\end{lemma}

\begin{proof}
    Splitting $\{Z_2 \le x\}$ along $\{|Z_1 - Z_2| < \eta\}$ and its complement, and likewise for $Z_1$, gives
    \begin{equation} \label{eq:VaR-level-shift}
        \Prob(Z_2 \le x) \le \Prob(Z_1 \le x + \eta) + \nu, \qquad \Prob(Z_1 \le x - \eta) \le \Prob(Z_2 \le x) + \nu, \qquad x \in \R.
    \end{equation}
    Since $2\nu \le \gamma$ and $\mathcal{B}$ is an interval, $\VaR_u[Z_1] \in \mathcal{B}$ for $|u - \alpha| \le 2\nu$, and the density of $Z_1$ on $\mathcal{B}$ gives $\Prob(Z_1 \le \VaR_u[Z_1]) = u$ there. Let $x_0 := \VaR_{\alpha-2\nu}[Z_1] - \eta$ and $x < x_0$. Since $x + \eta < \VaR_{\alpha - 2\nu}[Z_1]$, the definition of the quantile gives $\Prob(Z_1 \le x + \eta) < \alpha - 2\nu$, and the first inequality in \eqref{eq:VaR-level-shift} yields $\Prob(Z_2 \le x) < \alpha - \nu \le \alpha$. As this holds for every $x < x_0$, $\VaR_\alpha[Z_2] \ge \VaR_{\alpha - 2\nu}[Z_1] - \eta$. Taking $x \coloneqq \VaR_{\alpha+\nu}[Z_1] + \eta$ in the second inequality of~\eqref{eq:VaR-level-shift} gives $\Prob(Z_2 \le x) \ge \Prob(Z_1 \le \VaR_{\alpha + \nu}[Z_1]) - \nu = \alpha$, hence $\VaR_\alpha[Z_2] \le \VaR_{\alpha+\nu}[Z_1] + \eta$. Finally, $f \ge c$ on $\mathcal{B}$ makes the quantile function of $Z_1$ Lipschitz with constant $1/c$ on $[\alpha - \gamma, \alpha + \gamma]$, so that $|\VaR_\alpha[Z_1] - \VaR_u[Z_1]| \le 2\nu/c$ for $u \in \{\alpha - 2\nu, \alpha + \nu\}$, which concludes.
\end{proof}

\begin{proof}[Proof of Theorem~\ref{th:var-bound}]
    \emph{(i)} By Lemma~\ref{lem:VaR-dominated-by-ES}, $|\VaR_\alpha|[Z_1 - Z_2] \le |\ES_\alpha|[Z_1 - Z_2]$, and Theorem~\ref{th:coherent-bound} applies to $\ES_\alpha$ with $q = \infty$. \emph{(ii)} Fix $p \ge 1$ with $p/(p+1) \ge \theta$. If $\varrho = 0$, then $Z_1 = Z_2$ almost surely and there is nothing to prove. Otherwise, set $\eta \coloneqq \lambda\varrho^{p/(p+1)}$ for some $\lambda > 0$. By Markov's inequality,
    \[
        \nu \coloneqq \Prob(|Z_1 - Z_2| \ge \eta) \le \eta^{-p}\varrho^p = \lambda^{-p}\varrho^{p/(p+1)} \le \lambda^{-p}R^{p/(p+1)},
    \]
    so that the choice $\lambda \coloneqq (2R^{p/(p+1)}/\gamma)^{1/p}$ ensures $2\nu \le \gamma$ for every $t$, $M$ and $N$. Lemma~\ref{lem:VaR-stability} then gives
    \[
        \big|\VaR_\alpha[Z_1] - \VaR_\alpha[Z_2]\big| \le \Big(\lambda + \frac{2}{c\,\lambda^p}\Big)\varrho^{p/(p+1)} \le \Big(\lambda + \frac{2}{c\,\lambda^p}\Big)(1 \vee R)\,\varrho^\theta,
    \]
    where we used $\varrho^{p/(p+1)} \le \varrho^\theta$ if $\varrho \le 1$ and $\varrho^{p/(p+1)} \le R \le R\varrho^\theta$ otherwise. Corollary~\ref{cor:total-error} and $\mathcal{E}_a^\theta \le \mathcal{E}_{\theta a}$ conclude.
\end{proof}
\section{Proof of the pricing application} \label{sec:pricing-proof}

We prove the claims of Section~\ref{sec:pricing} in a slightly more general model. Fix $\phibold \in \mathcal{L}^1$ and set $\sigma_t \coloneqq \textbrackets{\phibold}{\WH_t}$ and $\sigma_t(M) \coloneqq \textbrackets{\phibold^{\le M}}{\WH_t}$. For $\alpha \in \R^{d+1}$, let $Y \coloneqq \sum_{i=0}^d \alpha_i \Wh^i$, so that $[Y]_t = \alpha^\star t$ with $\alpha^\star \coloneqq \sum_{i=1}^d \alpha_i^2$, and let $S$ and $S(M)$ solve
\begin{subnumcases}{\label{eq:pricing-sde}}
    \diff S_t = S_t \, \sigma_t \, \diff Y_t, \label{eq:pricing-sde-exact}\\
    \diff S_t(M) = S_t(M) \, \sigma_t(M) \, \diff Y_t. \label{eq:pricing-sde-trunc}
\end{subnumcases}
with $S_0(M) = S_0 > 0$. The model of Section~\ref{sec:pricing} is the case $d = 2$, $\alpha = (0, \rho, \sqrt{1-\rho^2})$, so that $Y = W^S$, and $\phibold = \lbold$.

\paragraph{Itō integrals in tensor form}
Define, for $i \in \intset{0}{d}$ and $\lbold \in T((\R^{d+1}))$, 
\[
    \pi_i(\lbold) \coloneqq \lbold \blold{i} - \frac{1}{2}\indic{i \ne 0} \, (\lbold|_{\blold{i}}) \blold{0}.
\]

\begin{lemma}[Itō integrals in tensor form] \label{lem:Ito-integral-signature}
    For every $\lbold \in \mathcal{L}^1$, $\pi_i(\lbold) \in \mathcal{L}^1$ and
    \[
        \int_0^t \brackets{\lbold}{\WH_s} \, \diff \Wh^i_s = \brackets{\pi_i(\lbold)}{\WH_t}, \qquad t \in [0,T].
    \]
\end{lemma}

\begin{proof}
    For a word $\blold{v}$, the Stratonovich integral of $\WH^{\blold{v}}$ against $\Wh^i$ is $\WH^{\blold{vi}}$ by \eqref{eq:signature}, and the Itō correction is $-\frac12[\WH^{\blold{v}}, \Wh^i]$. This bracket vanishes if $\blold{v} = \empt$ or $i = 0$; otherwise, writing $\blold{v} = \blold{uj}$, only the martingale part of $\WH^{\blold{v}}$ contributes, and $[\WH^{\blold{v}}, \Wh^i]_t = \indic{i=j} \WH^{\blold{u0}}_t = \textbrackets{(\blold{v}|_{\blold{i}})  \blold{0}}{\WH_t}$. This proves the identity for words, hence for $\lbold \in T(\R^{d+1})$ by linearity. For $\lbold \in \mathcal{L}^1$, the partial sums $Y^M_t \coloneqq \textbrackets{\lbold^{\le M}}{\WH_t}$ are continuous and adapted, and $Y^M_t \to \textbrackets{\lbold}{\WH_t}$ almost surely for every $t$ (Corollary~\ref{cor:well-posed}(i)). Hence $\bar{Y}_t \coloneqq \limsup_{M} Y^M_t$ is a progressively measurable version of $t \mapsto \textbrackets{\lbold}{\WH_t}$, with respect to which the stochastic integral is defined. By the Itō isometry and Corollary~\ref{cor:well-posed}(iii), $\int_0^t Y^M_s\,\diff \Wh^i_s \to \int_0^t \bar{Y}_s\,\diff \Wh^i_s$ in $L^2$. On the other side, $\pi_i(\lbold) - \pi_i(\lbold^{\le M}) = \lbold^{>M}\blold{i} - \frac12\indic{i \ne 0} (\lbold^{>M}|_{\blold{i}})\blold{0}$ by linearity, so that the first term has words of length at least $M+2$, and in the second, the words of $\lbold^{>M}$ ending in $\blold{i}$ have length at least $M+1$, which is preserved by deleting $\blold{i}$ and appending $\blold{0}$. The difference is then supported on words of length larger than $M$ and dominated by $\pi_i^+(|\lbold|) \in \mathcal{L}^1$ where $\pi_i^+$ is defined as $\pi_i$ with the minus sign replaced by a plus sign, so that $\textbrackets{\pi_i(\lbold^{\le M})}{\WH_t} \to \textbrackets{\pi_i(\lbold)}{\WH_t}$ in $L^2$ by Lemma~\ref{lem:truncation}.
\end{proof}

\paragraph{The log-price is a linear functional of the signature}
By Itō's formula applied to \eqref{eq:pricing-sde-exact}, $X \coloneqq \log S$ satisfies $X_t = \log S_0 - \frac{\alpha^\star}{2}\int_0^t \sigma_s^2 \diff s + \int_0^t \sigma_s \diff Y_s$. Lemma~\ref{lem:Ito-integral-signature} handles the stochastic integral. For the Lebesgue one, $\sigma_t^2 = \textbrackets{\phibold^{\shuffle 2}}{\WH_t}$ almost surely for every $t$, by the shuffle identity on $\mathcal{L}^1$
(Lemma~\ref{lem:shuffle-L1}), hence $\diff t \otimes \diff \Prob$-almost everywhere by Fubini's theorem, and integrating in time amounts to concatenating $\blold{0}$. Hence
\begin{equation} \label{eq:log-price-representation}
    X_t = \brackets{\Phi(\phibold)}{\WH_t}, \qquad \Phi(\lbold) \coloneqq \log(S_0)\, \empt - \frac{\alpha^\star}{2}\, \lbold^{\shuffle 2} \blold{0} + \sum_{i=0}^d \alpha_i \pi_i(\lbold),
\end{equation}
and $\Phi$ maps $\mathcal{L}^1$ into itself by Proposition~\ref{prop:L1-stability}.

\paragraph{Error bounds}
The key observation is that $\Phi$ commutes with truncation.

\begin{lemma}[Truncation commutes with $\Phi$]\label{lem:phi-trunc}
    For every $\lbold \in T((\R^{d+1}))$ and $M \in \N$, the tensor $\Phi(\lbold) - \Phi(\lbold^{\le M})$ is supported on words of length larger than $M$. In particular, $\Phi(\lbold)^{\le M} = \Phi(\lbold^{\le M})^{\le M}$.
\end{lemma}

\begin{proof}
    Write $\lbold = \lbold^{\le M} + \lbold^{>M}$. The constant term cancels. By the arguments of the proof of Lemma~\ref{lem:Ito-integral-signature}, $\pi_i(\lbold) - \pi_i(\lbold^{\le M})$ have words of length at least $M+1$. Finally, $\lbold^{\shuffle 2} - (\lbold^{\le M})^{\shuffle 2} = 2\,\lbold^{\le M} \shuffle \lbold^{>M} + (\lbold^{>M})^{\shuffle 2}$, and the shuffle product adds lengths, so this difference, and a fortiori its product with $\blold{0}$, has words of length at least $M+1$.
\end{proof}

Fix $M \in \N$ and set $\Psi_M \coloneqq \Phi(\phibold^{\le M})$, so that $\check{X}_t \coloneqq \textbrackets{\Psi_M}{\WH_t} = \log S_t(M)$ is the log-price of the model with truncated volatility~\eqref{eq:pricing-sde-trunc}. By Lemma~\ref{lem:phi-trunc}, the scheme of Section~\ref{sec:pricing} is
\[
  \check{X}_t(M,N) \coloneqq \brackets{\Psi_M^{\le M}}{\WT_t(N)} = \brackets{\Phi(\phibold)^{\le M}}{\WT_t(N)} = X_t(M,N) \quad \big(\check{X}_t(\infty, N) \coloneqq X_t(\infty, N) \big),
\]
the computable approximation of Section~\ref{sec:main-results} for the fixed coefficient family $\Phi(\phibold)$.

\begin{proposition}[Error bounds for the pricing scheme]\label{prop:pricing}
    Let $\phibold \in \mathcal{L}^1$ with summability exponent $\varepsilon$, let $\varepsilon' \in (0,\varepsilon)$, let $g \in \mathcal{G}^1$ with regularity exponent $a$, and let $p \ge 1$. Write $r'_M$ and $\mathcal{E}'_a(M,N)$ for $r_M$ and $\mathcal{E}_a(M,N)$ with $\varepsilon$ replaced by $\varepsilon'$. There exist $B, C > 0$, independent of $t$, $M$ and $N$, such that:
    \begin{enumerate}[label=(\roman*)]
        \item for $t \in [0,T]$, $M \in \Nbar$ and $N \in \Nbar^*$,
        \[
            \norm{g(X_t) - g\big(\check{X}_t(M,N)\big)}_p \le C\, \mathcal{E}'_a(M,N),
        \]
        and, if $g$ is polynomial, $\big|\E[g(X_t) - g(\check{X}_t(M,N))]\big| \le C\big(r'_M + N^{-1}\big)$;
        \item if $g$ is $a$-Hölder, the model with truncated volatility and its discretisation satisfy, for $t \in [0,T]$, $M \in \N$ and $N \in \Nbar^*$,
        \[
            \norm{g(X_t) - g\big( \brackets{\Psi_M}{\WT_t(N)} \big)}_p \le C\, \mathcal{E}'_a(M,N);
        \]
        for $N = \infty$, this is the error $\norm{g(\log S_t) - g(\log S_t(M))}_p$.
    \end{enumerate}
\end{proposition}

\begin{proof}
    Let $\Phi^+$ be defined as $\Phi$ with every coefficient replaced by its absolute value and the minus signs by plus signs. By Proposition~\ref{prop:L1-stability}, $\Phi$ and $\Phi^+$ map $\mathcal{L}^1$ into itself, and $\Phi(\phibold), \Phi^+(|\phibold|)$ admit $\varepsilon'$ as a summability exponent. \emph{(i)} By Lemma~\ref{lem:phi-trunc}, $\check X_t(M,N) = X_t(M,N)$ for the coefficient family $\Phi(\phibold)$, and Corollaries~\ref{cor:total-error} and~\ref{cor:total-error-weak} apply with exponent $\varepsilon'$. \emph{(ii)} Fix $N \in \Nbar^*$. By Lemma~\ref{lem:phi-trunc}, the family $\lbold(M) \coloneqq \Psi_M - \Phi(\phibold)$ is supported on words of length larger than $M$, and $\lbold(M) \preceq 2\Phi^+(|\phibold|)$ uniformly in $M$. Hence $\textbrackets{\Psi_M}{\WT_t(N)} - X_t(\infty,N) = \textbrackets{\lbold(M)}{\WT^{>M}_t(N)}$, and Lemma~\ref{lem:truncation} with $q = 1$ bounds its $L^p$ norm by $C r'_M$, uniformly in $t$ and $N$. Moreover, $\norm{X_t(\infty,N) - X_t}_p \le C N^{-1/2}$ by Theorem~\ref{th:Lp-discretisation-bound} with $g = \mathrm{Id}$ and $M = \infty$. The triangle inequality and $\norm{g(Z_1) - g(Z_2)}_p \le K_g \norm{Z_1 - Z_2}_p^a$ conclude, the power $a$ being distributed over the two terms by subadditivity.
\end{proof}

\begin{remark}[Why the logarithm] \label{rem:why-log}
    A representation $S = \textbrackets{\chi}{\WH}$ would require $\chi$ to solve the triangular system $\chi = S_0\, \empt + \sum_i \alpha_i \pi_i(\chi \shuffle \phibold)$, which is numerically solvable but for which we do not know non-trivial $\phibold \in \mathcal{L}^1$ ensuring $\chi \in \mathcal{L}^1$. The logarithm avoids the question, $\Phi$ mapping $\mathcal{L}^1$ into itself by construction, and requires no equation to be solved.
\end{remark}

\appendix

\section{General results on word algebra} \label{app:general-results}

\begin{lemma} \label{lem:sup-phi}
    Let $n \in \N^*$ and take $\varepsilon_1, \dots, \varepsilon_n \in (0, 1/2]$, then, for every $\blold{v} \in V$,
    \[
        \sup_{\substack{\blold{v_1}, \dots, \blold{v_n} \in V \text{ s.t.} \\ \textbrackets{\blold{v}}{\bigshuffle_{i = 1}^n \blold{v_i}} \ne 0}} \,\, \prod_{i=1}^n \varphi_{\varepsilon_i}(\blold{v_i}) \le \varphi_{\min_{1 \le i \le n} \varepsilon_i}(\blold{v}).
    \]
\end{lemma}

\begin{proof}
    Let $\varepsilon_\wedge \coloneqq \min_{1 \le i \le n} \varepsilon_i$. Since $\varphi_\varepsilon(\blold{v})$ is non-increasing in $\varepsilon$ and $\prod_{i=1}^n \norm{\blold{v_i}}! \le \norm{\blold{v}}!$,
    \[
        \prod_{i=1}^n \varphi_{\varepsilon_i}(\blold{v_i}) \le \prod_{i=1}^n \varphi_{\varepsilon_\wedge}(\blold{v_i}) = \left[ \prod_{i=1}^n \norm{\blold{v_i}}! \right]^{\frac{1}{2} - \varepsilon_\wedge} \le \varphi_{\varepsilon_\wedge}(\blold{v}).
    \]
\end{proof}

\begin{lemma} \label{lem:count-length-m-k-zeros}
    Let $k \le m$ be two non-negative integers. Then
    \[
        \#\{ \blold{v} \in V : |\blold{v}| = m, \ n_0(\blold{v}) = k \} = \binom{m}{k}d^{m-k}.
    \]
\end{lemma}

\begin{proof}
    A word $\blold{v}$ of length $m$ with exactly $k$ zeros is uniquely determined by:
    \begin{itemize}[noitemsep, topsep=0pt]
        \item the choice of the $k$ positions among $\intset{1}{m}$ where the letter $\blold{0}$ appears;
        \item the choice of the remaining $m-k$ letters, each taken from $\intset{1}{d}$.
    \end{itemize}
    There are $\binom{m}{k}$ ways to choose the positions of the zeros, and $d^{m-k}$ ways to fill the remaining positions. The result follows.
\end{proof}

\begin{lemma}\label{lem:count-degree}
    For every $n \in \N$,
    \[
        \#\{\blold{v} \in V : \norm{\blold{v}} = n\} \le (d+1)^n .
    \]
\end{lemma}

\begin{proof}
    Write $a_n$ for the left-hand side. The empty word is the only word of degree $0$, and the words of degree $1$ are the letters $\blold{1}, \dots, \blold{d}$, so that $a_0 = 1$ and $a_1 = d$. For $n \ge 2$, a word of degree $n$ ends either with a letter $\blold{i} \in \{\blold{1}, \dots, \blold{d}\}$, which contributes $1$ to the degree, or with the letter $\blold{0}$, which contributes $2$. Removing the last letter is therefore a bijection onto $d$ copies of the words of degree $n-1$ and one copy of the words of degree $n-2$, so that $a_n = d\, a_{n-1} + a_{n-2}$. The bound holds for $n \in \{0,1\}$, and if it holds up to $n-1$, then $a_n \le d(d+1)^{n-1} + (d+1)^{n-2} \le (d+1)^{n-1}(d+1) = (d+1)^n$.
\end{proof}

\begin{lemma} \label{lem:count-v-shuffle_p-elements}
    Let $\blold{v} \in V$ and $p \in \N^*$. Then
    \[
        \sum_{\blold{u}} \brackets{\blold{u}}{\blold{v}^{\shuffle p}} = \frac{(p |\blold{v}|)!}{(|\blold{v}|!)^p} \le p^{p|\blold{v}|}.
    \]
\end{lemma}

\begin{proof}
    The quantity $\sum_{\blold{u}} \brackets{\blold{u}}{\blold{v}^{\shuffle p}}$ counts the number of ways of interleaving $p$ copies of the word $\blold{v}$. An interleaving is uniquely determined by choosing, for each of the $p$ copies, the positions of its $|\blold{v}|$ letters among the $p|\blold{v}|$ available positions, in such a way that the order of letters within each copy is preserved. Thus, we partition the set $\intset{1}{p|\blold{v}|}$ into $p$ ordered subsets of size $|\blold{v}|$. The number of such choices is the multinomial coefficient
    \[
        \binom{p|\blold{v}|}{\underbrace{|\blold{v}|, \dots, |\blold{v}|}_{p \text{ times}}} = \frac{(p|\blold{v}|)!}{(|\blold{v}|!)^p},
    \]
    where $\binom{n}{k_1, \dots, k_m} = \frac{n!}{k_1! \cdots k_m!}$. Bounding it by the sum of all multinomial coefficients
    \[
        \frac{(p|\blold{v}|)!}{(|\blold{v}|!)^p} \le \sum_{j_1 + \dots + j_p = p|\blold{v}|} \binom{p|\blold{v}|}{j_1, \dots, j_p} = p^{p|\blold{v}|}
    \]
    concludes the proof.
\end{proof}

\begin{lemma} \label{lem:count-shuffle-decompositions}
    Let $\blold{v} \in V$ and $n \in \N^*$. Then
    \[
        \sum_{\blold{v_1}, \dots, \blold{v_n}} \brackets{\blold{v}}{\bigshuffle_{k=1}^n \blold{v_k}} = n^{|\blold{v}|}.
    \]
\end{lemma}

\begin{proof}
    The quantity $\sum_{\blold{v_1}, \dots, \blold{v_n}} \brackets{\blold{v}}{\bigshuffle_{k=1}^n \blold{v_k}}$ counts the number of ways to obtain $\blold{v}$ as a shuffle of $n$ words. Constructing such a decomposition amounts to assigning each position of $\blold{v}$ to one of the $n$ words $\blold{v_1}, \dots, \blold{v_n}$, while preserving the order of letters within each word. Thus, for each index $i \in \intset{1}{|\blold{v}|}$, we choose an integer in $\intset{1}{n}$ indicating to which word the $i$-th letter of $\blold{v}$ belongs. This completely determines the words $\blold{v_1}, \dots, \blold{v_n}$. Since there are $n$ choices for each of the $|\blold{v}|$ positions, the total number of such decompositions is $n^{|\blold{v}|}$.
\end{proof}

\section{Additional proofs}

\subsection{Proof of Proposition~\ref{prop:expected-signature-bound}} \label{app:expected-signature-bound}

\begin{proof}[Proof of Proposition~\ref{prop:expected-signature-bound}]
    Fix $N \in \N^*$, which we omit, and $0 \le s \le t \le T$. The case $s = t$ being immediate, we assume $s < t$, so in particular $s < T$. Choose $\ndown \le \nup$ in $\{1,\dots,N\}$ such that $s \in [t_{\ndown-1}, t_{\ndown})$ and $t \in (t_{\nup-1}, t_{\nup}]$, and set $\tilde t_i \coloneqq s$ for $i < \ndown$, $\tilde t_i \coloneqq t_i$ for $\ndown \le i < \nup$, and $\tilde t_i \coloneqq t$ for $i \ge \nup$. For $\ndown \le k \le \nup$, write
    \[
        r_k \coloneqq \frac{\tilde t_k - \tilde t_{k-1}}{t_k - t_{k-1}} \in [0,1],
        \qquad
        \hat{y}_k \coloneqq r_k\big(\Wh_{t_k} - \Wh_{t_{k-1}}\big) \in \R^{d+1},
    \]
    so that $\hat y_k^0 = \tilde t_k - \tilde t_{k-1}$ and $\hat y_k^i = r_k (W^i_{t_k} - W^i_{t_{k-1}})$ for $1 \le i \le d$.

    On $[\tilde t_{k-1}, \tilde t_k]$ the path $\Wt$ is linear with increment $\hat{y}_k$, so its signature there is $\exp_\otimes(\hat{y}_k)$, and Chen's identity~\eqref{eq:chen-identity} gives $\WT_{s,t} = \bigotimes_{k=\ndown}^{\nup} \exp_\otimes(\hat{y}_k)$. Each coordinate of this tensor product is a finite sum of products of factors depending on distinct, hence independent, Brownian increments, so that
    \begin{equation} \label{eq:scheme-factorisation}
        \E\big[\WT_{s,t}\big] = \bigotimes_{k=\ndown}^{\nup} \boldsymbol{e}_k,
        \qquad
        \boldsymbol{e}_k \coloneqq \E\big[\exp_\otimes(\hat{y}_k)\big].
    \end{equation}
    For a word $\blold{u} \in V$, $\boldsymbol{e}_k^{\blold{u}} = \frac{1}{|\blold{u}|!}(\tilde t_k - \tilde t_{k-1})^{n_0(\blold{u})}\prod_{i=1}^d r_k^{n_i(\blold{u})}\E\big[(W^i_{t_k} - W^i_{t_{k-1}})^{n_i(\blold{u})}\big] \ge 0$, moments of centred Gaussian variables being non-negative (odd ones vanish). Tensor products of coordinatewise non-negative tensors are coordinatewise non-negative, so \eqref{eq:scheme-factorisation} yields $\E[\WT^{\blold{u}}_{s,t}] \ge 0$ for every $\blold{u} \in V$.

    For $\lbold \in T((\R^{d+1}))$ coordinatewise non-negative and $z \ge 0$, set $\Gamma_z(\lbold) \coloneqq \sum_{\blold{u}} \lbold^{\blold{u}} z^{\norm{\blold{u}}} \in [0,\infty]$. Since $\norm{\blold{u}\blold{u'}} = \norm{\blold{u}} + \norm{\blold{u'}}$, Tonelli's theorem gives $\Gamma_z(\lbold \otimes \lbold') = \Gamma_z(\lbold)\,\Gamma_z(\lbold')$. For a deterministic $\hat{y} \in \R^{d+1}$, expanding the tensor exponential and using $\norm{\blold{0}} = 2$, $\norm{\blold{i}} = 1$ for $i \ge 1$,
    \[
        \sum_{\blold{u} \in V} \exp_\otimes^{\blold{u}}(\hat{y})\, z^{\norm{\blold{u}}}
        = \sum_{l \ge 0} \frac{1}{l!}\Big(\hat y^0 z^2 + \sum_{i=1}^d \hat y^i z\Big)^l
        = \exp\Big(\hat y^0 z^2 + \sum_{i=1}^d \hat y^i z\Big),
    \]
    the series converging absolutely, with absolute sum $\exp(|\hat y^0| z^2 + \sum_i |\hat y^i| z)$. The latter is integrable when $\hat{y} = \hat{y}_k$, so Fubini's theorem and the Gaussian Laplace transform give
    \[
        \Gamma_z(\boldsymbol{e}_k)
        = \E\Big[\exp\Big((\tilde t_k - \tilde t_{k-1}) z^2 + r_k z \sum_{i=1}^d \big(W^i_{t_k} - W^i_{t_{k-1}}\big)\Big)\Big]
        = e^{a_k z^2},
        \qquad
        a_k \coloneqq (\tilde t_k - \tilde t_{k-1})\Big(1 + \frac d2\, r_k\Big).
    \]
    By \eqref{eq:scheme-factorisation} and multiplicativity, $\Gamma_z(\E[\WT_{s,t}]) = e^{A_{s,t} z^2}$ with $A_{s,t} \coloneqq \sum_{k=\ndown}^{\nup} a_k$.

    Grouping words by degree, a finite operation for each degree, $\Gamma_z(\E[\WT_{s,t}]) = \sum_{n \ge 0} c_n z^n$ with $c_n \coloneqq \sum_{\norm{\blold{u}} = n} \E[\WT^{\blold{u}}_{s,t}] \ge 0$. This power series is finite for every $z \ge 0$, hence has infinite radius, and coincides on $[0,\infty)$ with $\sum_{l \ge 0} A_{s,t}^l z^{2l}/l!$. By uniqueness of power-series coefficients, $c_{2l} = A_{s,t}^l / l!$ and $c_{2l+1} = 0$. Finally, $r_k \le 1$ and $\sum_k (\tilde t_k - \tilde t_{k-1}) = t - s$ give $A_{s,t} \le c_d (t-s)$, with equality when $s$ and $t$ are grid points, since then $r_k = 1$ for every $k$. Bounding $\E[\WT^{\blold{v}}_{s,t}]$, which is non-negative by \eqref{eq:scheme-factorisation}, by the sum $c_{\norm{\blold{v}}}$ of its degree concludes the proof.
\end{proof}

\subsection{Proof of Proposition~\ref{prop:L1-stability}} \label{app:L1-closeness}

\begin{proof}[Proof of Proposition~\ref{prop:L1-stability}]
    Take $\varepsilon, \varepsilon_1, \varepsilon_2 \in (0, 1/2]$ such that
    \[
        |\lbold^{\blold{v}}| \lesssim \varphi_{\varepsilon}(\blold{v}), \quad |\lbold_1^{\blold{v}}| \lesssim \varphi_{\varepsilon_1}(\blold{v}), \quad |\lbold_2^{\blold{v}}| \lesssim \varphi_{\varepsilon_2}(\blold{v}), \qquad \blold{v} \in V.
    \]
    Linearity is immediate, and the tensor case is analogous to the shuffle. For the shuffle,
    \[
        |(\lbold_1 \shuffle \lbold_2)^{\blold{v}}| \lesssim \sum_{\blold{v_1}, \blold{v_2}} |\lbold_1^{\blold{v_1}}||\lbold_2^{\blold{v_2}}| \brackets{\blold{v}}{\blold{v_1} \shuffle \blold{v_2}}.
    \]
    By Lemma~\ref{lem:sup-phi}, $|\lbold_1^{\blold{v_1}}||\lbold_2^{\blold{v_2}}| \lesssim \varphi_{\varepsilon_1 \wedge \varepsilon_2}(\blold{v})$ for any $\blold{v} \in \blold{v_1} \shuffle \blold{v_2}$, and, by counting the shuffles, using Lemma~\ref{lem:count-shuffle-decompositions}, $\sum_{\blold{v_1}, \blold{v_2}} \brackets{\blold{v}}{\blold{v_1} \shuffle \blold{v_2}} = 2^{|\blold{v}|}$. Hence $|(\lbold_1 \shuffle \lbold_2)^{\blold{v}}| \lesssim 2^{|\blold{v}|} \varphi_{\varepsilon_1 \wedge \varepsilon_2}(\blold{v})$. Since factorial growth dominates exponentials, for any $\varepsilon^* < \varepsilon_1 \wedge \varepsilon_2$, $2^{|\blold{v}|} \lesssim (\norm{\blold{v}}!)^{\varepsilon_1 \wedge \varepsilon_2 - \varepsilon^*}$, so 
    \[
        |(\lbold_1 \shuffle \lbold_2)^{\blold{v}}| \lesssim \varphi_{\varepsilon^*}(\blold{v}), \quad \text{i.e. } \lbold_1 \shuffle \lbold_2 \in \mathcal{L}^1.
    \]
    For the projection, $|(\lbold|_{\blold{u}})^{\blold{v}}| = |\lbold^{\blold{vu}}| \lesssim \varphi_{\varepsilon}(\blold{vu})$. Using the inequality
    \[
        \norm{\blold{vu}}! = (\norm{\blold{v}} + \norm{\blold{u}})! \le \norm{\blold{v}}! \left[ (\norm{\blold{v}} + \norm{\blold{u}}) \right]^{\norm{\blold{u}}},
    \]
    we obtain $\varphi_{\varepsilon}(\blold{vu}) \le \varphi_{\varepsilon}(\blold{v}) \left[ (\norm{\blold{v}} + \norm{\blold{u}}) \right]^{\left( \frac{1}{2} - \varepsilon \right)\norm{\blold{u}}}$. Again the factorials dominate polynomials, so $|(\lbold|_{\blold{u}})^{\blold{v}}| \lesssim \varphi_{\varepsilon^*}(\blold{v})$ for some $\varepsilon^* < \varepsilon$, which yields $\lbold|_{\blold{u}} \in \mathcal{L}^1$.    
\end{proof}

\subsection{Complement for Section~\ref{sec:strong-discretisation}} \label{app:Lp-element-discretisation}

This section gathers the two auxiliary lemmas used in the proof of Proposition~\ref{prop:Lp-element-discretisation} and the proof of Lemma~\ref{lem:levy-area-words}.

\begin{lemma}[Products of moment bounds]\label{lem:prod-psi}
    Let $p \ge 1$, $m \in \N^*$ and $\blold{v_1}, \dots, \blold{v_m} \in V$ with $\blold{v_1} \cdots \blold{v_m} = \blold{v}$. Then
    \begin{equation}\label{eq:prod-phi-bound}
        \prod_{i=1}^m\psi_{p, t_i - s_i}(\blold{v_i})\le m^{\norm{\blold{v}}/2}\Big(1+\frac{\norm{\blold{v}}}2\Big)^{m/2}\psi_{p,T}(\blold{v}), \qquad 0 \le s_i \le t_i \le T.
    \end{equation}
\end{lemma}

\begin{proof}
    First, set $a_i\coloneqq\lfloor\norm{\blold{v_i}}/2\rfloor$, $A\coloneqq\sum_i a_i$, $B\coloneqq\lfloor\norm{\blold{v}}/2\rfloor$. By additivity of $\norm{\, \cdot \,}$, $(\norm{\blold{v}}-m)/2\le A\le B$, so $B-A\le\lfloor m/2\rfloor$. The multinomial bound $A!/\prod_i a_i!\le m^A$ and $B!/A!\le\max(1,B)^{B-A}$ give
    \[
        \prod_{i=1}^m\frac1{a_i!}\le\frac{m^A}{A!}\le\frac{m^{\norm{\blold{v}}/2}\,(1+\norm{\blold{v}}/2)^{m/2}}{B!}.
    \]
    The non-factorial factors of $\psi_{p,T}$ are multiplicative and $\psi_{p,t-s} \le \psi_{p,T}$, which, combined with the previous display, gives \eqref{eq:prod-phi-bound}.
\end{proof}

The proof of Lemma~\ref{lem:levy-area-words} relies on Rosenthal’s inequality~\citep*{Rosenthal1970}, which we recall in the form stated in the introduction of \citet*{ChenSung2020}. Let $n \in \N^*$ and $p \ge 2$. For a sequence $\{Z_i : 1 \le i \le n\}$ of independent, centred random variables with $\E[|Z_i|^p] < \infty$, there exists a positive constant $C_p$ such that
\begin{equation} \label{th:Rosenthal-inequality}
    \E\left[ \Big| \sum_{i=1}^n Z_i \Big|^p \right] \le C_p \max \left\{ \sum_{i=1}^n \E\big[|Z_i|^p\big] \ ; \left( \sum_{i=1}^n \E[Z_i^2] \right)^{p/2}\right\}.
\end{equation}
For $p=2$, \eqref{th:Rosenthal-inequality} holds with $C_2=1$, since the variance of a sum of independent centred variables is the sum of their variances.

\begin{lemma}[A combinatorial sum]\label{lem:combinatorial-sum}
    Let $p \in 2\N^*$. There exists $C_p > 0$ such that, for every $N \in \N^*$ and $n \in \intset{1}{N}$,
    \[
        S_{n,N} \coloneqq \sum_{\substack{\alpha \in \N^n \\ |\alpha| = p}} \binom{p}{\alpha_1, \dots, \alpha_n} N^{- \frac{1}{2} \#\{ l \in \intset{1}{n} \ : \ \alpha_l \in 2\N + 1\} } \le C_p\, N^{p/2}.
    \]
\end{lemma}

\begin{proof}
    Introducing
    \[
        a_k = \begin{cases}
            1, &\text{if $k$ is even}\\
            N^{-1/2}, &\text{if $k$ is odd}
        \end{cases}, \quad \text{we can rewrite } S_{n,N} = p! \sum_{\substack{\alpha \in \N^n \\ |\alpha| = p}} \prod_{k = 1}^{n} \frac{ a_{\alpha_k} }{ \alpha_k! }.
    \]
    Let $f$ denote the power series with coefficients $a$. Then, $f(x) \coloneqq \sum_k a_k \frac{x^k}{k!} = \cosh(x) + N^{-1/2}\sinh(x) = p_N e^x + q_N e^{-x}$ with $p_N \coloneqq \frac{1}{2}\big(1 + N^{-1/2}\big)$ and $q_N \coloneqq 1-p_N$. We can relate $f$ to $S_{n,N}$ through the identity $F(x) \coloneqq \left[ f(x) \right]^n = \sum_{k=0}^\infty x^k \sum_{|\alpha| = k} \prod_{i = 1}^{n} \frac{ a_{\alpha_i} }{ \alpha_i! }$. The quantity $S_{n,N}$ is $p!$ times the coefficient of order $p$ in the power series $F$, recovered via $F^{(p)}(0) = p! \sum_{|\alpha| = p} \prod_{i = 1}^{n} \frac{ a_{\alpha_i} }{ \alpha_i! } = S_{n,N}$. Using the binomial theorem, $F(x) = (p_N e^x + q_N e^{-x})^n = \sum_{k = 0}^n \binom{n}{k} p_N^k q_N^{n-k} e^{(2k - n)x}$, and noting that $\frac{d^p}{dx^p}e^{(2k-n)x} = (2k-n)^p e^{(2k-n)x}$, we deduce $F^{(p)}(0) = \sum_{k=0}^n\binom{n}{k} p_N^k q_N^{n-k}(2k-n)^p = \E\left[ (2Z_N - n)^p \right]$ where $Z_N$ is a binomial random variable of parameters $n$ and $p_N \in [0, 1]$. Next, using the identity $1 = 2p_N - N^{-1/2}$ together with the convexity inequality $|x + y|^p \le 2^{p-1}(|x|^p + |y|^p)$, we deduce that
    \[
        \E\left[ (2Z_N - n)^p \right] = \E\left[ \left\{2(Z_N - p_N n) +N^{-\frac{1}{2}}n\right\}^p \right] \lesssim \E\left[|Z_N - p_N n|^p \right] + N^{-\frac{p}{2}}n^p.
    \]
    Since $Z_N - p_N n \overset{d}{=} \sum_{i = 1}^n x_i$, where $x_i = \xi_i - p_N$, and $(\xi_i)$ are i.i.d. Bernoulli random variables with parameter $p_N$ such that $Z_N \overset{d}{=} \sum_{i=1}^n \xi_i$, we can apply Rosenthal's inequality~\eqref{th:Rosenthal-inequality} and obtain, for some constant $C_p > 0$, $\E\big[|Z_N - p_N n|^p \big] \le C_p \max \big\{ \sum_{i=1}^n \E[|x_i|^p]; \big(\sum_{i=1}^n \E[ x_i^2]\big)^{p/2} \big\} \lesssim n^{p / 2}$ as $\E[|x_i|^p] \le 1$ and $\big(\sum_{i=1}^n \E[ x_i^2]\big)^{p/2} \le n^{p/2}$. Combining these estimates, $\E\left[ (2Z_N - n)^p \right] \lesssim n^{p/2} + N^{-p/2}n^p \lesssim N^{p/2}$. Hence $S_{n,N} \lesssim N^{p/2}$, which concludes.
\end{proof}

We can now prove Lemma~\ref{lem:levy-area-words}. The local errors on the Lévy-area words are the Lévy areas of the successive steps; we expand the $p$-th moment and show that only terms with even multiplicities contribute at leading order.

\begin{proof}[Proof of Lemma~\ref{lem:levy-area-words}]
    We exploit the symmetry of the Lévy area, defined by
    \[
        A_{\blold{ij}}(s, t) \coloneqq \frac{1}{2}\left(\int_s^t W^i_{s, \tau} \diff W^j_\tau - \int_s^t W^j_{s, \tau} \diff W^i_\tau\right),
    \]
    where $W_{s, t} = W_t - W_s$. One has $\WH_{t_{k-1}, t_k}^{\blold{ij}} - \WT_{t_{k-1}, t_k}^{\blold{ij}} = A_{\blold{ij}}(t_{k-1}, t_k)$. Moreover, for $q \in \N$, 
    \[
        \E\left[ A_{\blold{ij}}(t_{k-1}, t_k)^{2q} \right] \lesssim h^{2q} \lesssim N^{-2q} \quad \text{and} \quad \E\left[ A_{\blold{ij}}(t_{k-1}, t_k)^{2q + 1} \right] = 0,
    \]
    the first bound follows from Proposition~\ref{prop:signature-moment-bound}. Using the triangle inequality and the fact that $\tilde{t}_k = t_k$ whenever $k \le n$, we decompose
    \begin{align*}
        \textnorm{\Delta^{\blold{v_1}, \blold{ij}, \blold{v_3}}(t)}_p 
        &\le \norm{\sum_{k = 1}^{n}\WH_{t_{k-1}}^{\blold{v_1}} \left(\WH_{t_{k-1}, t_k}^{\blold{ij}} - \WT_{t_{k-1}, t_k}^{\blold{ij}}\right) \WT_{t_{k}, t}^{\blold{v_3}}}_p\\
        & \quad + \norm{\WH_{\tilde{t}_{n}}^{\blold{v_1}} \left(\WH_{\tilde{t}_{n}, \tilde{t}_{n+1}}^{\blold{ij}} - \WT_{\tilde{t}_{n}, \tilde{t}_{n+1}}^{\blold{ij}}\right) \WT_{\tilde{t}_{n+1}, t}^{\blold{v_3}}}_p.
    \end{align*}
    The second term is a single local error. By independence and Proposition~\ref{prop:signature-moment-bound}, as in the proof of Lemma~\ref{lem:local-errors}, and since $\norm{\blold{ij}}=2$ and $\tilde{t}_{n+1} - \tilde{t}_n \le h$, its $L^p$ norm is $\lesssim N^{-1/2}$. We now focus on the first term. Assume that $n \ge 1$ (the case $n = 0$ being immediate). Expanding the norm yields
    \begin{equation} \label{eq:decomposition_2}
        \norm{\sum_{k = 1}^{n}\WH_{t_{k-1}}^{\blold{v_1}} \left(\WH_{t_{k-1}, t_k}^{\blold{ij}} - \WT_{t_{k-1}, t_k}^{\blold{ij}}\right) \WT_{t_{k}, t}^{\blold{v_3}}}_p^p = \sum_{1 \le n_1, \dots, n_p \le n}\E\left[ \Delta^{\blold{v_1}, \blold{ij}, \blold{v_3}}_{n_1, \dots, n_p}(t) \right],
    \end{equation}
    where $\Delta^{\blold{v_1}, \blold{ij}, \blold{v_3}}_{n_1, \dots, n_p}(t) \coloneqq \prod_{i=1}^p \WH_{t_{n_i-1}}^{\blold{v_1}} A_{\blold{ij}}(t_{n_i-1}, t_{n_i}) \WT_{t_{n_i}, t}^{\blold{v_3}}$. Studying the sum of $\Delta^{\blold{v_1}, \blold{ij}, \blold{v_3}}_{n_1, \dots, n_p}(t)$, we show below that
    \[
        \norm{\sum_{k = 1}^{n}\WH_{t_{k-1}}^{\blold{v_1}} \left(\WH_{t_{k-1}, t_k}^{\blold{ij}} - \WT_{t_{k-1}, t_k}^{\blold{ij}}\right) \WT_{t_{k}, t}^{\blold{v_3}}}_p^p \lesssim N^{-p/2}.
    \]
    We first treat the case $p = 2$, which contains the main idea, and then the general case.

    \medskip

    \noindent\emph{The case $p = 2$.}
    The terms $\Delta^{\blold{v_1}, \blold{ij}, \blold{v_3}}_{n_1, n_2}(t)$ can be divided into diagonal terms ($n_1 = n_2$) and off-diagonal terms ($n_1 \ne n_2$). Using the independence of the Brownian increments together with Proposition~\ref{prop:signature-moment-bound}, as in the proof of Lemma~\ref{lem:local-errors}, we can bound the $n$ diagonal terms
    \[
        \E\left[ \Delta^{\blold{v_1}, \blold{ij}, \blold{v_3}}_{k, k} \right] = \E\left[ \left( \WH_{t_{k-1}}^{\blold{v_1}} A_{\blold{ij}}(t_{k-1}, t_{k}) \WT_{t_{k}, t}^{\blold{v_3}} \right)^2 \right] \lesssim \E\left[ \big(A_{\blold{ij}}(t_{k-1}, t_{k})\big)^2 \right] \lesssim N^{-2}.
    \]
    For the off-diagonal terms, we first decompose the signatures in order to isolate contributions coming from independent Brownian increments. This is the key step of the general proof, illustrated here in the case $p = 2$. Let $1 \le k < l \le n$. We consider the $n(n-1)$ terms of the form
    \[
        \E\left[ \Delta^{\blold{v_1}, \blold{ij}, \blold{v_3}}_{k, l} \right] = \E\left[ \left( \WH_{t_{k-1}}^{\blold{v_1}} A_{\blold{ij}}(t_{k-1}, t_{k}) \WT_{t_{k}, t}^{\blold{v_3}} \right)\left(\WH_{t_{l-1}}^{\blold{v_1}} A_{\blold{ij}}(t_{l-1}, t_{l}) \WT_{t_{l}, t}^{\blold{v_3}} \right) \right].
    \]
    Since $\WT_{t_{k}, t}^{\blold{v_3}}$ contains contributions from the interval $[t_{l-1}, t_l]$, we extract this part to combine it with $A_{\blold{ij}}(t_{l-1}, t_{l})$. Similarly, $\WH_{t_{l-1}}^{\blold{v_1}}$ contains contributions from $[t_{k-1}, t_k]$, which we isolate to combine with $A_{\blold{ij}}(t_{k-1}, t_{k})$. To achieve this decomposition, we apply Chen's identity~\eqref{eq:chen-identity} and reorganise the terms as illustrated in Figure~\ref{fig:decomposition}.
    
    \begin{figure}[ht]
        \centering
        \begin{tikzpicture}
        
        \definecolor{colorA}{RGB}{225,225,225}
        \definecolor{colorB}{RGB}{180,180,180}
        
        \foreach \i/\col/\val in {
          0/colorA/$\WH_{t_{k-1}}^{\blold{v_1}}$,
          1/colorB/$A_{\blold{ij}}(t_{k-1}{,}\ t_{k})$,
          2/colorA/$\WT_{t_{k}, t_{l-1}}^{\blold{v_{31}}}$,
          3/colorA/$\WT_{t_{l-1}, t_{l}}^{\blold{v_{32}}}$,
          4/colorA/$\WT_{t_{l}, t}^{\blold{v_{33}}}$
        } {
          \draw[fill=\col] (3*\i,0) rectangle ++(3,1);
          \node at (3*\i+1.5,0.5) {\val};
        }
        
        \draw[decorate,decoration={brace,amplitude=5pt}]
          (6,1.2) -- (15,1.2)
          node[midway,above=6pt] {$\WT_{t_{k}, t}^{\blold{v_3}} =\sum_{\blold{v_{31}v_{32}v_{33}} = \blold{v_3}} \WT_{t_{k}, t_{l-1}}^{\blold{v_{31}}} \WT_{t_{l-1}, t_{l}}^{\blold{v_{32}}} \WT_{t_{l}, t}^{\blold{v_{33}}} $};
        
        \foreach \i/\col/\val in {
          0/colorA/$\WH_{t_{k-1}}^{\blold{v_{11}}}$,
          1/colorA/$\WH_{t_{k-1}, t_{k}}^{\blold{v_{12}}}$,
          2/colorA/$\WH_{t_{k}, t_{l-1}}^{\blold{v_{13}}}$,
          3/colorB/$A_{\blold{ij}}(t_{l-1}{,}\ t_{l})$,
          4/colorA/$\WT_{t_{l}, t}^{\blold{v_3}}$
        } {
          \draw[fill=\col] (3*\i,-1.5) rectangle ++(3,1);
          \node at (3*\i+1.5,-1) {\val};
        }
        
        \draw[decorate,decoration={brace,amplitude=5pt}]
          (9,-1.7) -- (0,-1.7)
          node[midway,below=6pt] {$\WH_{t_{l-1}}^{\blold{v_1}} =\sum_{\blold{v_{11}v_{12}v_{13}} = \blold{v_1}} \WH_{t_{k-1}}^{\blold{v_{11}}} \WH_{t_{k-1}, t_{k}}^{\blold{v_{12}}} \WH_{t_{k}, t_{l-1}}^{\blold{v_{13}}}$};
          
        \end{tikzpicture}
        \caption{Decomposition of the signatures into independent increments}
        \label{fig:decomposition}
    \end{figure}
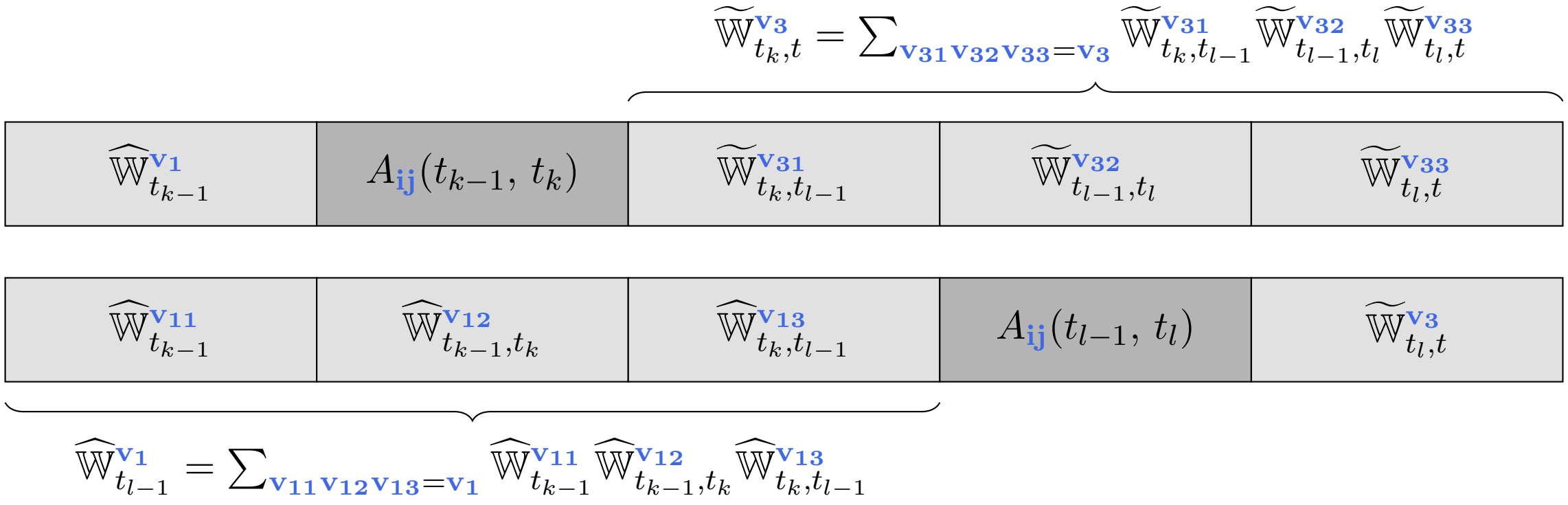
    
    \noindent Rearranging the resulting terms and using independence of increments over disjoint intervals, we obtain
    \begin{align*}
        &\E\left[ \left( \WH_{t_{k-1}}^{\blold{v_1}} A_{\blold{ij}}(t_{k-1}, t_{k}) \WT_{t_{k}, t}^{\blold{v_3}} \right)\left(\WH_{t_{l-1}}^{\blold{v_1}} A_{\blold{ij}}(t_{l-1}, t_{l}) \WT_{t_{l}, t}^{\blold{v_3}} \right) \right]\\[0.2cm]
        &= \sum_{\substack{\blold{v_{31}v_{32}v_{33}} = \blold{v_3} \\ \blold{v_{11}v_{12}v_{13}} = \blold{v_1}}} \E \left[ \left( \WH_{t_{k-1}}^{\blold{v_1}} A_{\blold{ij}}(t_{k-1}, t_{k}) \WT_{t_{k}, t_{l-1}}^{\blold{v_{31}}}\WT_{t_{l-1}, t_{l}}^{\blold{v_{32}}}\WT_{t_{l}, t}^{\blold{v_{33}}} \right)\right.\\[-0.5cm]
        & \hspace{4cm} \left.\left(\WH_{t_{k-1}}^{\blold{v_{11}}}\WH_{t_{k-1}, t_{k}}^{\blold{v_{12}}}\WH_{t_{k}, t_{l-1}}^{\blold{v_{13}}} A_{\blold{ij}}(t_{l-1}, t_{l}) \WT_{t_{l}, t}^{\blold{v_3}} \right) \right]\\[0.2cm]
        &= \sum_{\substack{\blold{v_{31}v_{32}v_{33}} = \blold{v_3} \\ \blold{v_{11}v_{12}v_{13}} = \blold{v_1}}} \E\left[ \WH_{t_{k-1}}^{\blold{v_1}} \WH_{t_{k-1}}^{\blold{v_{11}}} \right] \E\left[ A_{\blold{ij}}(t_{k-1}, t_{k}) \WH_{t_{k-1}, t_{k}}^{\blold{v_{12}}} \right] \E\left[ \WT_{t_{k}, t_{l-1}}^{\blold{v_{31}}} \WH_{t_{k}, t_{l-1}}^{\blold{v_{13}}} \right]\\[-0.5cm]
        & \hspace{4cm} \E\left[ \WT_{t_{l-1}, t_{l}}^{\blold{v_{32}}} A_{\blold{ij}}(t_{l-1}, t_{l}) \right] \E\left[ \WT_{t_{l}, t}^{\blold{v_{33}}} \WT_{t_{l}, t}^{\blold{v_3}} \right].
    \end{align*}
    Using Proposition~\ref{prop:signature-moment-bound}, we deduce
    \begin{align*}
        &\E\left[ \left( \WH_{t_{k-1}}^{\blold{v_1}} A_{\blold{ij}}(t_{k-1}, t_{k}) \WT_{t_{k}, t}^{\blold{v_3}} \right)\left(\WH_{t_{l-1}}^{\blold{v_1}} A_{\blold{ij}}(t_{l-1}, t_{l}) \WT_{t_{l}, t}^{\blold{v_3}} \right) \right]\\[0.2cm]
        &\qquad \qquad \lesssim \sum_{\substack{\blold{v_{31}v_{32}v_{33}} = \blold{v_3} \\ \blold{v_{11}v_{12}v_{13}} = \blold{v_1}}} \absE{A_{\blold{ij}}(t_{k-1}, t_{k}) \WH_{t_{k-1}, t_{k}}^{\blold{v_{12}}}} \absE{\WT_{t_{l-1}, t_{l}}^{\blold{v_{32}}} A_{\blold{ij}}(t_{l-1}, t_{l})}.
    \end{align*}
    Then,
    \begin{itemize}
        \item If $\blold{v_{12}} = \empt$, then $\E\left[ A_{\blold{ij}}(t_{k-1}, t_{k}) \WH_{t_{k-1}, t_{k}}^{\blold{v_{12}}} \right] = \E\left[ A_{\blold{ij}}(t_{k-1}, t_{k})\right] = 0$.
        \item Otherwise, $\norm{\blold{v_{12}}} \ge 1$ and, by Cauchy--Schwarz,
        \[
            \absE{A_{\blold{ij}}(t_{k-1}, t_{k}) \WH_{t_{k-1}, t_{k}}^{\blold{v_{12}}}} \le \norm{A_{\blold{ij}}(t_{k-1}, t_{k})}_2 \norm{\WH_{t_{k-1}, t_{k}}^{\blold{v_{12}}}}_2 
            \lesssim h^{1 + \norm{\blold{v_{12}}}/2} \lesssim N^{-1-\norm{\blold{v_{12}}}/2}.
        \]
    \end{itemize}
    Thus, combining both cases, we have $\textabsE{A_{\blold{ij}}(t_{k-1}, t_{k}) \WH_{t_{k-1}, t_{k}}^{\blold{v_{12}}}} \lesssim N^{-3/2}$. $\textabsE{\WT_{t_{l-1}, t_{l}}^{\blold{v_{32}}} A_{\blold{ij}}(t_{l-1}, t_{l})} \lesssim N^{-3/2}$ is obtained similarly, which yields $\E[ \Delta^{\blold{v_1}, \blold{ij}, \blold{v_3}}_{k, l} ] \lesssim N^{-3}$. Summing the $n$ diagonal estimates with the $n(n-1)$ off-diagonal ones, we finally obtain
    \[
        \norm{\sum_{k = 1}^{n}\WH_{t_{k-1}}^{\blold{v_1}} \left(\WH_{t_{k-1}, t_k}^{\blold{ij}} - \WT_{t_{k-1}, t_k}^{\blold{ij}}\right) \WT_{t_{k}, t}^{\blold{v_3}}}_2^2 \lesssim nN^{-2} + n(n-1)N^{-3} \lesssim N^{-1}.
    \]

    \medskip

    \noindent\emph{General case.}
    We now bound the right-hand side of \eqref{eq:decomposition_2} for $p > 2$, using a similar approach. By the commutativity of the product, for any permutation $\sigma : \{1, \dots, p\} \to \{1, \dots, p\}$, we have $\Delta^{\blold{v_1}, \blold{ij}, \blold{v_3}}_{n_1, \dots, n_p}(t) = \Delta^{\blold{v_1}, \blold{ij}, \blold{v_3}}_{n_{\sigma(1)}, \dots, n_{\sigma(p)}}(t)$. Then, we restrict our attention to the case $n_1 \le \cdots \le n_p$. Let $q \le p$, and let $m_1, \dots, m_q$ and $\alpha_1, \dots, \alpha_q$ in $\N$ be such that $m_1 < \cdots < m_q$, $|\alpha| = p$ and $n_k = m_i$ for $\sum_{j < i} \alpha_j < k \le \sum_{j \le i} \alpha_j$. Define $m_0 = 0$, $m_{q+1} = n+1$ and $t_{-1} = 0$, while redefining $t_{n+1} = t$. Using Chen's identity \eqref{eq:chen-identity}, we obtain
    \begin{align*}
        \WH_{0, t_{n_k - 1}} &= \bigotimes_{l = 1}^{\kt} \WH_{t_{m_{l-1}-1}, t_{m_l - 1}} = \bigotimes_{l = 1}^{\kt} \left( \WH_{t_{m_{l-1} - 1}, t_{m_{l-1}}} \otimes \WH_{t_{m_{l-1}}, t_{m_l - 1}} \right)\\
        \WT_{t_{n_k}, t} &= \bigotimes_{l = \kt}^{q} \left( \WT_{t_{m_l}, t_{m_{l+1} - 1}} \otimes \WT_{t_{m_{l+1} - 1}, t_{m_{l+1}}} \right)
    \end{align*}
    where $\kt$ is defined by the condition $n_k = m_{\kt}$. This yields the same type of decomposition as in Figure~\ref{fig:decomposition}, now in the general setting. As the next steps are mostly computational, we introduce simplified notation $\blold{u} \coloneqq \blold{v_1}$, $\blold{\upsilon} \coloneqq \blold{v_3}$ and $A_k \coloneqq A_{\blold{ij}}(t_{m_k - 1}, t_{m_k})$. Expanding the signatures, we obtain
    \[
        \E\left[ \Delta^{\blold{u}, \blold{ij}, \blold{\upsilon}}_{n_1, \dots, n_p}(t) \right] 
        = \sum_{\substack{\bigotimes_{l=1}^{\kt} \udown_{\blue{l}}^{\blue{k}} \uup_{\blue{l}}^{\blue{k}} = \blold{u} \\ \bigotimes_{l=\kt}^{q} \updown_{\blue{l}}^{\blue{k}} \upup_{\blue{l}}^{\blue{k}} = \blold{\upsilon}}, \ 1 \le k \le p} \E\left[ \Delta\left( \udown, \uup, \updown, \upup \right) \right]
    \]
    where the sum runs over every $( \udown, \uup, \updown, \upup ) = \{ (\udown_{\blue{l}}^{\blue{k}})_{l=1}^{\kt}, (\uup_{\blue{l}}^{\blue{k}})_{l=1}^{\kt}, (\updown_{\blue{l}}^{\blue{k}})_{l=\kt}^{q}, (\upup_{\blue{l}}^{\blue{k}})_{l=\kt}^{q} : 1 \le k \le p\}$ such that $\bigotimes_{l=1}^{\kt} \udown_{\blue{l}}^{\blue{k}} \uup_{\blue{l}}^{\blue{k}} = \blold{u}$ and $\bigotimes_{l=\kt}^{q} \updown_{\blue{l}}^{\blue{k}} \upup_{\blue{l}}^{\blue{k}} = \blold{\upsilon}$, and where
    \begin{align*}
        \Delta\left( \udown, \uup, \updown, \upup \right) 
        &\coloneqq \prod_{k=1}^p \bigg\{ \prod_{l = 1}^{\kt} \left( \WH^{\udown_{\blue{l}}^{\blue{k}}}_{t_{m_{l-1} - 1}, t_{m_{l-1}}} \WH^{\uup_{\blue{l}}^{\blue{k}}}_{t_{m_{l-1}}, t_{m_l - 1}} \right) \bigg\}\\
        & \qquad \qquad A_{\kt} \bigg\{ \prod_{l = \kt}^q \left( \WT^{\updown_{\blue{l}}^{\blue{k}}}_{t_{m_l}, t_{m_{l+1} - 1}} \WT^{\upup_{\blue{l}}^{\blue{k}}}_{t_{m_{l+1} - 1}, t_{m_{l+1}}} \right) \bigg\}.
    \end{align*}
    Using independence of the signatures over disjoint time intervals, we can extract the terms $\WH^{\uup_{\blue{l}}^{\blue{k}}}_{t_{m_{l-1}}, t_{m_l - 1}}$ and $\WT^{\updown_{\blue{l}}^{\blue{k}}}_{t_{m_l}, t_{m_{l+1} - 1}}$, and applying Proposition~\ref{prop:signature-moment-bound}, we obtain
    \begin{align*}
        \left|\E\left[\Delta\left( \udown, \uup, \updown, \upup \right)\right]\right| \lesssim \left|\E\left[ \prod_{k=1}^p \bigg\{ \prod_{l = 1}^{\kt} \WH^{\udown_{\blue{l}}^{\blue{k}}}_{t_{m_{l-1} - 1}, t_{m_{l-1}}} \bigg\} A_{\kt} \bigg\{ \prod_{l = \kt}^q \WT^{\upup_{\blue{l}}^{\blue{k}}}_{t_{m_{l+1} - 1}, t_{m_{l+1}}} \bigg\} \right]\right|.
    \end{align*}
    Grouping the indices according to the sets $\{k : \kt = \kappa\}$, for $1 \le \kappa \le q$, and rearranging the products,
    \[
        \prod_{k=1}^p \prod_{l = 1}^{\kt} \WH^{\udown_{\blue{l}}^{\blue{k}}}_{t_{m_{l-1} - 1}, t_{m_{l-1}}} = \prod_{\kappa=1}^q \prod_{l = 1}^{\kappa} \prod_{\kt = \kappa} \WH^{\udown_{\blue{l}}^{\blue{k}}}_{t_{m_{l-1} - 1}, t_{m_{l-1}}} = \prod_{l=1}^q \prod_{\kappa = l}^{q} \prod_{\kt = \kappa} \WH^{\udown_{\blue{l}}^{\blue{k}}}_{t_{m_{l-1} - 1}, t_{m_{l-1}}}.
    \]
    Define the following groupings of terms associated with each time interval:
    \[
        \pih_l \coloneqq \prod_{\kappa = l}^{q} \prod_{\kt = \kappa} \WH^{\udown_{\blue{l}}^{\blue{k}}}_{t_{m_{l-1} - 1}, t_{m_{l-1}}} \quad \text{and} \quad \pit_l \coloneqq \prod_{\kappa = 1}^{l} \prod_{\kt = \kappa} \WT^{\upup_{\blue{l}}^{\blue{k}}}_{t_{m_{l+1} - 1}, t_{m_{l+1}}}
    \]
    for $l \in \intset{1}{q}$ such that
    \[
        \E\left[ \prod_{k=1}^p \bigg\{ \prod_{l = 1}^{\kt} \WH^{\udown_{\blue{l}}^{\blue{k}}}_{t_{m_{l-1} - 1}, t_{m_{l-1}}} \bigg\} A_{\kt} \bigg\{ \prod_{l = \kt}^q \WT^{\upup_{\blue{l}}^{\blue{k}}}_{t_{m_{l+1} - 1}, t_{m_{l+1}}} \bigg\} \right] = \E\left[ \prod_{l=1}^q \pih_l A_l^{\alpha_l} \pit_l \right]
    \]
    since we have $\prod_{k=1}^p A_{\kt} = \prod_{l=1}^q A_l^{\alpha_l}$. Hence, using independence of the different factors,
    \[
        \E\left[ \prod_{l=1}^q \pih_l A_l^{\alpha_l} \pit_l \right] = \E\left[ \pih_1 \right] \left(\prod_{l=1}^q \E\left[  \pih_{l+1} A_l^{\alpha_l} \pit_{l-1} \right]\right) \E\left[ \pit_q \right]
    \]
    where $\pih_{q+1} = \pit_{0} = 1$. Using Proposition~\ref{prop:signature-moment-bound}, we deduce $\absE{\Delta\left( \udown, \uup, \updown, \upup \right)} \lesssim \textstyle\prod_{l=1}^q \left|\E\left[  \pih_{l+1} A_l^{\alpha_l} \pit_{l-1} \right]\right|$. Define $\blold{w}_{\blue{l}} = \big(\bigotimes_{\kappa = l+1}^q \bigotimes_{\kt = \kappa} \udown_{\blue{l+1}}^{\blue{k}}\big)\big(\bigotimes_{\kappa = 1}^{l-1} \bigotimes_{\kt = \kappa} \upup_{\blue{l-1}}^{\blue{k}}\big)$. Tracking the coordinates of the signatures in $\pih_{l+1}$ and $\pit_{l-1}$ and applying again Proposition~\ref{prop:signature-moment-bound}, one can show that
    \[
        \absE{ \pih_{l+1} A_l^{\alpha_l} \pit_{l-1} } \lesssim h^{\alpha_l + \frac{1}{2}\norm{\blold{w}_{\blue{l}}}} \lesssim N^{-\alpha_l - \frac{1}{2}\norm{\blold{w}_{\blue{l}}}}.
    \]
    Moreover, since the Lévy area is symmetric, if $\norm{\blold{w}_{\blue{l}}} = 0$, that is, $\blold{w}_{\blue{l}} = \empt$, meaning that $\pih_{l+1} = \pit_{l-1} = 1$, then for odd $\alpha_l$, we have $\textabsE{\pih_{l+1} A_l^{\alpha_l} \pit_{l-1}} = 0$, which implies $\absE{ \pih_{l+1} A_l^{\alpha_l} \pit_{l-1} } \lesssim N^{-\alpha_l - \frac{1}{2}\norm{\blold{w}_{\blue{l}}}}\indic{\norm{\blold{w}_{\blue{l}}} > 0}$ $\le N^{-\alpha_l - \frac{1}{2}}$. When $\alpha_l$ is even, it suffices to use the bound $\textabsE{ \pih_{l+1} A_l^{\alpha_l} \pit_{l-1} } \lesssim N^{-\alpha_l}$. Thus, defining
    \[
        I(n_1, \dots, n_p) \coloneqq \big\{ l \in \{1, \dots, q \} \ : \ \alpha_l \in 2\N + 1 \big\}
    \]
    and using $\sum_{l=1}^q \alpha_l = p$, we obtain $\textabsE{\Delta\left( \udown, \uup, \updown, \upup \right)} \lesssim N^{-p - \frac{1}{2}\# I(n_1, \dots, n_p)}$, and more generally
    \[
        \textabsE{ \Delta^{\blold{u}, \blold{ij}, \blold{\upsilon}}_{n_1, \dots, n_p}(t) } \lesssim N^{-p - \frac{1}{2}\# I(n_1, \dots, n_p)}.
    \]
    Writing $S \coloneqq \sum_{1 \le n_1, \dots, n_p \le n}  N^{-\frac{1}{2}\# I(n_1, \dots, n_p)}$, since the values of $n_1, \dots, n_p$ only intervene through the multiplicities of repeated indices, $S = S_{n,N}$, and Lemma~\ref{lem:combinatorial-sum} gives $S \lesssim N^{p/2}$. Hence, \eqref{eq:decomposition_2} gives
    
    \[
        \norm{\sum_{k = 1}^{n}\WH_{t_{k-1}}^{\blold{v_1}} \left(\WH_{t_{k-1}, t_k}^{\blold{ij}} - \WT_{t_{k-1}, t_k}^{\blold{ij}}\right) \WT_{t_{k}, t}^{\blold{v_3}}}_p^p \lesssim N^{-p} \ S \lesssim N^{-\frac{p}{2}}.
    \]

    We finally track the constant: we improve our estimate of $\Delta^{\blold{u}, \blold{ij}, \blold{\upsilon}}_{n_1, \dots, n_p}(t)$, by collecting the fragments of the constant hidden through the general case and using Hölder inequality. We obtain
    \begin{itemize}
        \item $\textabsE{\prod_{k=1}^p \big\{ \prod_{l = 1}^{\kt} \WH^{\uup_{\blue{l}}^{\blue{k}}}_{t_{m_{l-1}}, t_{m_l - 1}} \big\}\big\{ \prod_{l = \kt}^q \WT^{\updown_{\blue{l}}^{\blue{k}}}_{t_{m_l}, t_{m_{l+1} - 1}} \big\}} \le \prod_{k=1}^p\big\{\prod_{l=1}^{\kt} \psi(\uup_{\blue{l}}^{\blue{k}})\big\}\big\{\prod_{l=\kt}^{q} \psi(\updown_{\blue{l}}^{\blue{k}})\big\}$,
        \item $\textabsE{\pih_1} \le \prod_{k=1}^p\psi(\udown_{\blue{1}}^{\blue{k}})$ and $\textabsE{\pit_q} \le \prod_{k=1}^p\psi(\upup_{\blue{q}}^{\blue{k}})$,
        \item $\absE{ \pih_{l+1} A_l^{\alpha_l} \pit_{l-1} } \le \big(2 \psi(\blold{ij}) \big)^{\alpha_l} \big\{ \prod_{k : \kt \ge l+1} \psi(\udown_{\blue{l+1}}^{\blue{k}}) \big\} \big\{ \prod_{k : \kt \le l-1} \psi(\upup_{\blue{l-1}}^{\blue{k}}) \big\} N^{-\alpha_l - \frac{1}{2}\indic{\alpha_l \notin 2\N^*}}$.
    \end{itemize}
    The last inequality follows from Hölder's inequality applied to the $p$ factors supported on $[t_{m_l-1}, t_{m_l}]$, from $\textnorm{A_l}_p \le 2\psi(\blold{ij})N^{-1}$ and $\psi_{p,h}(\blold{u}) \le \psi(\blold{u})N^{-\norm{\blold{u}}/2}$ (Proposition~\ref{prop:signature-moment-bound}), and from the parity argument above. Combining all the fragments,
    \[
        \left|\E\left[\Delta\left( \udown, \uup, \updown, \upup \right)\right]\right| \le 2^p N^{-p - \frac{1}{2}\# I(n_1, \dots, n_p)} \textstyle{\prod_{k=1}^p} \Psi_k,
    \]
    with $\Psi_k \coloneqq \big\{ \prod_{l=1}^{\kt} \psi(\udown_{\blue{l}}^{\blue{k}})\psi(\uup_{\blue{l}}^{\blue{k}}) \big\} \psi(\blold{ij}) \big\{ \prod_{l=\kt}^{q} \psi(\updown_{\blue{l}}^{\blue{k}})\psi(\upup_{\blue{l}}^{\blue{k}}) \big\}$. These $2q + 3 \le 2p + 3$ subwords concatenate to $\blold{v}$, so Lemma~\ref{lem:prod-psi} gives $\Psi_k \le (2p + 3)^{\norm{\blold{v}}/2} \big(1 + \norm{\blold{v}}/2\big)^{p + 3/2} \psi(\blold{v})$. Summing over the at most $(|\blold{v}|+1)^{2p^2}$ decompositions of $\blold{u}$ and $\blold{\upsilon}$, and over $(n_1,\dots,n_p)$ (Lemma~\ref{lem:combinatorial-sum}), then taking the $p$-th root,
    \[
        \textnorm{\Delta^{\blold{v_1},\blold{ij},\blold{v_3}}(t)}_p\le P(\norm{\blold{v}})(2p+3)^{\norm{\blold{v}}/2}\psi(\blold{v})N^{-1/2}
    \]
    for a polynomial $P$, which concludes.
\end{proof}

\subsection{Complement for the Proof of Lemma~\ref{lem:cdf-bound}}

\begin{lemma} \label{lem:cdf-comparison}
    Let $Z_1, Z_2$ be real random variables with distribution functions $F_1, F_2$. Then, for all $x \in \R$ and $\eta > 0$, $\big| F_1(x) - F_2(x) \big| \le F_1(x + \eta) - F_1(x - \eta) + \Prob\big( |Z_1 - Z_2| > \eta \big)$.
\end{lemma}

\begin{proof}
    Let $E \coloneqq \{ |Z_1 - Z_2| \le \eta \}$. On $E$, $\{Z_1 \le x - \eta\} \subset \{Z_2 \le x\}$, so that $F_1(x - \eta) \le F_2(x) + \Prob(E^c)$, that is, $F_1(x) - F_2(x) \le F_1(x) - F_1(x - \eta) + \Prob(E^c)$. Similarly, $F_2(x) - F_1(x) \le F_1(x + \eta) - F_1(x) + \Prob(E^c)$. Since $F_1$ is non-decreasing, both right-hand sides are bounded by $F_1(x + \eta) - F_1(x - \eta) + \Prob(E^c)$.
\end{proof}


\renewcommand{\doi}[1]{\url{https://doi.org/#1}}
\bibliography{_references}

\end{document}